\documentclass[a4paper,twoside,11pt,english,fleqn]{article}

\usepackage[all]{xy}
\usepackage[intlimits]{amsmath}
\usepackage{amssymb,amsthm,graphicx}
\usepackage[dvips,dvipdf,pdftex]{epsfig}
\usepackage{babel}
\usepackage[T1]{fontenc}
\usepackage{fancyhdr,stmaryrd,color,times,euler,euscript,eufrak,array}

\CompileMatrices

  \newtheorem{thm*}{Theorem}[subsection]
  \newtheorem{cor*}[thm*]{Corollary}
  \newtheorem{lem*}[thm*]{Lemma}
  \newtheorem{prop*}[thm*]{Proposition}
\theoremstyle{definition}
  \newtheorem{defn*}[thm*]{Definition}
  \newtheorem{ex*}[thm*]{Example}
  \newtheorem{rem*}[thm*]{Remark}
  \newtheorem{ntt*}[thm*]{Notation}
\theoremstyle{remark}
  \newtheorem*{dem}{Proof}

\newcommand{\Cr}{\EuScript{C}} 
\newcommand{\Gr}{\EuScript{G}} 
 
\newcommand{\Sr}{\EuScript{S}} 
\newcommand{\Nb}{\mathbb{N}} 
\newcommand{\Tb}{\mathbb{T}} 
\newcommand{\Zb}{\mathbb{Z}}

\DeclareMathAlphabet{\mathscr}{T1}{pzc}{m}{it}

\newcommand{\red}[1]{\rightarrow\!\!_{{\scriptscriptstyle #1}}}
\newcommand{\mred}[1]{\twoheadrightarrow\!\!_{{\scriptscriptstyle #1}}}
\newcommand{\mder}[1]{\twoheadleftarrow\!\!_{{\scriptscriptstyle #1}}}
\newcommand{\equi}[1]{\equiv\!\!_{{\scriptscriptstyle #1}}}

\newcommand{\fl}{\rightarrow}
\newcommand{\lfl}{\longrightarrow}

\renewcommand{\phi}{\varphi}
\renewcommand{\epsilon}{\varepsilon}
\newcommand{\tens}{\otimes}
\newcommand{\mon}[1]{\langle #1 \rangle}
\newcommand{\ens}[1]{\{#1\}}
\newcommand{\findem}{\hfill $\diamondsuit$}
\newcommand{\ol}[1]{\overline{#1}}
\newcommand{\et}{\quad\text{and}\quad}
\newcommand{\etll}{\binampersand}
\newcommand{\parll}{\bindnasrepma}

\newdir{ >}{{}*!/-10pt/@{>}}

\DeclareMathOperator{\id}{id}

\newcommand{\bib}[3]
{ 
  \noindent\textsc{#1}
  \newline\indent\emph{#2}, #3.
}

\newcommand{\bibdeux}[5]
{
  \noindent\textsc{#1}
	\newline\indent\emph{#2}, #3. 
	\newline\indent\emph{#4}, #5.
}

\newcommand{\bibtrois}[7]
{
  \noindent\textsc{#1}
	\newline\indent\emph{#2}, #3.
	\newline\indent\emph{#4}, #5.
	\newline\indent\emph{#6}, #7.
}

\newcommand{\emptysectionmark}[1]{\markboth{\textbf{#1}}{\textbf{#1}}}

\fancypagestyle{plain}{
  \fancyhf{}\fancyfoot[LC,RC]{}
  \fancyfoot[LE,RO]{\textbf{\thepage}}
  
  }

\newcommand{\titre}[2]{
\thispagestyle{plain}
\hfill {\large \textbf{#2}} \\

\hrule height 1.5pt
\bigskip
{\LARGE \textbf{#1}}

\bigskip
{\LARGE \textbf{Yves Guiraud}\footnote{Institut de mathématiques de Luminy, Marseille, France - http://iml.univ-mrs.fr/$\sim$guiraud} }
\bigskip
\hrule height 1.5pt
\bigskip
}

\begin{document}
\titre{THE THREE DIMENSIONS OF PROOFS}{22nd July 2005 - Modified 18th November 2005}

\begin{quote}
\textbf{Abstract:} In this document, we study a $3$-polygraphic translation for the proofs of SKS, a formal system for classical propositional logic. We prove that the free $3$-category generated by this $3$-polygraph describes the proofs of classical propositional logic \emph{modulo} structural bureaucracy. We give a $3$-dimensional generalization of Penrose diagrams and use it to provide several pictures of a proof. We sketch how local transformations of proofs yield a non contrived example of $4$-dimensional rewriting.
\end{quote}

\section*{Outline}
\emptysectionmark{Outline}

\noindent In the first section of this paper, we give a $2$-dimensional translation of the formulas of system SKS, a formal system for propositional classical logic [Brünnler 2004] expressed in the style of the calculus of structures [Guglielmi 2004]. The idea consists in the replacement of formulas by circuit-like objects organized in a $2$-polygraph [Burroni 1993]. This construction is formalized in theorem~\ref{th:convergence-gestion-ressources}.

We proceed to section~\ref{sec:preuves}, whose purpose is to translate the proofs of SKS into $3$-dimensional objects that form a $3$-polygraph. There we note that every inference rule can be interpreted as a directed $3$-cell between two circuits. We prove theorem~\ref{thm:1} stating that the $3$-polygraph we have built can be equipped with a proof theory which is the same as the SKS one. Section~\ref{sec:3d} is where the $3$-dimensional nature of proofs happens to be useful: theorem~\ref{thm:2} states that the structural bureaucracy of SKS [Guglielmi~2004] corresponds to topological moves of $3$-cells, called exchange relations.

In section~\ref{sec:dessins-3d} we draw several $3$-dimensional representations of a given proof. Section~\ref{sec:4d} is an informal discussion about the $4$-dimensional nature of local transformations of $3$-dimensional proofs. The final section~\ref{sec:slls} describes how to adapt the work done here to SLLS, the calculus of structures-style formalism for linear logic [Stra\ss{}burger~2003].

\section{The two dimensions of formulas}\label{sec:formules}

\noindent This section gives a $2$-dimensional translation of SKS formulas, heavily inspired by the one already known for terms, studied in [Burroni 1993], [Lafont 2003] and [Guiraud 2004].

After having described the SKS formulas (\ref{sub:formules}), we give the intuition behind their translation into circuit-like objects (\ref{sub:idee-formules-2d}): this works by replacing variables with explicit local resources management operators. This construction requires some theoretical material which is recalled at this moment (\ref{sub:2-poly}). Then we formalize the translation and study its properties (\ref{sub:formules-2d}): the main purpose of this technical part, that can be skipped on a first approach, is to prove that we can compute a canonical representative for circuits corresponding to the same SKS formula (theorem~\ref{th:convergence-gestion-ressources}). Finally we translate the structural congruence on SKS formulas into a congruence on the corresponding circuits (\ref{sub:congruence-2d}).

\subsection{The formulas of SKS}\label{sub:formules}

\noindent System SKS is a formal system for proofs of propositional logic [Br\"unnler 2004]. It is one of the formalisms expressed in the calculus of structures-style, an alternative to sequent calculus where inference rules can be applied at any depth inside formulas [Guglielmi 2004]. Here an alternative definition is used, with a term rewriting vocabulary, such as in [Baader Nipkow 1998].

\begin{defn*}\label{def:formules}
Let us consider two countable sets $V_A$ and $V_F$, which elements are respectively denoted by $a_1$, $a_2$, etc. and $x_1$, $x_2$, etc. The \emph{set of SKS terms} is the set $T$ defined as the disjoint union of the two sets of the pair $(A,F)$ freely generated by the following signature $\Sr$ on the pair $(V_A,V_F)$:
$$
\xymatrix{&& A \ar@(ul,ur)^{\nu} \ar[d]^{\iota} \\ \ast \ar@<0.5ex>[rr]^{\top} \ar@<-0.5ex>[rr]_{\bot} && F && F\times F. \ar@<-0.5ex>[ll]_{\wedge} \ar@<0.5ex>[ll]^{\vee}}
$$

\noindent Terms of sort $A$ are called \emph{SKS atoms} and terms of sort $F$ are called \emph{SKS formulas}. The binary relation denoted by $\equi{S}$ is defined as the congruence on SKS terms generated by the following rewriting rules:
$$
\begin{array}{r c l c r c l}
(x_1\wedge x_2)\wedge x_3 &\lfl& x_1\wedge (x_2\wedge x_3) &&
(x_1\vee x_2)\vee x_3 &\lfl& x_1\vee (x_2\vee x_3) \\
x_1\wedge x_2 &\lfl& x_2\wedge x_1 &&
x_1\vee x_2 &\lfl& x_2\vee x_1 \\
\top\wedge x_1 &\lfl& x_1 &&
\bot\vee x_1 &\lfl& x_1 \\
\bot\wedge\bot &\lfl& \bot &&
\top\vee\top &\lfl& \top \\
&& \hfill\nu(\nu(a_1)) &\lfl& a_1.\hfill\strut
\end{array}
$$
\end{defn*}

\begin{rem*}
The binary relation $\equi{S}$ is defined in three steps:
\begin{enumerate}
\item[1.] One defines the reduction relation $\red{S}$ on terms by $u\red{S}v$ if there exist a context $C$, a substitution~$\sigma$ and one of the nine above rules $\alpha:s(\alpha)\fl t(\alpha)$ such that $u=C[s(\alpha)\cdot\sigma]$ and $v=C[t(\alpha)\cdot\sigma]$. As usual, $C[u]$ denotes the application of a context $C$ to a term $u$, while $u\cdot\sigma$ stands for the application of a substitution $\sigma$ to a term $u$.
\item[2.] Then, one defines the relation $\mred{S}$ from $\red{S}$ by $u\mred{S}v$ if $u=v$ or if there exists a possibly empty family $(u_1,\dots,u_n)$ of terms such that:
$$
u\red{S}u_1\red{S}u_2\red{S}\dots\red{S}u_n\red{S}v.
$$

\item[3.] Finally, one defines the relation $\equi{S}$ by $u\equi{S}v$ if there exists a possibly empty family $(u_1,\dots,u_{2n})$ of terms such that:
$$
u\mred{S}u_1\mder{S}u_2\mred{S}\dots\mder{S}u_{2n}\mred{S}v.
$$
\end{enumerate}

\noindent Let us note that, \emph{modulo $\equi{S}$}, the pairs $(\wedge,\top)$ and $(\vee,\bot)$ are commutative monoid structures on the set of SKS terms and that the map $\nu$ is an involution.
\end{rem*}

\begin{rem*}
As they are defined here, the SKS terms are more general than the original SKS formulas of [Br\"unnler 2004]. It is straightforward to check that the original formulas are the closed SKS terms of sort $F$, \emph{modulo} the relation~$\equi{S}$.

The SKS terms described here are more convenient for many reasons, among which the possibility to reduce the inference rules to a finite number. However, this generalization allows non-linear terms: this is where we need results from [Burroni 1993] to translate terms into circuits, as described in the rest of this section.

Another choice could have been made: replacing variables and their negations by a countable number of constants. This would simplify the translations of terms, since one would need only one sort ($F$) and one resources management operator ($\tau_{FF}$, defined thereafter). The main drawback of this choice is that it requires a countable number of $3$-dimensional cells to translate the inference rules, in addition to the countable number of $2$-dimensional cells for variables.
\end{rem*}

\subsection{From formulas to circuits: the informal idea}\label{sub:idee-formules-2d}

\noindent The translation of terms into $2$-dimensional objects has been developped troughout [Burroni 1993], [Lafont 2003] and [Guiraud 2004]. The idea is to replace each (family of) term(s) with a circuit: it is built with the tree-part of the term with, plugged in the leaves, an additional part replacing variables and consisting of local resources management operators.

Before any formalization, let us give a few examples:
\begin{center}
\begin{picture}(0,0)%
\includegraphics{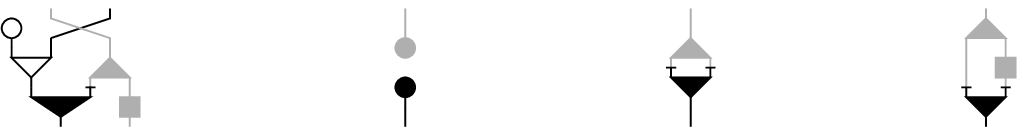}%
\end{picture}%
\setlength{\unitlength}{4144sp}%
\begingroup\makeatletter\ifx\SetFigFont\undefined%
\gdef\SetFigFont#1#2#3#4#5{%
  \reset@font\fontsize{#1}{#2pt}%
  \fontfamily{#3}\fontseries{#4}\fontshape{#5}%
  \selectfont}%
\fi\endgroup%
\begin{picture}(4655,961)(-7,-200)
\put(4636,-151){\makebox(0,0)[b]{\smash{{\SetFigFont{10}{12.0}{\rmdefault}{\mddefault}{\updefault}{\color[rgb]{0,0,0}$\iota(a_1)\wedge\iota(\nu(a_1))$}%
}}}}
\put(316,-151){\makebox(0,0)[b]{\smash{{\SetFigFont{10}{12.0}{\rmdefault}{\mddefault}{\updefault}{\color[rgb]{0,0,0}$\big((\bot\vee x_1)\wedge\iota(a_1),\nu(a_1)\big)$}%
}}}}
\put(1846,-151){\makebox(0,0)[b]{\smash{{\SetFigFont{10}{12.0}{\rmdefault}{\mddefault}{\updefault}{\color[rgb]{0,0,0}$\top$}%
}}}}
\put(3151,-151){\makebox(0,0)[b]{\smash{{\SetFigFont{10}{12.0}{\rmdefault}{\mddefault}{\updefault}{\color[rgb]{0,0,0}$\iota(a_1)\wedge\iota(a_1)$}%
}}}}
\end{picture}%
\end{center}

\noindent These circuits are built using two kinds of wires (one for formulas, in black, and one for atoms, in grey) and the following fourteen components (six corresponding to the terms constructors and eight for explicit resources management):
\begin{center}
\begin{picture}(0,0)%
\includegraphics{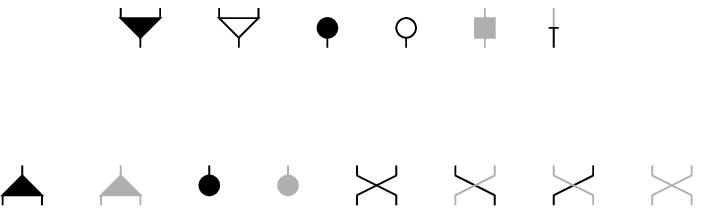}%
\end{picture}%
\setlength{\unitlength}{4144sp}%
\begingroup\makeatletter\ifx\SetFigFont\undefined%
\gdef\SetFigFont#1#2#3#4#5{%
  \reset@font\fontsize{#1}{#2pt}%
  \fontfamily{#3}\fontseries{#4}\fontshape{#5}%
  \selectfont}%
\fi\endgroup%
\begin{picture}(3174,1231)(-461,-470)
\put(181,299){\makebox(0,0)[b]{\smash{{\SetFigFont{10}{12.0}{\rmdefault}{\mddefault}{\updefault}{\color[rgb]{0,0,0}$\wedge$}%
}}}}
\put(2611,-421){\makebox(0,0)[b]{\smash{{\SetFigFont{10}{12.0}{\rmdefault}{\mddefault}{\updefault}{\color[rgb]{0,0,0}$\tau_{AA}$}%
}}}}
\put(2161,-421){\makebox(0,0)[b]{\smash{{\SetFigFont{10}{12.0}{\rmdefault}{\mddefault}{\updefault}{\color[rgb]{0,0,0}$\tau_{AF}$}%
}}}}
\put(1711,-421){\makebox(0,0)[b]{\smash{{\SetFigFont{10}{12.0}{\rmdefault}{\mddefault}{\updefault}{\color[rgb]{0,0,0}$\tau_{FA}$}%
}}}}
\put(1261,-421){\makebox(0,0)[b]{\smash{{\SetFigFont{10}{12.0}{\rmdefault}{\mddefault}{\updefault}{\color[rgb]{0,0,0}$\tau_{FF}$}%
}}}}
\put( 91,-421){\makebox(0,0)[b]{\smash{{\SetFigFont{10}{12.0}{\rmdefault}{\mddefault}{\updefault}{\color[rgb]{0,0,0}$\delta_A$}%
}}}}
\put(856,-421){\makebox(0,0)[b]{\smash{{\SetFigFont{10}{12.0}{\rmdefault}{\mddefault}{\updefault}{\color[rgb]{0,0,0}$\epsilon_A$}%
}}}}
\put(-359,-421){\makebox(0,0)[b]{\smash{{\SetFigFont{10}{12.0}{\rmdefault}{\mddefault}{\updefault}{\color[rgb]{0,0,0}$\delta_F$}%
}}}}
\put(496,-421){\makebox(0,0)[b]{\smash{{\SetFigFont{10}{12.0}{\rmdefault}{\mddefault}{\updefault}{\color[rgb]{0,0,0}$\epsilon_F$}%
}}}}
\put(2071,299){\makebox(0,0)[b]{\smash{{\SetFigFont{10}{12.0}{\rmdefault}{\mddefault}{\updefault}{\color[rgb]{0,0,0}$\iota$}%
}}}}
\put(1756,299){\makebox(0,0)[b]{\smash{{\SetFigFont{10}{12.0}{\rmdefault}{\mddefault}{\updefault}{\color[rgb]{0,0,0}$\nu$}%
}}}}
\put(1396,299){\makebox(0,0)[b]{\smash{{\SetFigFont{10}{12.0}{\rmdefault}{\mddefault}{\updefault}{\color[rgb]{0,0,0}$\bot$}%
}}}}
\put(1036,299){\makebox(0,0)[b]{\smash{{\SetFigFont{10}{12.0}{\rmdefault}{\mddefault}{\updefault}{\color[rgb]{0,0,0}$\top$}%
}}}}
\put(631,299){\makebox(0,0)[b]{\smash{{\SetFigFont{10}{12.0}{\rmdefault}{\mddefault}{\updefault}{\color[rgb]{0,0,0}$\vee$}%
}}}}
\end{picture}%
\end{center}

\noindent Two operations are allowed to build the circuits, one for each dimension (note that the second one is only defined if the circuits inputs/outputs match):
\begin{center}
\begin{picture}(0,0)%
\includegraphics{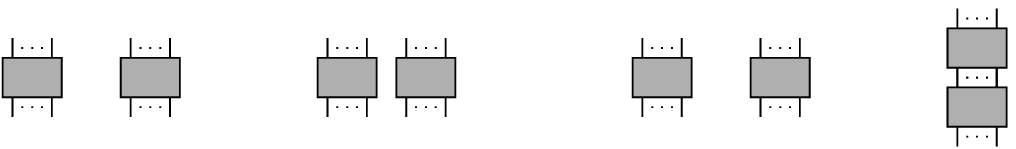}%
\end{picture}%
\setlength{\unitlength}{4144sp}%
\begingroup\makeatletter\ifx\SetFigFont\undefined%
\gdef\SetFigFont#1#2#3#4#5{%
  \reset@font\fontsize{#1}{#2pt}%
  \fontfamily{#3}\fontseries{#4}\fontshape{#5}%
  \selectfont}%
\fi\endgroup%
\begin{picture}(4614,654)(79,-568)
\put(3376,-286){\makebox(0,0)[b]{\smash{{\SetFigFont{10}{12.0}{\rmdefault}{\mddefault}{\updefault}{\color[rgb]{0,0,0}$\circ$}%
}}}}
\put(226,-286){\makebox(0,0)[b]{\smash{{\SetFigFont{10}{12.0}{\rmdefault}{\mddefault}{\updefault}{\color[rgb]{0,0,0}$f$}%
}}}}
\put(1666,-286){\makebox(0,0)[b]{\smash{{\SetFigFont{10}{12.0}{\rmdefault}{\mddefault}{\updefault}{\color[rgb]{0,0,0}$f$}%
}}}}
\put(2026,-286){\makebox(0,0)[b]{\smash{{\SetFigFont{10}{12.0}{\rmdefault}{\mddefault}{\updefault}{\color[rgb]{0,0,0}$g$}%
}}}}
\put(4546,-151){\makebox(0,0)[b]{\smash{{\SetFigFont{10}{12.0}{\rmdefault}{\mddefault}{\updefault}{\color[rgb]{0,0,0}$f$}%
}}}}
\put(4546,-421){\makebox(0,0)[b]{\smash{{\SetFigFont{10}{12.0}{\rmdefault}{\mddefault}{\updefault}{\color[rgb]{0,0,0}$g$}%
}}}}
\put(766,-286){\makebox(0,0)[b]{\smash{{\SetFigFont{10}{12.0}{\rmdefault}{\mddefault}{\updefault}{\color[rgb]{0,0,0}$g$}%
}}}}
\put(3646,-286){\makebox(0,0)[b]{\smash{{\SetFigFont{10}{12.0}{\rmdefault}{\mddefault}{\updefault}{\color[rgb]{0,0,0}$f$}%
}}}}
\put(3106,-286){\makebox(0,0)[b]{\smash{{\SetFigFont{10}{12.0}{\rmdefault}{\mddefault}{\updefault}{\color[rgb]{0,0,0}$g$}%
}}}}
\put(496,-286){\makebox(0,0)[b]{\smash{{\SetFigFont{10}{12.0}{\rmdefault}{\mddefault}{\updefault}{\color[rgb]{0,0,0}$\tens$}%
}}}}
\put(1216,-286){\makebox(0,0)[b]{\smash{{\SetFigFont{10}{12.0}{\rmdefault}{\mddefault}{\updefault}{\color[rgb]{0,0,0}$=$}%
}}}}
\put(4096,-286){\makebox(0,0)[b]{\smash{{\SetFigFont{10}{12.0}{\rmdefault}{\mddefault}{\updefault}{\color[rgb]{0,0,0}$=$}%
}}}}
\end{picture}%
\end{center}

\noindent Usual alternative notations include $f\star_0 g$ for $f\tens g$, $f\star_1 g$ for $g\circ f$. The circuits are seen as topological objects and, as such, considered \emph{modulo} homeomorphic deformation. This means that wires can be lengthened or shortened and that components can be moved, provided no crossing of wires is created, such as in the following:
\begin{center}
\begin{picture}(0,0)%
\includegraphics{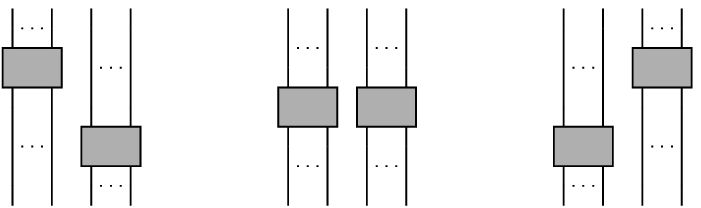}%
\end{picture}%
\setlength{\unitlength}{4144sp}%
\begingroup\makeatletter\ifx\SetFigFont\undefined%
\gdef\SetFigFont#1#2#3#4#5{%
  \reset@font\fontsize{#1}{#2pt}%
  \fontfamily{#3}\fontseries{#4}\fontshape{#5}%
  \selectfont}%
\fi\endgroup%
\begin{picture}(3174,924)(79,-73)
\put(2296,344){\makebox(0,0)[b]{\smash{{\SetFigFont{10}{12.0}{\rmdefault}{\mddefault}{\updefault}{\color[rgb]{0,0,0}$=$}%
}}}}
\put(226,524){\makebox(0,0)[b]{\smash{{\SetFigFont{10}{12.0}{\rmdefault}{\mddefault}{\updefault}{\color[rgb]{0,0,0}$f$}%
}}}}
\put(586,164){\makebox(0,0)[b]{\smash{{\SetFigFont{10}{12.0}{\rmdefault}{\mddefault}{\updefault}{\color[rgb]{0,0,0}$g$}%
}}}}
\put(2746,164){\makebox(0,0)[b]{\smash{{\SetFigFont{10}{12.0}{\rmdefault}{\mddefault}{\updefault}{\color[rgb]{0,0,0}$f$}%
}}}}
\put(3106,524){\makebox(0,0)[b]{\smash{{\SetFigFont{10}{12.0}{\rmdefault}{\mddefault}{\updefault}{\color[rgb]{0,0,0}$g$}%
}}}}
\put(1486,344){\makebox(0,0)[b]{\smash{{\SetFigFont{10}{12.0}{\rmdefault}{\mddefault}{\updefault}{\color[rgb]{0,0,0}$f$}%
}}}}
\put(1846,344){\makebox(0,0)[b]{\smash{{\SetFigFont{10}{12.0}{\rmdefault}{\mddefault}{\updefault}{\color[rgb]{0,0,0}$g$}%
}}}}
\put(1036,344){\makebox(0,0)[b]{\smash{{\SetFigFont{10}{12.0}{\rmdefault}{\mddefault}{\updefault}{\color[rgb]{0,0,0}$=$}%
}}}}
\end{picture}%
\end{center}

\noindent In [Burroni 1993], this kind of collection of circuits was given a name: a \emph{$\mathit{2}$-polygraph}.

\subsection{Two-polygraphs and two-categories}\label{sub:2-poly}

In order to define this structure, we recall some notions about graphs and free categories.

\begin{ntt*}
If $G$ is a graph, its set of objects is denoted by $G_0$ and its set of arrows going from an object~$x$ to another object $y$ is denoted by $G(x,y)$; for such an arrow $f$, $s_0(f)$ is the source $x$ of $f$ and~$t_0(f)$ its target $y$. The set of all arrows of $G$ is denoted by $G_1$ and $G$ itself is often abusively denoted by~$(G_0,G_1)$ only, assuming that the source and target mappings are given with $G_1$.
\end{ntt*}

\begin{defn*}\label{def:categorie_libre}
Let $G=(G_0,G_1)$ be a graph. The \emph{free category generated by $G$}, denoted by $\mon{G}$, is the following (small) category:
\begin{enumerate}
\item[0.] The objects of $\mon{G}$ are the objects of $G$.
\item[1.] The arrows of $\mon{G}$ from $x$ to $y$, are all the finite paths in $G$ going from $x$ to $y$, with concatenation~$\circ$ as composition and empty paths as local identities. The set of all arrows of $\mon{G}$ is denoted by~$\mon{G}_1$.
\end{enumerate}
\end{defn*}

\begin{defn*}\label{def:2_polygraphe}
A \emph{$\mathit{2}$-polygraph} $\Sigma$ is given by:
\begin{enumerate}
\item[0.] A set $\Sigma_0$ of \emph{$\mathit{0}$-cells}.
\item[1.] A set $\Sigma_1$ of \emph{$\mathit{1}$-cells}, together with two maps $s_0,t_0:\Sigma_1\fl\Sigma_0$, called \emph{$\mathit{0}$-source} and \emph{$\mathit{0}$-target}. The arrows of the free category $(\Sigma_0,\mon{\Sigma}_1)$ are called \emph{$\mathit{1}$-arrows}. The composition of $f$ followed by $g$ is denoted by $f\star_0g$ in the general case and by $f\tens g$ when $\Sigma_0$ has only one element.
\item[2.] A set $\Sigma_2$ of \emph{$\mathit{2}$-cells}, together with two maps $s_1,t_1:\Sigma_2\fl\mon{\Sigma}_1$, called \emph{$\mathit{1}$-source} and \emph{$\mathit{1}$-target}, and such that $s_0\circ s_1=s_0\circ t_1$ and $t_0\circ s_1 = t_0\circ t_1$. The first equality gives a map $s_0:\Sigma_2\fl\Sigma_0$ and the second one yields $t_0:\Sigma_2\fl\Sigma_0$.
\end{enumerate}
\end{defn*}

\noindent Thus, in order to translate formulas, we build a $2$-polygraph $\Sigma^F$ with one $0$-cell $\ast$ (this one can be seen as the background color in the graphical representations), two $1$-cells $A$ and $F$ (the two colors of wires) and fourteen $2$-cells (the circuit components).

All the circuits that can be built with the elementary bricks given by $\Sigma^F$, equipped with their two compositions, considered \emph{modulo} homeomorphic deformation, form the \emph{$\mathit{2}$-arrows of the free $\mathit{2}$-category~$\mon{\Sigma^F}$ generated by $\Sigma^F$}. The set of all the $2$-arrows of $\mon{\Sigma^F}$ is denoted by $\mon{\Sigma^F}_2$.

\begin{rem*}
We do not give here a complete definition of this notion, which can be found in either of [Burroni 1993], [Métayer 2003] or [Guiraud 2005]. Other sources of information about this topic include [MacLane 1998] for $2$-categories, [Baez Dolan 1998] for a certain kind of higher-dimensional categories and [Chang Lauda 2004] for a whole zoo of them. For this document, let us say that every $2$-category we are interested in can be seen as the quotient of a free $2$-category (generated by some $2$-polygraph) by equations between parallel $2$-arrows ($2$-arrows that have the same $1$-source and the same $1$-target).
\end{rem*}

\subsection{From formulas to circuits: the formal construction}\label{sub:formules-2d}

\noindent In this paragraph, we build translations between terms and circuits. We follow the same path as in [Guiraud 2004]: the results we seek are the same as in that document, except for the generalization to the two-sorted case. We start with the construction of the $2$-category $\Tb$ of terms, built from the set $T$ of SKS terms. First of all, we give some useful notations:

\begin{ntt*}
Let $X=X_1\tens\dots\tens X_n$ be a $1$-arrow in $\mon{\Sigma}$: each $X_i$ is either $A$ or $F$. We denote by~$\sharp X$ the pair $(\sharp_A X,\sharp_F X)$ of natural numbers such that $\sharp_A X$ (resp. $\sharp_F X$) is the number of $A$ (resp. $F$) appearing in $X$. If $u=(u_1,\dots,u_n)$ is a family of $n$ terms in $T$, we denote by $u_A$ (resp. $u_F$) the subfamily of $u$ consisting only of the $u_i$ in $A$ (resp. in $F$), appearing in the same order as in $u$. We denote by $\sharp u$ the pair~$(\sharp_A u,\sharp_F u)$ of natural numbers defined by: $\sharp_A u$ (resp. $\sharp_F u$) is the greatest of the natural numbers $k$ such that the variable $a_k$ (resp. $x_k$) appears in at least one of the terms $u_1$, $\dots$, $u_n$. Two pairs of natural numbers are compared with the product order given by the natural one on $\Nb$.
\end{ntt*}

\begin{defn*}
Let us define the \emph{$\mathit{2}$-category of terms}, denoted by $\Tb$, as follows:
\begin{enumerate}
\item[0.] It contains one $0$-arrow, denoted by $\ast$.
\item[1.] Its $1$-arrows are the elements $X_1\tens\dots\tens X_n$ of the free monoid generated by $\ens{A,F}$.
\item[2.] If $X$ and $Y=Y_1\tens\dots\tens Y_n$ are two $1$-arrows, then the $2$-arrows of $\Tb$ from $X$ to $Y$ are all the families $u=(u_1,\dots,u_n)$ of $n$ terms such that each $u_i$ is in $Y_i$ and such that $\sharp u\leq\sharp X$.
\end{enumerate}

\noindent The two compositions are given by:
\begin{enumerate}
\item[-] If $u=(u_1,\dots,u_n)$ is a $2$-arrow from $X$ to $Y$ and $v=(v_1,\dots,v_q)$ is a $2$-arrow from $X'$ to $Y'$, then their product $u\tens v$ is the $2$-arrow from $X\tens X'$ to $Y\tens Y'$ defined by:
$$
u\tens v \:=\: \big(u_1,\dots,u_n,v_1\cdot\rho_{\sharp u},\dots,v_q\cdot\rho_{\sharp u}\big),
$$

\noindent where $\rho_{\sharp u}$ is the substitution that sends each $a_i$ onto $a_{i+\sharp_A u}$ and each $x_i$ onto $x_{i+\sharp_F u}$.

\item[-] If $u=(u_1,\dots,u_n)$ is a $2$-arrow from $X$ to $Y$ and $v=(v_1,\dots,v_p)$ is a $2$-arrow from $Y$ to $Z$, then their composite $v\circ u$ is the $2$-arrow $(w_1,\dots,w_p)$ from $X$ to $Z$ such that $w_i$ is $v_i$ where each $a_k$ (resp. $x_k$) is replaced by the $k^{\text{th}}$ element of $u_A$ (resp. $u_F$).
\end{enumerate}
\end{defn*}

\begin{rem*}
One must check that the operations $\circ$ and $\tens$ are well-defined and that they satisfy the axioms for the structure of $2$-category, as given in [MacLane 1998] for example.
\end{rem*}

\noindent A family $u=(u_1,\dots,u_n)$ of terms can be seen as many $2$-arrows in $\Tb$. Indeed, let us assume that $\sharp u=(m,n)$. Then, for any $p\geq m$ and $q\geq n$, $u$ can be seen as a $2$-arrow with source $A^p\times F^q$: this means seeing $u$ as using more variables than it seems (these are dummy variables). Furthermore, one can also shuffle the source $A^p\times F^q$ and still see $u$ as a $2$-arrow with source the result of this shuffle. On the other hand, the target of all these $2$-arrows is always the same: it is entirely and uniquely fixed by the sorts of each $u_i$.

\begin{ex*}
Let us consider $u=(a_3,x_2\wedge x_3,\nu a_1)$, seen as a $2$-arrow from $A^3\tens F^3$ to $A\tens F\tens A$, and $v=(\iota a_2\wedge x_1, x_1)$, seen as a $2$-arrow from $A^2\tens F$ to $F^2$. Then $u\tens v$ and $v\circ u$ are:
$$
u\tens v \:=\: (a_3,x_2\wedge x_3,\nu a_1, \iota a_5\wedge x_4,x_4) \et v\circ u \:=\: (\iota\nu a_1 \wedge (x_2\wedge x_3),x_2\wedge x_3).
$$

\noindent Note that, if $u$ was considered as an arrow with one dummy variable of type $F$, for example from $A^3\tens F^4$ to $A\tens F\tens A$, then the result of $v\circ u$ would not be changed (except from its source), while $u\tens v$ would become:
$$
u\tens v \:=\: (a_3,x_2\wedge x_3,\nu a_1, \iota a_5\wedge x_5, x_5).
$$

\noindent On the other hand, the result would not change if only the source $A^3\tens F^3$ was shuffled, into the $1$-arrow $A\tens F^2\tens A\tens F\tens A$ for example.
\end{ex*}

\noindent Now we want to prove that $\Tb$ has a graphical presentation as a quotient of a free $2$-category. We use a result from [Burroni 1993] which requires the following notations:

\begin{ntt*}
We recall that $\Sigma^F$ is the following $2$-polygraph:
\begin{center}
\begin{picture}(0,0)%
\includegraphics{deux_cellules_sks.eps}%
\end{picture}%
\setlength{\unitlength}{4144sp}%
\begingroup\makeatletter\ifx\SetFigFont\undefined%
\gdef\SetFigFont#1#2#3#4#5{%
  \reset@font\fontsize{#1}{#2pt}%
  \fontfamily{#3}\fontseries{#4}\fontshape{#5}%
  \selectfont}%
\fi\endgroup%
\begin{picture}(3174,1231)(-461,-470)
\put(181,299){\makebox(0,0)[b]{\smash{{\SetFigFont{10}{12.0}{\rmdefault}{\mddefault}{\updefault}{\color[rgb]{0,0,0}$\wedge$}%
}}}}
\put(2611,-421){\makebox(0,0)[b]{\smash{{\SetFigFont{10}{12.0}{\rmdefault}{\mddefault}{\updefault}{\color[rgb]{0,0,0}$\tau_{AA}$}%
}}}}
\put(2161,-421){\makebox(0,0)[b]{\smash{{\SetFigFont{10}{12.0}{\rmdefault}{\mddefault}{\updefault}{\color[rgb]{0,0,0}$\tau_{AF}$}%
}}}}
\put(1711,-421){\makebox(0,0)[b]{\smash{{\SetFigFont{10}{12.0}{\rmdefault}{\mddefault}{\updefault}{\color[rgb]{0,0,0}$\tau_{FA}$}%
}}}}
\put(1261,-421){\makebox(0,0)[b]{\smash{{\SetFigFont{10}{12.0}{\rmdefault}{\mddefault}{\updefault}{\color[rgb]{0,0,0}$\tau_{FF}$}%
}}}}
\put( 91,-421){\makebox(0,0)[b]{\smash{{\SetFigFont{10}{12.0}{\rmdefault}{\mddefault}{\updefault}{\color[rgb]{0,0,0}$\delta_A$}%
}}}}
\put(856,-421){\makebox(0,0)[b]{\smash{{\SetFigFont{10}{12.0}{\rmdefault}{\mddefault}{\updefault}{\color[rgb]{0,0,0}$\epsilon_A$}%
}}}}
\put(-359,-421){\makebox(0,0)[b]{\smash{{\SetFigFont{10}{12.0}{\rmdefault}{\mddefault}{\updefault}{\color[rgb]{0,0,0}$\delta_F$}%
}}}}
\put(496,-421){\makebox(0,0)[b]{\smash{{\SetFigFont{10}{12.0}{\rmdefault}{\mddefault}{\updefault}{\color[rgb]{0,0,0}$\epsilon_F$}%
}}}}
\put(2071,299){\makebox(0,0)[b]{\smash{{\SetFigFont{10}{12.0}{\rmdefault}{\mddefault}{\updefault}{\color[rgb]{0,0,0}$\iota$}%
}}}}
\put(1756,299){\makebox(0,0)[b]{\smash{{\SetFigFont{10}{12.0}{\rmdefault}{\mddefault}{\updefault}{\color[rgb]{0,0,0}$\nu$}%
}}}}
\put(1396,299){\makebox(0,0)[b]{\smash{{\SetFigFont{10}{12.0}{\rmdefault}{\mddefault}{\updefault}{\color[rgb]{0,0,0}$\bot$}%
}}}}
\put(1036,299){\makebox(0,0)[b]{\smash{{\SetFigFont{10}{12.0}{\rmdefault}{\mddefault}{\updefault}{\color[rgb]{0,0,0}$\top$}%
}}}}
\put(631,299){\makebox(0,0)[b]{\smash{{\SetFigFont{10}{12.0}{\rmdefault}{\mddefault}{\updefault}{\color[rgb]{0,0,0}$\vee$}%
}}}}
\end{picture}%
\end{center}

\noindent We denote by $E_{\Delta}$ the union of the following two families of relations on parallel $2$-arrows of the free $2$-category $\mon{\Sigma^F}$:
\begin{enumerate}
\item[1.] The first family is made of 26 relations, given by all the possible colorations of wires of the following diagrams:
\begin{center}
\begin{picture}(0,0)%
\includegraphics{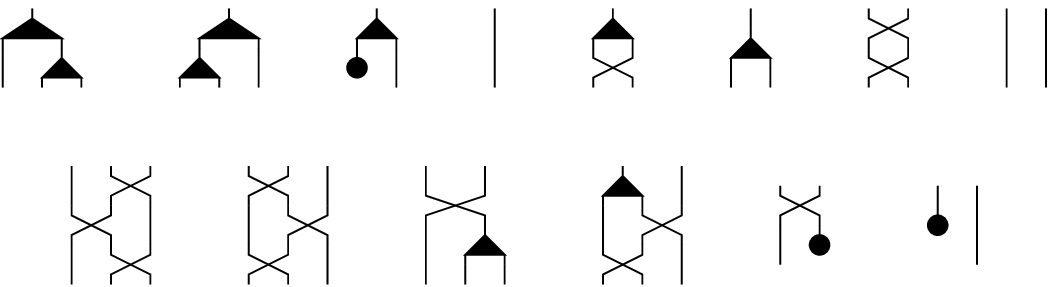}%
\end{picture}%
\setlength{\unitlength}{4144sp}%
\begingroup\makeatletter\ifx\SetFigFont\undefined%
\gdef\SetFigFont#1#2#3#4#5{%
  \reset@font\fontsize{#1}{#2pt}%
  \fontfamily{#3}\fontseries{#4}\fontshape{#5}%
  \selectfont}%
\fi\endgroup%
\begin{picture}(4794,1284)(7639,-4933)
\put(12016,-3886){\makebox(0,0)[b]{\smash{{\SetFigFont{12}{14.4}{\familydefault}{\mddefault}{\updefault}{\color[rgb]{0,0,0}$=$}%
}}}}
\put(11656,-4696){\makebox(0,0)[b]{\smash{{\SetFigFont{12}{14.4}{\familydefault}{\mddefault}{\updefault}{\color[rgb]{0,0,0}$=$}%
}}}}
\put(10126,-4696){\makebox(0,0)[b]{\smash{{\SetFigFont{12}{14.4}{\familydefault}{\mddefault}{\updefault}{\color[rgb]{0,0,0}$=$}%
}}}}
\put(8551,-4696){\makebox(0,0)[b]{\smash{{\SetFigFont{12}{14.4}{\familydefault}{\mddefault}{\updefault}{\color[rgb]{0,0,0}$=$}%
}}}}
\put(8236,-3886){\makebox(0,0)[b]{\smash{{\SetFigFont{12}{14.4}{\familydefault}{\mddefault}{\updefault}{\color[rgb]{0,0,0}$=$}%
}}}}
\put(9676,-3886){\makebox(0,0)[b]{\smash{{\SetFigFont{12}{14.4}{\familydefault}{\mddefault}{\updefault}{\color[rgb]{0,0,0}$=$}%
}}}}
\put(10756,-3886){\makebox(0,0)[b]{\smash{{\SetFigFont{12}{14.4}{\familydefault}{\mddefault}{\updefault}{\color[rgb]{0,0,0}$=$}%
}}}}
\end{picture}%
\end{center}

\vfill\pagebreak
\item[2.] The second family is made of 24 relations, four for each of $\wedge$, $\vee$, $\top$, $\bot$, $\iota$, $\nu$:
\begin{center}
\begin{picture}(0,0)%
\includegraphics{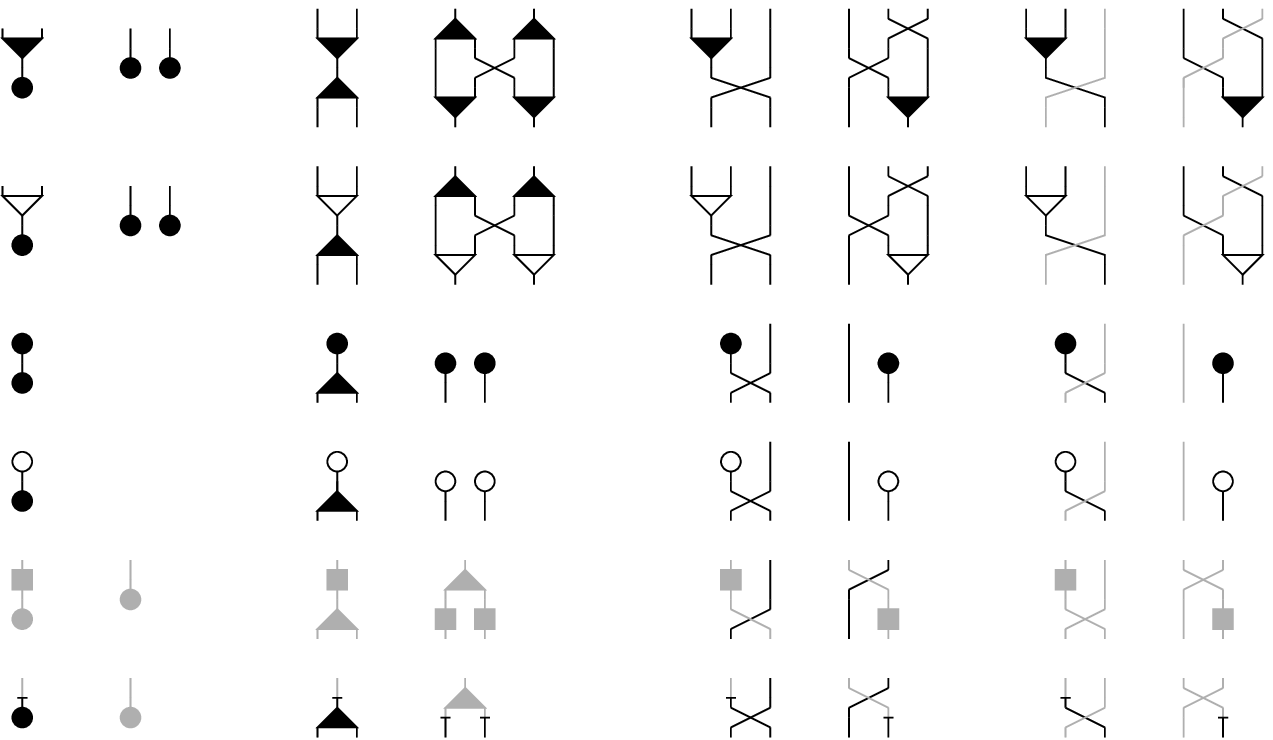}%
\end{picture}%
\setlength{\unitlength}{4144sp}%
\begingroup\makeatletter\ifx\SetFigFont\undefined%
\gdef\SetFigFont#1#2#3#4#5{%
  \reset@font\fontsize{#1}{#2pt}%
  \fontfamily{#3}\fontseries{#4}\fontshape{#5}%
  \selectfont}%
\fi\endgroup%
\begin{picture}(5784,3354)(169,-3223)
\put(5401,-2626){\makebox(0,0)[b]{\smash{{\SetFigFont{10}{12.0}{\rmdefault}{\mddefault}{\updefault}{\color[rgb]{0,0,0}$=$}%
}}}}
\put(541,-196){\makebox(0,0)[b]{\smash{{\SetFigFont{10}{12.0}{\rmdefault}{\mddefault}{\updefault}{\color[rgb]{0,0,0}$=$}%
}}}}
\put(1981,-196){\makebox(0,0)[b]{\smash{{\SetFigFont{10}{12.0}{\rmdefault}{\mddefault}{\updefault}{\color[rgb]{0,0,0}$=$}%
}}}}
\put(1981,-916){\makebox(0,0)[b]{\smash{{\SetFigFont{10}{12.0}{\rmdefault}{\mddefault}{\updefault}{\color[rgb]{0,0,0}$=$}%
}}}}
\put(3871,-916){\makebox(0,0)[b]{\smash{{\SetFigFont{10}{12.0}{\rmdefault}{\mddefault}{\updefault}{\color[rgb]{0,0,0}$=$}%
}}}}
\put(3871,-196){\makebox(0,0)[b]{\smash{{\SetFigFont{10}{12.0}{\rmdefault}{\mddefault}{\updefault}{\color[rgb]{0,0,0}$=$}%
}}}}
\put(5401,-196){\makebox(0,0)[b]{\smash{{\SetFigFont{10}{12.0}{\rmdefault}{\mddefault}{\updefault}{\color[rgb]{0,0,0}$=$}%
}}}}
\put(5401,-916){\makebox(0,0)[b]{\smash{{\SetFigFont{10}{12.0}{\rmdefault}{\mddefault}{\updefault}{\color[rgb]{0,0,0}$=$}%
}}}}
\put(541,-3121){\makebox(0,0)[b]{\smash{{\SetFigFont{10}{12.0}{\rmdefault}{\mddefault}{\updefault}{\color[rgb]{0,0,0}$=$}%
}}}}
\put(1981,-3121){\makebox(0,0)[b]{\smash{{\SetFigFont{10}{12.0}{\rmdefault}{\mddefault}{\updefault}{\color[rgb]{0,0,0}$=$}%
}}}}
\put(3871,-3121){\makebox(0,0)[b]{\smash{{\SetFigFont{10}{12.0}{\rmdefault}{\mddefault}{\updefault}{\color[rgb]{0,0,0}$=$}%
}}}}
\put(5401,-3121){\makebox(0,0)[b]{\smash{{\SetFigFont{10}{12.0}{\rmdefault}{\mddefault}{\updefault}{\color[rgb]{0,0,0}$=$}%
}}}}
\put(541,-916){\makebox(0,0)[b]{\smash{{\SetFigFont{10}{12.0}{\rmdefault}{\mddefault}{\updefault}{\color[rgb]{0,0,0}$=$}%
}}}}
\put(541,-1546){\makebox(0,0)[b]{\smash{{\SetFigFont{10}{12.0}{\rmdefault}{\mddefault}{\updefault}{\color[rgb]{0,0,0}$=$}%
}}}}
\put(541,-2086){\makebox(0,0)[b]{\smash{{\SetFigFont{10}{12.0}{\rmdefault}{\mddefault}{\updefault}{\color[rgb]{0,0,0}$=$}%
}}}}
\put(541,-2626){\makebox(0,0)[b]{\smash{{\SetFigFont{10}{12.0}{\rmdefault}{\mddefault}{\updefault}{\color[rgb]{0,0,0}$=$}%
}}}}
\put(1981,-1546){\makebox(0,0)[b]{\smash{{\SetFigFont{10}{12.0}{\rmdefault}{\mddefault}{\updefault}{\color[rgb]{0,0,0}$=$}%
}}}}
\put(1981,-2086){\makebox(0,0)[b]{\smash{{\SetFigFont{10}{12.0}{\rmdefault}{\mddefault}{\updefault}{\color[rgb]{0,0,0}$=$}%
}}}}
\put(1981,-2626){\makebox(0,0)[b]{\smash{{\SetFigFont{10}{12.0}{\rmdefault}{\mddefault}{\updefault}{\color[rgb]{0,0,0}$=$}%
}}}}
\put(3871,-1546){\makebox(0,0)[b]{\smash{{\SetFigFont{10}{12.0}{\rmdefault}{\mddefault}{\updefault}{\color[rgb]{0,0,0}$=$}%
}}}}
\put(3871,-2086){\makebox(0,0)[b]{\smash{{\SetFigFont{10}{12.0}{\rmdefault}{\mddefault}{\updefault}{\color[rgb]{0,0,0}$=$}%
}}}}
\put(3871,-2626){\makebox(0,0)[b]{\smash{{\SetFigFont{10}{12.0}{\rmdefault}{\mddefault}{\updefault}{\color[rgb]{0,0,0}$=$}%
}}}}
\put(5401,-1546){\makebox(0,0)[b]{\smash{{\SetFigFont{10}{12.0}{\rmdefault}{\mddefault}{\updefault}{\color[rgb]{0,0,0}$=$}%
}}}}
\put(5401,-2086){\makebox(0,0)[b]{\smash{{\SetFigFont{10}{12.0}{\rmdefault}{\mddefault}{\updefault}{\color[rgb]{0,0,0}$=$}%
}}}}
\end{picture}%
\end{center}
\end{enumerate}

\noindent Finally, we denote by $\equi{\Delta}$ the congruence relation on the free $2$-category $\mon{\Sigma^F}$ generated by the family~$E_{\Delta}$: this is the smallest equivalence relation on parallel $2$-arrows of $\mon{\Sigma^F}$ which contains the relations of~$E_{\Delta}$.
\end{ntt*}

\begin{rem*}
We recall the following definitions from [Guiraud 2004]. Let us assume that $R$ is a family of rewriting rules on parallel $2$-arrows generated by a $2$-polygraph $\Sigma$. If $\alpha:s_2(\alpha)\fl t_2(\alpha)$ is in $R$, then the \emph{reduction relation} $\red{\alpha}$ it generates is the smallest binary relation on parallel $2$-arrows of $\mon{\Sigma}$ which contains $\alpha$ and which is compatible with the two compositions of $\mon{\Sigma}$:
\begin{enumerate}
\item[-] We have $s_2(\alpha)\red{\alpha}t_2(\alpha)$.
\item[-] If $f\red{\alpha}g$ and if $h$ is a $2$-arrow of $\mon{\Sigma}$, then the following relations hold whenever their left (or right) side is defined:
$$
f\tens h\red{\alpha}g\tens h, \quad h\tens f\red{\alpha}h\tens g, \quad f\circ h\red{\alpha}g\circ h, \quad h\circ f\red{\alpha}h\circ g.
$$
\end{enumerate}
\noindent The \emph{reduction relation} $\red{R}$ generated by the whole of $R$ is the union of all the $\red{\alpha}$, for $\alpha$ in $R$. The relations $\mred{\alpha}$ and $\mred{R}$ are the reflexive-transitive closures of $\red{\alpha}$ and $\red{R}$. The relations $\equi{\alpha}$ and $\equi{R}$ are the reflexive-symmetric-transitive closures of $\red{\alpha}$ and $\red{R}$.
\end{rem*}

\begin{thm*}[Burroni]\label{thm:burroni}
The $2$-category $\Tb$ is isomorphic to the quotient $2$-category $\mon{\Sigma^F}/\equi{\Delta}$.
\end{thm*}

\begin{rem*}
The proof of theorem \ref{thm:burroni} is detailed in [Burroni 1993] in the one-sorted case and, as noted there, generalizes to the many-sorted case. It consists in the following steps:

\begin{enumerate}
\item[1.] One defines a $2$-functor $\pi$ from $\mon{\Sigma^F}$ to $\Tb$ as the unique $2$-functor such that:
\begin{enumerate}
\item[-] $\pi(\delta_A)=(a_1,a_1)$ and $\pi(\delta_F)=(x_1,x_1)$, respectively seen as a $2$-arrows from $A$ to $A^2$ and from $F$ to $F^2$.
\item[-] $\pi(\epsilon_A)=\ast(A)$ and $\pi(\epsilon_F)=\ast(F)$, where $\ast(A)$ (resp. $\ast(F)$) is the empty family of terms, seen as a $2$-arrow from $A$ (resp. $F$) to $\ast$, the empty family of wires.
\item[-] $\pi(\tau_{A,A})=(a_2,a_1)$, $\pi(\tau_{A,F})=(x_1,a_1)$, $\pi(\tau_{F,A})=(a_1,x_1)$ and $\pi(\tau_{F,F})=(x_2,x_1)$, seen respectively as $2$-arrows from $A^2$ to $A^2$, $A\tens F$ to $F\tens A$, $F\tens A$ to $A\tens F$ and $F^2$ to $F^2$.
\item[-] $\pi(\wedge)=x_1\wedge x_2$ and $\pi(\vee)=x_1\vee x_2$, both seen as $2$-arrows from $F^2$ to $F$.
\item[-] $\pi(\top)=\top$ and $\pi(\bot)=\bot$, both seen as $2$-arrows from $\ast$ to $F$.
\item[-] $\pi(\iota)=\iota(a_1)$, seen as a $2$-arrow from $A$ to $F$.
\item[-] $\pi(\nu)=\nu(a_1)$, seen as a $2$-arrow from $A$ to itself.
\end{enumerate}

\item[2.] Then one proves that $\pi$ is compatible with the relations of $E_{\Delta}$. This means that, for every $f\equiv g$ in~$E_{\Delta}$, $\pi(f)=\pi(g)$ holds. For example, let us prove this equality for the first relation, with the wires colored with $A$:
\begin{align*}
\pi\big((\delta_A\tens A)\circ\delta_A\big) \:
&=\: \big(\pi(\delta_A)\tens\pi(A)\big)\circ\pi(\delta_A) \\
&=\: \big((a_1,a_1)\tens a_1\big)\circ(a_1,a_1) \\
&=\: (a_1,a_1,a_2)\circ(a_1,a_1) \\
&=\: (a_1,a_1,a_1) \\
&=\: (a_1,a_2,a_2)\circ(a_1,a_1) \\
&=\: \big(a_1\tens(a_1,a_1)\big)\circ(a_1,a_1) \\
&=\: \big(\pi(A)\tens\pi(\delta_A)\big)\circ\pi(\delta_A) \\
&=\: \pi\big((A\tens\delta_A)\circ\delta_A\big).
\end{align*}

\item[3.] This proves that $\pi$ yields a $2$-functor from $\mon{\Sigma^F}/E_{\Delta}$ to $\Tb$. In order to prove that $\pi$ has an inverse, one starts with the construction of a decomposition of every $2$-arrow of $\Tb$ in elementary $2$-cells, all of the form $\pi(\phi)$, where $\phi$ is any $2$-cell of $\Sigma$. Let us consider a family $u=(u_1,\dots,u_n)$ of terms, seen as an arrow from $X$ to $Y$.

\begin{enumerate}

\item[-] The first layer is built only from the six operators of the terms signature $\Sr$, as the juxtaposition of the tree-parts of the terms $u_1$, $\dots$, $u_n$. For example, if $n=2$, $u_1=x_1\wedge\iota(\nu(a_2))$ and $u_2=x_1\vee\bot$, one gets:
\begin{center}
\includegraphics{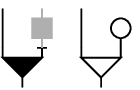}
\end{center}

\item[-] Then, the second layer is built from the eight operators of resources management. One takes the concatenation of the variables remaining from the first layer: in our example, $(x_1,a_2)$ and~$x_1$ remain, giving the family $(x_1,a_2,x_1)$. Then, one makes a diagram, using the resources management operators to link this family to the one corresponding to $X$. In our example, the following possibilities exist, among others, when $X=F^2\tens A^2$:
\begin{center}
\includegraphics{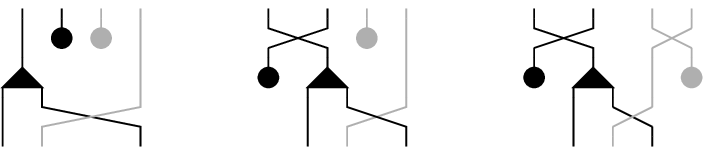}
\end{center}

\item[-] The seeked decomposition is built from the first layer, composed with any possible second layer on its top. In our example, we can get the following decompositions (note that we will make sure that the chosen one is the former):
\begin{center}
\includegraphics{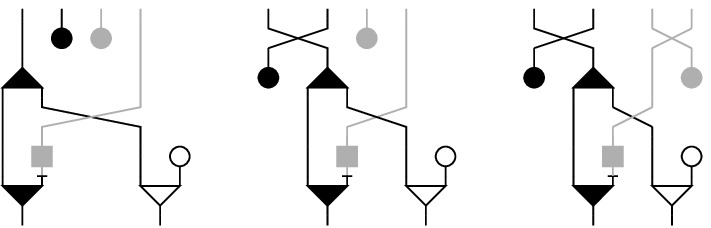}
\end{center}
\end{enumerate}

\item[4.] The final and most difficult part, fully detailed in [Burroni 1993], consists in proving that two decompositions of the same $2$-arrow are equal \emph{modulo} the relation $\equi{\Delta}$. This result comes from a polygraphic presentation of the $2$-category of finite sets. This step yields a $2$-functor from $\Tb$ to~$\mon{\Sigma^F}/E_{\Delta}$, that is checked to be inverse to $\pi$.
\end{enumerate}
\end{rem*}

\begin{rem*}
The family $E_{\Delta}$ of relations is minimal [Lafont 2003]: there is no other family with less elements that generates $\equi{\Delta}$. The result from [Burroni 1993] adapts to a general case, where the formal system to be translated into circuits is made of $n$ sorts and $m$ constructors: in this situtation, the first family would consist of $n(n^2+3n+3)$ relations, while the second one would have $m(n+2)$ relations.
\end{rem*}

\noindent For the moment, we have a translation from circuits into families of terms. In order to build translations going the reverse way, we prove that the family $E_{\Delta}$ can be extended into a finite, equivalent and convergent family of rewriting rules. The rules were given in [Lafont 2003] then proved to be convergent in [Guiraud 2004].

\begin{ntt*}
We denote by $R_{\Delta}$ the union of the two following families of rewriting rules on the $2$-category $\mon{\Sigma^F}$:
\begin{enumerate}
\item[1.] The first family consists of 42 rules, given by the following twelve schemes, with every possible colorations of wires:
\begin{center}
\includegraphics{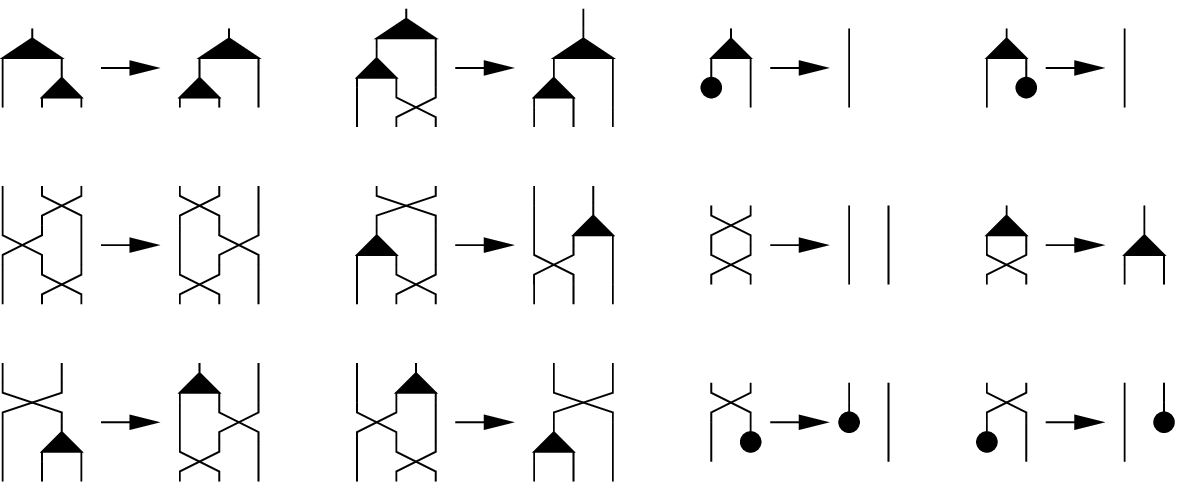}
\end{center}

\vfill\pagebreak
\item[2.] The second family consists of 36 rules, given by the following twelve schemes, sorted by arity of each constructor of $\Sr$, with every possible coloration of wires:
\begin{center}
\includegraphics{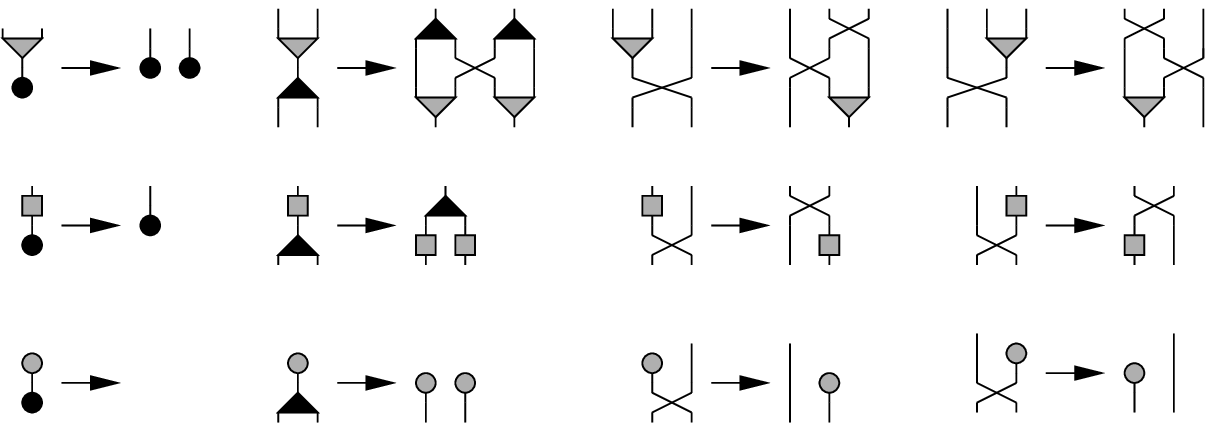}
\end{center}
\end{enumerate}
\end{ntt*}

\begin{rem*}
This definition extends to the case with $n$ sorts and $m$ constructors: the associated $2$-polygraph would have one cell in dimension $0$, $n$ cells in dimension $1$ and $m+2n+n^2$ cells in dimension~$2$. The set of rewriting rules on the $2$-polygraph would consist of $n(n^2+6n+5)$ rules in the first family and of $2m(n+1)$ rules in the second one.
\end{rem*}

\begin{lem*}
The families $R_{\Delta}$ and $E_{\Delta}$ are equivalent on $\mon{\Sigma^F}$.
\end{lem*}

\begin{dem}
We want to prove that the two families generate the same congruence relation on $\mon{\Sigma^F}$. Since~$E_{\Delta}$ is a subfamily of $R_{\Delta}$, it is sufficient to prove that each extra rule of $R_{\Delta}$ is derivable from $E_{\Delta}$. This means that, for each extra rule $f\fl g$, the relation $f\equi{\Delta}g$ holds. Let us consider, for example, the second scheme colored with $A$:
\begin{align*}
(A\tens\tau_{AA})\circ(\delta_A\tens A)\circ\delta_A
&\equi{\Delta} (A\tens\tau_{AA})\circ(A\tens\delta_A)\circ\delta_A \\
&= \big(A\tens(\tau_{AA}\circ\delta_A)\big)\circ\delta_A \\
&\equi{\Delta} (A\tens\delta_A)\circ\delta_A \\
&\equi{\Delta} (\delta_A\tens A)\circ\delta_A.
\end{align*}

\noindent We can also prove this fact graphically:
\begin{center}
\begin{picture}(0,0)%
\includegraphics{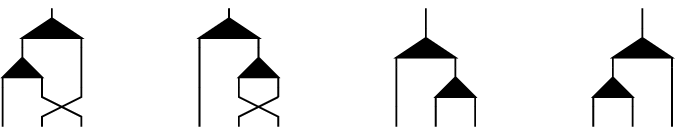}%
\end{picture}%
\setlength{\unitlength}{4144sp}%
\begingroup\makeatletter\ifx\SetFigFont\undefined%
\gdef\SetFigFont#1#2#3#4#5{%
  \reset@font\fontsize{#1}{#2pt}%
  \fontfamily{#3}\fontseries{#4}\fontshape{#5}%
  \selectfont}%
\fi\endgroup%
\begin{picture}(3084,564)(79,197)
\put(2521,389){\makebox(0,0)[b]{\smash{{\SetFigFont{10}{12.0}{\rmdefault}{\mddefault}{\updefault}{\color[rgb]{0,0,0}$\equi{\Delta}$}%
}}}}
\put(721,389){\makebox(0,0)[b]{\smash{{\SetFigFont{10}{12.0}{\rmdefault}{\mddefault}{\updefault}{\color[rgb]{0,0,0}$\equi{\Delta}$}%
}}}}
\put(1621,389){\makebox(0,0)[b]{\smash{{\SetFigFont{10}{12.0}{\rmdefault}{\mddefault}{\updefault}{\color[rgb]{0,0,0}$\equi{\Delta}$}%
}}}}
\end{picture}%
\end{center}

\noindent Let us make another graphical proof:
\begin{center}
\begin{picture}(0,0)%
\includegraphics{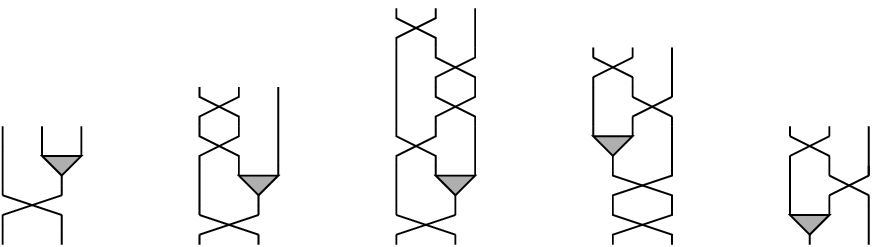}%
\end{picture}%
\setlength{\unitlength}{4144sp}%
\begingroup\makeatletter\ifx\SetFigFont\undefined%
\gdef\SetFigFont#1#2#3#4#5{%
  \reset@font\fontsize{#1}{#2pt}%
  \fontfamily{#3}\fontseries{#4}\fontshape{#5}%
  \selectfont}%
\fi\endgroup%
\begin{picture}(3984,1104)(79,197)
\put(3421,434){\makebox(0,0)[b]{\smash{{\SetFigFont{10}{12.0}{\rmdefault}{\mddefault}{\updefault}{\color[rgb]{0,0,0}$\equi{\Delta}$}%
}}}}
\put(721,434){\makebox(0,0)[b]{\smash{{\SetFigFont{10}{12.0}{\rmdefault}{\mddefault}{\updefault}{\color[rgb]{0,0,0}$\equi{\Delta}$}%
}}}}
\put(1621,434){\makebox(0,0)[b]{\smash{{\SetFigFont{10}{12.0}{\rmdefault}{\mddefault}{\updefault}{\color[rgb]{0,0,0}$\equi{\Delta}$}%
}}}}
\put(2521,434){\makebox(0,0)[b]{\smash{{\SetFigFont{10}{12.0}{\rmdefault}{\mddefault}{\updefault}{\color[rgb]{0,0,0}$\equi{\Delta}$}%
}}}}
\end{picture}%
\end{center}

\findem\end{dem}

\medskip
\noindent Now, we want to prove that the family $R_{\Delta}$ is convergent on $\mon{\Sigma^F}$. In [Guiraud 2004], the same set of rules was proved to be convergent on a monochromatic version of $\mon{\Sigma^F}$, which we denote here by $\mon{\Omega^F}$. Here, instead of doing the proof again, we can use this result to prove the convergence of $R_{\Delta}$. However, we need some extra notations.

In the $2$-category $\mon{\Omega^F}$, there is only one $1$-cell, denoted by $1$, and nine $2$-cells: the six from $\Sr$, with only their arity kept, together with $\epsilon$, $\delta$ and $\tau$. We define the $2$-functor $\gamma$ from $\mon{\Sigma^F}$ to $\mon{\Omega^F}$ as the only one which sends $A$ and $F$ onto $1$, each constructor of $\Sr$ onto itself, $\epsilon_X$ onto $\epsilon$, $\delta_X$ onto $\delta$ and $\tau_{XY}$ onto $\tau$.

Then let us consider a $2$-arrow $f$ in $\mon{\Omega^F}$ and a $1$-arrow $X$ in $\mon{\Sigma^F}$ such that $\gamma(X)=s_1(f)$ and such that there exists a $2$-arrow $f_X$ in $\mon{\Delta}$ with $1$-source $X$ and with $\gamma(f_X)=f$. In that case, $X$ is seen as a coloration of the input wires of the colorless $f$; then, this coloration is propagated throughout $f$, giving a label to each of $\epsilon$, $\delta$ and $\tau$ and yielding an arrow $f_X$.

Let us define these notions formally by induction on the \emph{size} of $2$-arrows: this is the least number of compositions $\circ$ and $\tens$ required to build them from the generators (the cells in each dimension); this notion is well defined because $\mon{\Sigma^F}$ is free.

\begin{defn*}
From now on, if $X$ is a $1$-arrow in any $2$-category, its identity $2$-arrow is also denoted by $X$. The set $\Gamma(f)$ of \emph{admissible colorations} for a $2$-arrow $f$ and the $2$-arrow $f_X$ are inductively defined as follows:
\begin{enumerate}
\item[-] If $n\in\Nb$, then $\Gamma(n)=\ens{A,F}^n$. If $X\in\ens{A,F}^n$, then $n_X=X$.
\item[-] If $\phi$ is in $\Sr$ with $1$-source $X$, then $\Gamma(\phi)=\ens{X}$ and $\phi_X=\phi$.
\item[-] The $2$-cells $\delta$ and $\epsilon$ satisfy $\Gamma(\delta)=\Gamma(\epsilon)=\ens{A,F}$. If $X\in\ens{A,F}$, then $\delta_X$ and $\epsilon_X$ are the $2$-cells of~$\Sigma$ with the same notations.
\item[-] The $2$-cell $\tau$ satisfies $\Gamma(\tau)=\ens{A,F}^2$. If $X,Y\in\ens{A,F}$, then $\tau_{X\tens Y}=\tau_{XY}$.
\item[-] If $f$ and $g$ are $2$-arrows of $\mon{\Omega^F}$, then $\Gamma(f\tens g)=\Gamma(f)\tens\Gamma(g)$. If $X$ is in $\Gamma(f)$ and $Y$ is in $\Gamma(g)$, then~$(f\tens g)_{X\tens Y}=f_X\tens g_Y$.
\item[-] If $f$ and $g$ are $2$-arrows of $\mon{\Omega^F}$ such that $t_1(f)=s_1(g)$, then $\Gamma(g\circ f)=\Gamma(f)$. If $X$ is an element of~$\Gamma(f)$, then~$(g\circ f)_X=g_{t_1(f_X)}\circ f_X$.
\end{enumerate}
\end{defn*}

\noindent By induction on the size of $f$, we get:

\begin{lem*}
For every $2$-arrow $f$ in $\mon{\Sigma^F}$, the $1$-arrow $s_1(f)$ is in $\Gamma(\gamma(f))$ and $(\gamma(f))_{s_1(f)}=f$.
\end{lem*}

\noindent We extend the constructions $\gamma$ and $(\cdot)_X$ on rules of $R_{\Delta}$ this way: for each rule $\alpha:f\fl g$ in $R_{\Delta}$, we denote by $\gamma(\alpha)$ the rule $\gamma(f)\fl\gamma(g)$ on $\mon{\Omega^F}$ and by $\gamma(R_{\Delta})$ the family of all rules $\gamma(\alpha)$.

Conversely, each rule $\alpha:f\fl g$ in $\gamma(R_{\Delta})$ yields one, two, four or eight rules in $R_{\Delta}$, each one of the form $\alpha_X:f_X\fl g_X$.

We prove the following result, using the definition of $\red{\alpha}$ and the functoriality of $\gamma$:

\begin{lem*}
For every rule $\alpha$ in $R_{\Delta}$ and every $2$-arrows $f$ and $g$ in $\mon{\Sigma^F}$ such that $f\red{\alpha}g$, then the property $\gamma(f)\red{\gamma(\alpha)}\gamma(g)$ holds in $\mon{\Omega^F}$.

Conversely, if $f$ and $g$ are $2$-arrows in $\mon{\Omega^F}$ and $\alpha$ is a rule in $\gamma(R_{\Delta})$ such that $f\red{\alpha}g$, then every $1$-arrow~$X$ in $\Gamma(f)$ is in $\Gamma(g)$ and $f_X\red{\alpha_Y}g_X$ holds for some $Y$ in~$\mon{\Sigma^F}$.
\end{lem*}

\noindent Then, we can prove:

\begin{thm*}\label{th:convergence-gestion-ressources}
The family of rules $R_{\Delta}$ is convergent on $\mon{\Sigma^F}$.
\end{thm*}

\begin{dem}
In order to prove the termination, let us assume that there exists an infinite reduction path $(f_n)_{n\in\Nb}$ in $\mon{\Sigma^F}$ generated by $R_{\Delta}$: this means that the $f_n$ are parallel $2$-arrows of $\mon{\Sigma^F}$ such that, for every $n$, there exists a rule $\alpha_n$ in $R_{\Delta}$ with $f_n\red{\alpha_n}f_{n+1}$. From the previous lemma, we deduce that, for every $n$, the reduction $\gamma(f_n)\red{\gamma(\alpha_n)}\gamma(f_{n+1})$ holds in $\mon{\Omega^F}$. Hence $(\gamma(f_n))_n$ is an infinite reduction path in $\mon{\Omega^F}$ generated by~$\gamma(R_{\Delta})$. However, we know since [Guiraud 2004] that $\gamma(R_{\Delta})$ terminates on $\mon{\Omega^F}$: this prevents the existence of such an infinite path. Hence $R_{\Delta}$ terminates on $\mon{\Sigma^F}$.

Now, let us consider a branching $(f,g,h)$ generated by $R_{\Delta}$ in $\mon{\Sigma^F}$: this means that $f$, $g$ and $h$ are parallel $2$-arrows such that there exist two reduction paths in $\mon{\Sigma^F}$ with the following shapes:
$$
f\red{\alpha_1}g_1\red{\alpha_2}\dots\red{\alpha_m} g \et f\red{\beta_1}h_1\red{\beta_2}\dots\red{\beta_n} h,
$$

\noindent with all the $\alpha_i$ and $\beta_j$ in $R_{\Delta}$. Then, an application of $\gamma$ on both paths proves that the triple of $2$-arrows $(\gamma(f),\gamma(g),\gamma(h))$ in $\mon{\Omega^F}$ is a branching generated by $\gamma(R_{\Delta})$. Indeed, from the previous lemma, we get:
$$
\gamma(f)\red{\gamma(\alpha_1)}\gamma(g_1)\red{\gamma(\alpha_2)}\dots\red{\gamma(\alpha_m)} \gamma(g) \et \gamma(f)\red{\gamma(\beta_1)}\gamma(h_1)\red{\gamma(\beta_2)}\dots\red{\gamma(\beta_n)} \gamma(h).
$$

\noindent We know that $\gamma(R_{\Delta})$ is confluent from [Guiraud 2004]. Hence, the branching $(\gamma(f),\gamma(g),\gamma(h))$ can be closed with a $2$-arrow $k$ in $\mon{\Omega^F}$, together with two reduction paths generated by $\gamma(R_{\Delta})$:
$$
\gamma(g)\red{\alpha'_1}g'_1\red{\alpha'_2}\dots\red{\alpha'_p} k \et \gamma(h)\red{\beta'_1}h'_1\red{\beta'_2}\dots\red{\beta'_q} k.
$$

\noindent Let us consider $X=s_1(f)$ in $\mon{\Sigma^F}$. Since all considered arrows in these paths are parallel, $X$ is the $1$-source of all of them and, in particular, admissible for all of them. Then, an application of $(\cdot)_X$ yields $1$-arrows denoted by $Y_1$, $\dots$, $Y_p$ and $Z_1$, $\dots$, $Z_q$ such that:
$$
g \red{(\alpha'_1)_{Y_1}}(g'_1)_X\red{(\alpha'_2)_{Y_2}}\dots\red{(\alpha'_p)_{Y_p}} k_X \et h\red{(\beta'_1)_{Z_1}}(h'_1)_X\red{(\beta'_2)_{Z_2}}\dots\red{(\beta'_q)_{Z_q}} k_X.
$$

\noindent Hence, there exist reduction paths generated by $R_{\Delta}$ from $g$ to $k_X$ and from $h$ to $k_X$, so that $k_X$ closes the branching $(f,g,h)$. Thus $R_{\Delta}$ is also confluent and, finally, convergent.
\findem\end{dem}

\begin{ntt*}
If $f$ is a $2$-arrow in $\mon{\Sigma^F}$, we denote by $R_{\Delta}(f)$ its unique normal form with respect to the congruence~$\equi{\Delta}$.
\end{ntt*}

\noindent Now, we can define translations from families of terms to circuits:

\begin{ntt*}
Let $u=(u_1,\dots,u_n)$ be a family of terms, let $X$ be a $1$-arrow of $\mon{\Sigma^F}$ such that $\sharp u\leq\sharp X$ and let $Y=Y_1\tens\dots\tens Y_n$ be the $1$-arrow such that $u_i\in Y_i$. Then we denote by $\Phi^X(u)$ the unique $2$-arrow from $X$ to $Y$ in $\mon{\Sigma^F}$ which is in normal form with respect to $R_{\Delta}$ and such that $\pi(\Phi^X(u))=u$. In the special case where $X$ is of the form $A^k\tens F^l$, $\Phi^X(u)$ is denoted by $\Phi^{(k,l)}(u)$.
\end{ntt*}

\noindent In order to conclude this technical part, we prove the following result:

\begin{lem*}
For every $2$-arrow $f$ in $\mon{\Sigma^F}$, we have $R_{\Delta}(f)=\Phi^{s_1(f)}(\pi(f))$.
\end{lem*}

\begin{dem}
\noindent By definition of $\Phi^{s_1(f)}(\pi(f))$, the following relations are satisfied:
$$
\begin{aligned}
s_1(\Phi^{s_1(f)}(\pi(f))) &\:=\: s_1(f), \\
t_1(\Phi^{s_1(f)}(\pi(f))) &\:=\: t_1(f), \\
\pi(\Phi^{s_1(f)}(\pi(f))) &\:=\: \pi(f), \\
R_{\Delta}(\Phi^{s_1(f)}(\pi(f))) &\:=\: \Phi^{s_1(f)}(\pi(f))).
\end{aligned}
$$

\noindent The first two equations tell us that $\Phi^{s_1(f)}(\pi(f))$ and $f$ are parallel $2$-arrows of $\mon{\Sigma^F}$. The third one gives that both $f$ and $\Phi^{s_1(f)}(\pi(f))$ have the same image through $\pi$. However, we already know that, for every parallel $2$-arrows $g$ and $h$ in $\mon{\Sigma^F}$, we have $g\equi{\Delta}h$ if and only if $\pi(g)$ and~$\pi(h)$ are equal. Thus, $f\equi{\Delta}\Phi^{s_1(f)}(\pi(f))$, which is equivalent, since $R_{\Delta}$ is a convergent presentation of $\equi{\Delta}$, to the fact that $R_{\Delta}(f)$ and $R_{\Delta}(\Phi^{s_1(f)}(\pi(f)))$ are equal. Finally, the fourth equation gives the result: $R_{\Delta}(f)=\Phi^{s_1(f)}(\pi(f)))$.

\findem\end{dem}

\subsection{Translation of the structural congruence}\label{sub:congruence-2d}

\noindent In this paragraph, we give one way to translate the relation of structural congruence from terms to circuits. We generalize a result from [Guiraud 2004], from the one-sorted to the two-sorted case. For that, we define translations from rewriting rules on terms to rewriting rules on circuits:

\begin{ntt*}
Let $\alpha:u\fl v$ be a rewriting rule on terms. We denote by $\sharp\alpha$ the pair $(\sharp_A\alpha,\sharp_F\alpha)$ of natural numbers that is the upper bound of $\sharp u$ and $\sharp v$. Then $\Phi(\alpha)$ is defined as the rewriting rule $\Phi^{\sharp\alpha}(u)\fl\Phi^{\sharp\alpha}(v)$ on $\mon{\Sigma^F}$. If $R$ is a family of rules on terms, then $\Phi(R)$ is the family made of the translations through $\Phi$ of each rule in $R$.
\end{ntt*}

\begin{rem*}
The definition of $\Phi(\alpha)$ is not restricted to rewriting rules: indeed, the left part can be a variable and the right part may contain more variables than the left one. This would create infinite reduction paths, but in what follows we are not really interested in the rewriting properties of paths, but rather in their classification.
\end{rem*}

\noindent Let us prove that redexes are preserved by the translations from terms to circuits.

\begin{lem*}
Let $u$ be a term, $C$ be a context and $\sigma$ be a substitution. Let $X$ be a $1$-arrow of $\mon{\Sigma^F}$ such that $\sharp X$ is greater than $\sharp u$. Then, there exists a $2$-arrow $f$ in $\mon{\Sigma^F}$ such that:
\begin{enumerate}
\item[-] There exist $1$-arrows $Y$ and $Z$ and $2$-arrows $h$ and $k$ in $\mon{\Sigma^F}$ with:
$$
f\:=\:k\circ(Y\tens\Phi^X(u)\tens Z)\circ h.
$$

\item[-] The relation $\pi(f)=C[u\cdot\sigma]$ holds.
\end{enumerate}
\end{lem*}

\begin{dem}
Let us denote by $(x_1,\dots,x_n)$ the family of variables $\pi(X)$. Let us denote by $(y_1,\dots,y_k)$ and~$(z_1,\dots,z_l)$ the two families of variables appearing from left to right in the context $C$, the first one at the left of the empty slot, the second one at its right. Let us consider any $1$-arrow $X'$ such that $\sharp X'\geq\sharp C[u\cdot\sigma]$ holds. Then, we denote by $h$ the arrow:
$$
h=\Phi^{X'}(y_1,\dots,y_k,x_1\cdot\sigma,\dots,x_n\cdot\sigma,z_1,\dots,z_l).
$$

\noindent Then, let us consider the term $C_0$ built from the tree-part of the context $C$ by putting variables on each leaf, with no repetition and in order from left to right. Let us denote by $U$ the sort of the term $u$, which is either $A$ or $F$. Let us denote by $Y_1$, $\dots$, $Y_k$ the sorts of the variables $y_1$, $\dots$, $y_k$ and by $Z_1$, $\dots$, $Z_l$ the sorts of the variables $z_1$, $\dots$, $z_l$. Finally, $Y$ is the product $Y_1\tens\dots\tens Y_k$ and $Z$ is the product $Z_1\tens\dots\tens Z_l$. We define $k$ as the arrow:
$$
k=\Phi^{Y\tens U\tens Z}(C_0).
$$

\noindent Then the $2$-arrow $f=k\circ(Y\tens\Phi^X(u)\tens Z)\circ h$ has been built to satisfy $\pi(f)=C[u\cdot\sigma]$.
\findem\end{dem}

\medskip
\noindent Now we can prove that reductions on terms can be lifted to reductions on the corresponding circuits.

\begin{prop*}
Let $\alpha$ be a rewriting rule on the set of terms. If $u$ and $v$ are terms such that $u\red{\alpha}v$, then for every $1$-arrow $X$ such that $\sharp X$ is greater than both $\sharp u$ and $\sharp v$, there exist $2$-arrows $f$ and $g$ in~$\mon{\Sigma^F}$ such that:
$$
\Phi^X(u)\equi{\Delta} f\red{\Phi(\alpha)}g\equi{\Delta}\Phi^X(v).
$$
\end{prop*}

\begin{dem}
Let us use the notations $\alpha:s(\alpha)\fl t(\alpha)$ and $\Phi(\alpha):s_2(\Phi(\alpha))\fl t_2(\Phi(\alpha))$. Since $u\red{\alpha}v$, there exist a context $C$ and a substitution $\sigma$ such that $u=C[s(\alpha)\cdot\sigma]$ and $v=C[t(\alpha)\cdot\sigma]$. From the previous lemma, this implies that there exist $2$-arrows $h$ and $k$ and $1$-arrows $Y$ and $Z$ in $\mon{\Sigma^F}$ such that the $2$-arrows defined thereafter satisfy $\pi(f)=u$ and $\pi(g)=v$:
$$
f=k\circ(Y\tens s_2(\Phi(\alpha))\tens Z)\circ h \et g=k\circ(Y\tens t_2(\Phi(\alpha)\tens Z)\circ h.
$$

\noindent Hence $f\red{\Phi(\alpha)}g$. Furthermore, since $\pi(f)=u=\pi(\Phi^X(u))$, we know that $\Phi^X(u)\equi{\Delta} f$ and, for the same reasons, $\Phi^X(v)\equi{\Delta}g$, which concludes the proof.
\findem\end{dem}

\begin{cor*}
Let $R$ be a family of relations or rewriting rules on terms, let $\equi{R}$ be the congruence it generates on terms and $\equi{\Delta R}$ the one on parallel circuits generated by the union of $R_{\Delta}$ and $\Phi(R)$. If $u$ and~$v$ are two terms such that $u\equi{R} v$, then $\Phi^X(u)\equi{\Delta R}\Phi^X(v)$ holds for every $X$ such that $\sharp X$ is greater than both $\sharp u$ and $\sharp v$. Conversely, if $f$ and $g$ are two $2$-arrows of $\mon{\Sigma^F}$ such that $f\equi{\Delta R}g$, then~$\pi(f)\equi{R}\pi(g)$.
\end{cor*}

\noindent We use this result on the example of the structural rules:

\begin{defn*}
The family $S$ is the following family of rules:
\begin{center}
\includegraphics{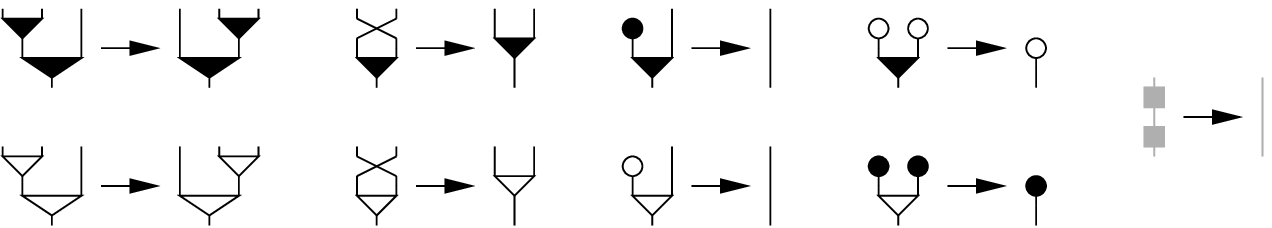}
\end{center}
\end{defn*}

\begin{rem*}
The rules for commutativity have been reversed, only for aesthetic and termination reasons. This choice does not change the congruence they generate on circuits. Furthermore, the rule $x_1\wedge x_2\fl x_2\wedge x_1$ and its converse generate the same reduction relation on terms, hence the same congruence.
\end{rem*}

\begin{rem*}
Informally, the set of circuits, equipped with the reduction relation $\red{S}$, is projected, through $\pi$, onto the set of families of terms, equipped with $\red{S}$, and the fiber of $\pi$ over each family $u$ of terms is an $\equi{\Delta}$-equivalence. One of the future objectives of higher-dimensional rewriting is to make this remark formal.
\end{rem*}

\begin{rem*}
The given set of structural rules is terminating but not confluent. However, it can be completed into a convergent one. One open question is to determine if the union of the resources management rules and of the structural rules can be completed into an equivalent, finite and convergent set of rules. The main direction towards this result consists in following the example of the rewriting system~$\mathrm{L}(\Zb_2)$, proposed in [Lafont 2003] as a finite presentation of the structure of $\Zb/2\Zb$-vector spaces, and proved to be convergent in [Guiraud 2004]. Such a result would provide canonical representatives of formulas (\emph{modulo} structural congruence) into circuits.
\end{rem*}

\noindent So far, we have translated the structural congruence from terms to circuits in such a way that, for any parallel $2$-arrows $f$ and $g$ in $\mon{\Sigma^F}$, we have $f\equi{\Delta S} g$ if and only if $\pi(f)\equi{S}\pi(g)$. However, all the relations between $2$-arrows can be given a name and a richer structure than a mere relational one: they have an intrinsic $3$-dimensional nature and so have the inference rules generating the proofs.

\section{The three dimensions of proofs}\label{sec:preuves}

\noindent After a presentation of the inference rules of SKS (\ref{sub:preuves}), we give the intuition leading to the construction we seek (\ref{sub:idee-preuves-3d}). Once again, this requires some theoretical notions (\ref{sub:3-poly}). Then we give the formal translation (\ref{sub:preuves-3d}) and prove that the $3$-dimensional object one gets can be equipped with a notion of proof that corresponds to the one of SKS (theorem~\ref{thm:1}).

\subsection{The SKS proofs}\label{sub:preuves}

\noindent In this paragraph, we recall definitions from [Brünnler 2003]. Once again, they are slightly adapted to our needs; in particular they are written in a term rewriting style.

\begin{defn*}
The \emph{SKS inference rules} are the following rewriting rules on the set $T$ of SKS terms:
$$
\begin{array}{r c l c r c l}
\top &\lfl& \iota(a_1)\vee\iota(\nu(a_1)) && \iota(a_1)\wedge\iota(\nu(a_1)) &\lfl& \bot \\
&&\hfill (x_1\vee x_2)\wedge x_3 &\lfl& x_1\vee (x_2\wedge x_3) \hfill\strut \\
&&\hfill (x_1\wedge x_2)\vee(x_3\wedge x_4) &\lfl& (x_1\vee x_3)\wedge(x_2\vee x_4) \hfill\strut \\
\bot &\lfl& \iota(a_1) && \iota(a_1) &\lfl& \top \\
\iota(a_1)\vee \iota(a_1) &\lfl& \iota(a_1) && \iota(a_1) &\lfl& \iota(a_1)\wedge \iota(a_1).
\end{array}
$$

\noindent The set of SKS inference rules is denoted by $R$. We denote by $S$ the set of structural rules on SKS terms, by $S^{-1}$ the same set with the rules reversed and by $\ol{S}$ the union of both sets.
\end{defn*}

\noindent Note that the generated congruences $\equi{S}$, $\equi{S^{-1}}$ and $\equi{\ol{S}}$ are the same relations. We define a graphical object associated to SKS in which arrows are formal proofs.

\vfill\pagebreak
\begin{defn*}
\noindent The \emph{reduction graph associated to SKS} is the graph $G^K$ defined as follows:
\begin{enumerate}
\item[0.] Its objects are the families of SKS terms.
\item[1.] If $u$ and $v$ are two objects, then there is an arrow in $G$ from $u$ to $v$ for each $\alpha$ in either of $R$ or $\ol{S}$ such that $u\red{\alpha}v$.
\end{enumerate}

\noindent A \emph{SKS proof from $u$ to $v$} is a finite path in the graph $G^K$, starting at $u$ and ending at $v$. A \emph{complete SKS proof of $u$} is a SKS proof from $\top$ to $u$.
\end{defn*}

\noindent Hence, the SKS proofs are the rewriting paths generated by the inference rules, together with the structural rules and their converse. In [Guiraud 2004], it was proved that any term rewriting system can be translated into a $3$-polygraph, alike what was done for structural rules in the previous section.

\subsection{From proofs to three-dimensional arrows: the informal idea}\label{sub:idee-preuves-3d}

The inference rules are rewriting rules on circuits: they transform one circuit into another one, with the same inputs and the same outputs. Let us consider a rewriting rule $\alpha:f\fl g$ on circuits and call $f$ the $2$-source and $g$ the $2$-target of $\alpha$. Then, the fact that $f$ and $g$ are parallel means that $f$ and $g$ have the same $1$-source and the same $1$-target. Equationally, $s_1(f)=s_1(g)$ and $t_1(f)=t_1(g)$. If we denote $f$ by $s_2(\alpha)$ and $g$ by $t_2(\alpha)$, then we get:
$$
s_1\circ s_2(\alpha) = s_1\circ t_2(\alpha) \et t_1\circ s_2(\alpha) = t_1\circ t_2(\alpha).
$$

\medskip
\noindent This means that $\alpha$ can be seen as a $3$-dimensional cell over the free $2$-category $\mon{\Sigma^F}$: a directed volume between two parallel directed surfaces. However, such an object is difficult to represent. For that reason, we use here another type of pictures, in order to give the intuition, where $3$-cells are drawn as blocks:
\begin{center}
\begin{picture}(0,0)%
\includegraphics{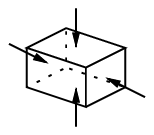}%
\end{picture}%
\setlength{\unitlength}{4144sp}%
\begingroup\makeatletter\ifx\SetFigFont\undefined%
\gdef\SetFigFont#1#2#3#4#5{%
  \reset@font\fontsize{#1}{#2pt}%
  \fontfamily{#3}\fontseries{#4}\fontshape{#5}%
  \selectfont}%
\fi\endgroup%
\begin{picture}(1123,873)(928,-61)
\put(1486,704){\makebox(0,0)[b]{\smash{{\SetFigFont{10}{12.0}{\rmdefault}{\mddefault}{\updefault}{\color[rgb]{0,0,0}$x$}%
}}}}
\put(1036,479){\makebox(0,0)[b]{\smash{{\SetFigFont{10}{12.0}{\rmdefault}{\mddefault}{\updefault}{\color[rgb]{0,0,0}$f$}%
}}}}
\put(1936,119){\makebox(0,0)[b]{\smash{{\SetFigFont{10}{12.0}{\rmdefault}{\mddefault}{\updefault}{\color[rgb]{0,0,0}$g$}%
}}}}
\put(1486,-16){\makebox(0,0)[b]{\smash{{\SetFigFont{10}{12.0}{\rmdefault}{\mddefault}{\updefault}{\color[rgb]{0,0,0}$y$}%
}}}}
\end{picture}%
\end{center}

\noindent This represents a $3$-cell $\alpha$ going from a circuit $f$ to another one $g$. Both circuits must have the same inputs (number and color), here $x$, and the same ouputs, here $y$. Note that, although useful, this representation can be misleading: for example, $x$ is a $1$-dimensional cell but it is pictured as a $2$-dimensional object, like $f$.

This being noticed, we use this block representation for intuition, together with the following one, much more accurate though only $2$-dimensional, made of three vertical slices of the block - one before, one in the middle, one after:
\begin{center}
\begin{picture}(0,0)%
\includegraphics{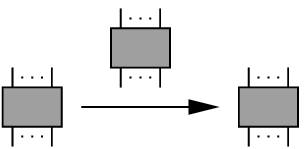}%
\end{picture}%
\setlength{\unitlength}{4144sp}%
\begingroup\makeatletter\ifx\SetFigFont\undefined%
\gdef\SetFigFont#1#2#3#4#5{%
  \reset@font\fontsize{#1}{#2pt}%
  \fontfamily{#3}\fontseries{#4}\fontshape{#5}%
  \selectfont}%
\fi\endgroup%
\begin{picture}(1374,654)(79,377)
\put(721,794){\makebox(0,0)[b]{\smash{{\SetFigFont{10}{12.0}{\rmdefault}{\mddefault}{\updefault}{\color[rgb]{0,0,0}$\alpha$}%
}}}}
\put(1306,524){\makebox(0,0)[b]{\smash{{\SetFigFont{10}{12.0}{\rmdefault}{\mddefault}{\updefault}{\color[rgb]{0,0,0}$g$}%
}}}}
\put(226,524){\makebox(0,0)[b]{\smash{{\SetFigFont{10}{12.0}{\rmdefault}{\mddefault}{\updefault}{\color[rgb]{0,0,0}$f$}%
}}}}
\end{picture}%
\end{center}

\noindent Hence, giving a set of rewriting rules on a free $2$-category amounts at giving a family of $3$-cells over it: this is a \emph{$\mathit{3}$-polygraph}. Furthermore, this object generates a reduction graph which paths will be proved to be representatives of the SKS proofs.

In order to give the underlying idea, let us consider extensions of the two compositions of circuits on $3$-cells: with these operations, one can put circuits aside a block or plug another ones in its inputs and outputs. Let us give an example, with the sliced representation:
\begin{center}
\begin{picture}(0,0)%
\includegraphics{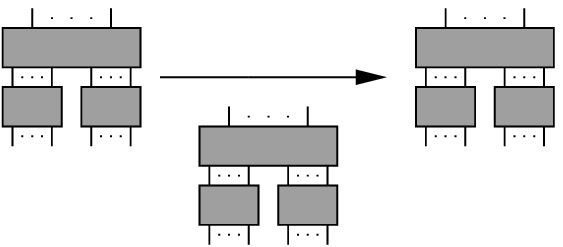}%
\end{picture}%
\setlength{\unitlength}{4144sp}%
\begingroup\makeatletter\ifx\SetFigFont\undefined%
\gdef\SetFigFont#1#2#3#4#5{%
  \reset@font\fontsize{#1}{#2pt}%
  \fontfamily{#3}\fontseries{#4}\fontshape{#5}%
  \selectfont}%
\fi\endgroup%
\begin{picture}(2544,1104)(79,-73)
\put(226,524){\makebox(0,0)[b]{\smash{{\SetFigFont{10}{12.0}{\rmdefault}{\mddefault}{\updefault}{\color[rgb]{0,0,0}$f$}%
}}}}
\put(586,524){\makebox(0,0)[b]{\smash{{\SetFigFont{10}{12.0}{\rmdefault}{\mddefault}{\updefault}{\color[rgb]{0,0,0}$h$}%
}}}}
\put(406,794){\makebox(0,0)[b]{\smash{{\SetFigFont{10}{12.0}{\rmdefault}{\mddefault}{\updefault}{\color[rgb]{0,0,0}$k$}%
}}}}
\put(2116,524){\makebox(0,0)[b]{\smash{{\SetFigFont{10}{12.0}{\rmdefault}{\mddefault}{\updefault}{\color[rgb]{0,0,0}$g$}%
}}}}
\put(2476,524){\makebox(0,0)[b]{\smash{{\SetFigFont{10}{12.0}{\rmdefault}{\mddefault}{\updefault}{\color[rgb]{0,0,0}$h$}%
}}}}
\put(2296,794){\makebox(0,0)[b]{\smash{{\SetFigFont{10}{12.0}{\rmdefault}{\mddefault}{\updefault}{\color[rgb]{0,0,0}$k$}%
}}}}
\put(1126, 74){\makebox(0,0)[b]{\smash{{\SetFigFont{10}{12.0}{\rmdefault}{\mddefault}{\updefault}{\color[rgb]{0,0,0}$\alpha$}%
}}}}
\put(1486, 74){\makebox(0,0)[b]{\smash{{\SetFigFont{10}{12.0}{\rmdefault}{\mddefault}{\updefault}{\color[rgb]{0,0,0}$h$}%
}}}}
\put(1306,344){\makebox(0,0)[b]{\smash{{\SetFigFont{10}{12.0}{\rmdefault}{\mddefault}{\updefault}{\color[rgb]{0,0,0}$k$}%
}}}}
\end{picture}%
\end{center}

\noindent In this diagram, we see an application of the rule $\alpha$ in the context formed of all the surrounding circuits: it transforms $(f\tens h)\circ k$ into $(g\tens h)\circ k$. This operation corresponds, \emph{modulo} some $\equi{\Delta}$ equivalences, to an application in context of a SKS rule.

If one considers the graph made of all applications (in context) of rules on circuits, its paths should have a strong link with the SKS proofs. This is what will be explored, after some formal definitions.

\subsection{Three-polygraphs and their reduction graphs}\label{sub:3-poly}

We start with the definition of $3$-cells over a $2$-category:

\begin{defn*}
Let $\Cr$ be a $2$-category. A \emph{family of $\mathit{3}$-cells over $\Cr$} is a triple $(\Sigma_3,s_2,t_2)$ made of a set~$\Sigma_3$ and two maps $s_2,t_2:\Sigma_3\fl\Cr_2$ such that the following two equations hold:
$$
s_1\circ s_2 = s_1\circ t_2 \et t_1\circ s_2 = t_1\circ t_2.
$$
\end{defn*}

\begin{ex*}
We have already encountered several families of $3$-cells, both over the $2$-category~$\mon{\Sigma^F}$: the resources management equations (seen as $3$-cells going from left to right), the resources management rules, the structural rules and their reverse rules. Furthermore, we have seen that any rule on terms generates a $3$-cell over $\mon{\Sigma^F}$.
\end{ex*}

\begin{defn*}
The \emph{family of inference rules} is the family of $3$-cells over the free $2$-category $\mon{\Sigma^F}$ given graphically as follows:
\begin{center}
\includegraphics{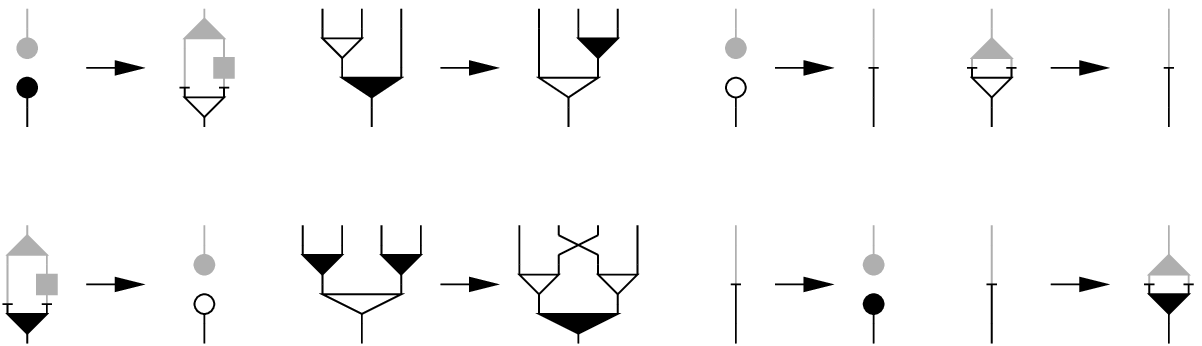}
\end{center}
\end{defn*}

\vfill\pagebreak
\noindent Such an extension of a $2$-polygraph is called a $3$-polygraph:

\begin{defn*}
A \emph{$\mathit{3}$-polygraph} is a data $(\Sigma,\Sigma_3,s_2,t_2)$ made of a $2$-polygraph $\Sigma$ and a family of $3$-cells $(\Sigma_3,s_2,t_2)$ over the free $2$-category generated by $\Sigma$. The elements of $\Sigma_3$ are the \emph{$\mathit{3}$-cells} of the $3$-polygraph and the maps $s_2$ and $t_2$ are respectively the \emph{$\mathit{2}$-source} and the \emph{$\mathit{2}$-target} maps.

A $3$-polygraph $\Sigma$ is often denoted by the family $(\Sigma_0,\Sigma_1,\Sigma_2,\Sigma_3)$ of its sets of $0$, $1$, $2$ and $3$-cells, assuming that the sources and targets are implicitely given with them.
\end{defn*}

\begin{ex*}\label{ex:familles-3-cellules}
The $2$-polygraph $\Sigma^F$ associated to the signature of SKS terms can be extended into a $3$-polygraph with any of the families encountered so far. For example, with the following ones:
\begin{enumerate}
\item[-] The family of 100 resources management $3$-cells, made from the equations of $E_{\Delta}$ and their converse (each equation is split into two rules, one going in one direction, one in the reverse direction).
\item[-] The family of 18 structural $3$-cells, made of the translation of structural rules and their converse.
\item[-] The family $R$ of 8 inference $3$-cells.
\end{enumerate}

\noindent All these $3$-polygraphs have the same cells in dimensions $0$, $1$, $2$. Hence, one can enrich $\Sigma^F$ with the union of any of all these families.
\end{ex*}

\begin{defn*}
The $3$-polygraph $\Sigma^K$ consists of the $2$-polygraph $\Sigma^F$ extended with the three families of $3$-cells from example \ref{ex:familles-3-cellules}: it has one cell in dimension $0$, two in dimension $1$, 14 in dimension $2$ and~126 in dimension $3$.
\end{defn*}

\noindent We define the reduction graph associated to a $3$-polygraph with only one $0$-cell: this is the case we need and this restriction makes graphical representations clearer. The idea behind this notion is that an arrow in this graph is an application of a $3$-cell, seen as a rewrite rule, inside a context. Note that [Guiraud 2004(T)] contains a formal categorical approach to contexts over a $2$-polygraph.

\begin{defn*}
Let $\Sigma=(\ast,\Sigma_1,\Sigma_2,\Sigma_3)$ be a $3$-polygraph with one $0$-cell. Its \emph{associated reduction graph} is the graph denoted by $G(\Sigma)$ defined this way:
\begin{enumerate}
\item[0.] The objects of $G(\Sigma)$ are the $2$-arrows of $\mon{\Sigma}_2$.
\item[1.] The arrows of $G(\Sigma)$ from $u$ to $v$ are all the triples $(f,\alpha,g)$, made of two $2$-arrows $f$ and $g$ of $\mon{\Sigma}_2$ and one $3$-cell $\alpha$ of $\Sigma_3$, such that the following two equalities are defined and hold:
\begin{center}
\begin{picture}(0,0)%
\includegraphics{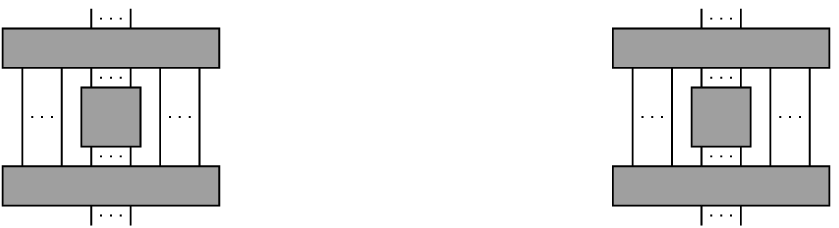}%
\end{picture}%
\setlength{\unitlength}{4144sp}%
\begingroup\makeatletter\ifx\SetFigFont\undefined%
\gdef\SetFigFont#1#2#3#4#5{%
  \reset@font\fontsize{#1}{#2pt}%
  \fontfamily{#3}\fontseries{#4}\fontshape{#5}%
  \selectfont}%
\fi\endgroup%
\begin{picture}(4447,1014)(79,-253)
\put(3376,209){\makebox(0,0)[b]{\smash{{\SetFigFont{8}{9.6}{\rmdefault}{\mddefault}{\updefault}{\color[rgb]{0,0,0}$t_2\alpha$}%
}}}}
\put(586,524){\makebox(0,0)[b]{\smash{{\SetFigFont{10}{12.0}{\rmdefault}{\mddefault}{\updefault}{\color[rgb]{0,0,0}$f$}%
}}}}
\put(3376,524){\makebox(0,0)[b]{\smash{{\SetFigFont{10}{12.0}{\rmdefault}{\mddefault}{\updefault}{\color[rgb]{0,0,0}$f$}%
}}}}
\put(586,-106){\makebox(0,0)[b]{\smash{{\SetFigFont{10}{12.0}{\rmdefault}{\mddefault}{\updefault}{\color[rgb]{0,0,0}$g$}%
}}}}
\put(3376,-106){\makebox(0,0)[b]{\smash{{\SetFigFont{10}{12.0}{\rmdefault}{\mddefault}{\updefault}{\color[rgb]{0,0,0}$g$}%
}}}}
\put(1351,209){\makebox(0,0)[b]{\smash{{\SetFigFont{10}{12.0}{\rmdefault}{\mddefault}{\updefault}{\color[rgb]{0,0,0}$=$}%
}}}}
\put(1621,209){\makebox(0,0)[b]{\smash{{\SetFigFont{10}{12.0}{\rmdefault}{\mddefault}{\updefault}{\color[rgb]{0,0,0}$u$}%
}}}}
\put(4141,209){\makebox(0,0)[b]{\smash{{\SetFigFont{10}{12.0}{\rmdefault}{\mddefault}{\updefault}{\color[rgb]{0,0,0}$=$}%
}}}}
\put(4411,209){\makebox(0,0)[b]{\smash{{\SetFigFont{10}{12.0}{\rmdefault}{\mddefault}{\updefault}{\color[rgb]{0,0,0}$v$}%
}}}}
\put(586,209){\makebox(0,0)[b]{\smash{{\SetFigFont{8}{9.6}{\rmdefault}{\mddefault}{\updefault}{\color[rgb]{0,0,0}$s_2\alpha$}%
}}}}
\end{picture}%
\end{center}

\noindent Each triple $(f,\alpha,g)$ is represented by the following diagram:
\begin{center}
\begin{picture}(0,0)%
\includegraphics{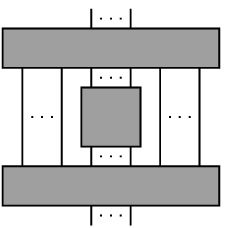}%
\end{picture}%
\setlength{\unitlength}{4144sp}%
\begingroup\makeatletter\ifx\SetFigFont\undefined%
\gdef\SetFigFont#1#2#3#4#5{%
  \reset@font\fontsize{#1}{#2pt}%
  \fontfamily{#3}\fontseries{#4}\fontshape{#5}%
  \selectfont}%
\fi\endgroup%
\begin{picture}(1014,1014)(79,-253)
\put(586,209){\makebox(0,0)[b]{\smash{{\SetFigFont{10}{12.0}{\rmdefault}{\mddefault}{\updefault}{\color[rgb]{0,0,0}$\alpha$}%
}}}}
\put(586,524){\makebox(0,0)[b]{\smash{{\SetFigFont{10}{12.0}{\rmdefault}{\mddefault}{\updefault}{\color[rgb]{0,0,0}$f$}%
}}}}
\put(586,-106){\makebox(0,0)[b]{\smash{{\SetFigFont{10}{12.0}{\rmdefault}{\mddefault}{\updefault}{\color[rgb]{0,0,0}$g$}%
}}}}
\end{picture}%
\end{center}

\noindent The triples are considered \emph{modulo} the following \emph{deformation} equations, given for every possible $2$-arrows~$f$,~$g$ and~$h$ and $3$-cell $\alpha$:
\begin{center}
\begin{picture}(0,0)%
\includegraphics{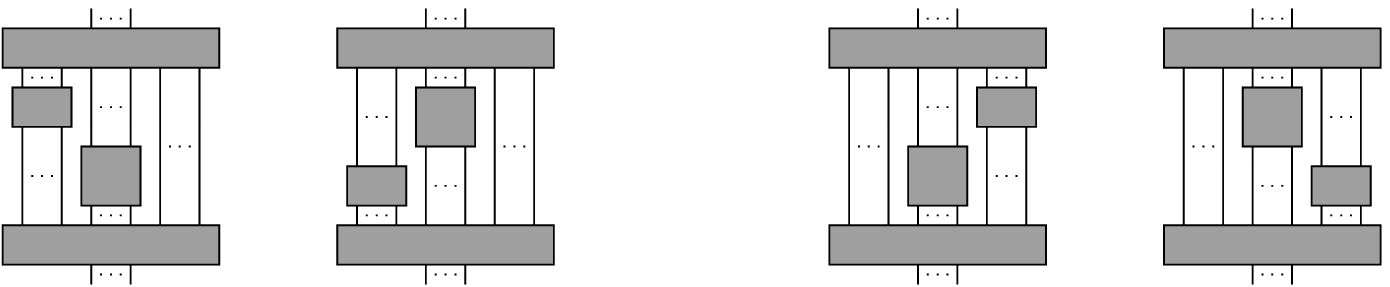}%
\end{picture}%
\setlength{\unitlength}{4144sp}%
\begingroup\makeatletter\ifx\SetFigFont\undefined%
\gdef\SetFigFont#1#2#3#4#5{%
  \reset@font\fontsize{#1}{#2pt}%
  \fontfamily{#3}\fontseries{#4}\fontshape{#5}%
  \selectfont}%
\fi\endgroup%
\begin{picture}(6324,1284)(79,-253)
\put(5896,479){\makebox(0,0)[b]{\smash{{\SetFigFont{10}{12.0}{\rmdefault}{\mddefault}{\updefault}{\color[rgb]{0,0,0}$\alpha$}%
}}}}
\put(1351,389){\makebox(0,0)[b]{\smash{{\SetFigFont{10}{12.0}{\rmdefault}{\mddefault}{\updefault}{\color[rgb]{0,0,0}$=$}%
}}}}
\put(5131,389){\makebox(0,0)[b]{\smash{{\SetFigFont{10}{12.0}{\rmdefault}{\mddefault}{\updefault}{\color[rgb]{0,0,0}$=$}%
}}}}
\put(586,794){\makebox(0,0)[b]{\smash{{\SetFigFont{10}{12.0}{\rmdefault}{\mddefault}{\updefault}{\color[rgb]{0,0,0}$f$}%
}}}}
\put(2116,794){\makebox(0,0)[b]{\smash{{\SetFigFont{10}{12.0}{\rmdefault}{\mddefault}{\updefault}{\color[rgb]{0,0,0}$f$}%
}}}}
\put(4366,794){\makebox(0,0)[b]{\smash{{\SetFigFont{10}{12.0}{\rmdefault}{\mddefault}{\updefault}{\color[rgb]{0,0,0}$f$}%
}}}}
\put(5896,794){\makebox(0,0)[b]{\smash{{\SetFigFont{10}{12.0}{\rmdefault}{\mddefault}{\updefault}{\color[rgb]{0,0,0}$f$}%
}}}}
\put(586,-106){\makebox(0,0)[b]{\smash{{\SetFigFont{10}{12.0}{\rmdefault}{\mddefault}{\updefault}{\color[rgb]{0,0,0}$g$}%
}}}}
\put(2116,-106){\makebox(0,0)[b]{\smash{{\SetFigFont{10}{12.0}{\rmdefault}{\mddefault}{\updefault}{\color[rgb]{0,0,0}$g$}%
}}}}
\put(4681,524){\makebox(0,0)[b]{\smash{{\SetFigFont{10}{12.0}{\rmdefault}{\mddefault}{\updefault}{\color[rgb]{0,0,0}$h$}%
}}}}
\put(6211,164){\makebox(0,0)[b]{\smash{{\SetFigFont{10}{12.0}{\rmdefault}{\mddefault}{\updefault}{\color[rgb]{0,0,0}$h$}%
}}}}
\put(271,524){\makebox(0,0)[b]{\smash{{\SetFigFont{10}{12.0}{\rmdefault}{\mddefault}{\updefault}{\color[rgb]{0,0,0}$h$}%
}}}}
\put(1801,164){\makebox(0,0)[b]{\smash{{\SetFigFont{10}{12.0}{\rmdefault}{\mddefault}{\updefault}{\color[rgb]{0,0,0}$h$}%
}}}}
\put(4366,-106){\makebox(0,0)[b]{\smash{{\SetFigFont{10}{12.0}{\rmdefault}{\mddefault}{\updefault}{\color[rgb]{0,0,0}$g$}%
}}}}
\put(5896,-106){\makebox(0,0)[b]{\smash{{\SetFigFont{10}{12.0}{\rmdefault}{\mddefault}{\updefault}{\color[rgb]{0,0,0}$g$}%
}}}}
\put(586,209){\makebox(0,0)[b]{\smash{{\SetFigFont{10}{12.0}{\rmdefault}{\mddefault}{\updefault}{\color[rgb]{0,0,0}$\alpha$}%
}}}}
\put(2116,479){\makebox(0,0)[b]{\smash{{\SetFigFont{10}{12.0}{\rmdefault}{\mddefault}{\updefault}{\color[rgb]{0,0,0}$\alpha$}%
}}}}
\put(4366,209){\makebox(0,0)[b]{\smash{{\SetFigFont{10}{12.0}{\rmdefault}{\mddefault}{\updefault}{\color[rgb]{0,0,0}$\alpha$}%
}}}}
\end{picture}%
\end{center}
\end{enumerate}
\end{defn*}

\noindent We have the following link between the reduction relation generated by a family of rewriting rules and the reduction graph generated by the corresponding $3$-polygraph:

\begin{rem*}
Let $R$ be a family of rewriting rules on the parallel $2$-arrows of a free $2$-category $\mon{\Sigma}$ generated by a $2$-polygraph $\Sigma$. Let us denote by $\Sigma_R$ the $3$-polygraph built from $\Sigma$ extended with the family $R$, which elements are seen as $3$-cells. Then, given $2$-arrows $u$ and $v$ in $\mon{\Sigma}$ and a rule $\alpha$ in $R$, one has $u\red{\alpha}v$ if and only if there exists an arrow of the form $(f,\alpha,g)$ in $G(\Sigma_R)$ from $u$ to $v$.
\end{rem*}

\noindent Hence, the reduction graph is almost the same as the graph of the reduction relation. However, in the former, we give names to reductions: we explicitely give the context of application of each rule, thus making a difference between two applications of the same rule on the same circuit but in different contexts.

Furthermore, this allows one to explicitely equip applications of rules with the structure of the circuits, instead of an implicit compatibility: this is a first step towards the naming of reductions, which will be of great help in order to deal with bureaucracy.

To conclude this paragraph, let us give some additional notations that will be useful in section~\ref{sec:3d}:

\begin{defn*}
Let $\Sigma=(\ast,\Sigma_1,\Sigma_2,\Sigma_3)$ be a $3$-polygraph with one $0$-cell. In $G(\Sigma)$, we denote by $\alpha$ the arrow $(s_2(\alpha),\alpha,t_2(\alpha))$. The operations $\tens$ and $\circ$ are extended this way between a $2$-arrow of $\Sigma$ and an arrow of $G(\Sigma)$:
\begin{enumerate}
\item[-] If $(f,\alpha,g)$ is an arrow of $G(\Sigma)$ and $h$ and $k$ are $2$-arrows of $\Sigma$ such that $t_1(h)=s_1(f)$ and $s_1(k)=t_1(g)$, then:
$$
(f,\alpha,g)\circ h\:=\:(f\circ h,\alpha,g) \et k\circ(f,\alpha,g)=(f,\alpha,k\circ g).
$$

\item[-] If $(f,\alpha,g)$ is an arrow of $G(\Sigma)$ and $h$ is a $2$-arrow of $\Sigma$, then:
$$
h\tens(f,\alpha,g)\:=\:(h\tens f,\alpha,t_1(h)\tens g) \et (f,\alpha,g)\tens h=(f\tens h,\alpha,g\tens t_1(h)).
$$
\end{enumerate}
\end{defn*}

\begin{rem*}
The extension of $\tens$ is not arbitrary since the deformation relations yield the following equalities in $G(\Sigma)$:
$$
(h\tens f,\alpha,t_1(h)\tens g)\:=\:(s_1(h)\tens f,\alpha,h\tens g) \et (f\tens h,\alpha,g\tens t_1(h))\:=\:(f\tens s_1(h),\alpha,g\tens h).
$$
\end{rem*}

\subsection{From proofs to three-dimensional arrows: the formal construction}\label{sub:preuves-3d}

If we apply the results from the first section concerning the structural rules to the inference rules, we get:

\begin{prop*}
The following definition extends $\pi$ into a surjective functor from $\mon{G(\Sigma^K)}$ to $\mon{G^K}$:
\begin{enumerate}
\item[-] If $\alpha$ is in $\Delta$, then $\pi(f,\alpha,g)=\id_{s(f,\alpha,g)}$.
\item[-] If $\alpha$ is in $R$ or $\ol{S}$, then $\pi(f,\alpha,g)$ is the arrow $\alpha:\pi(s(f,\alpha,g))\fl\pi(t(f,\alpha,g))$.
\end{enumerate}
\end{prop*}

\begin{dem}
Since $\pi$ is defined on objects and arrows of the graph $G(\Sigma^K)$, with values into the category~$\mon{G^K}$, it uniquely extends into a functor from $\mon{G(\Sigma^K)}$ to $\mon{G^K}$. Furthermore, we already know that $\pi$ is surjective on objects. Now, let us consider an arrow from $u$ to $v$ in the free category $\mon{G^K}$. Such an arrow is a sequence of reductions from $u$ to $v$, using the rules of either of $R$ or $\ol{S}$:
$$
u=u_0\red{\alpha_1}u_1\red{\alpha_2}\dots\red{\alpha_n}u_n=v.
$$

\noindent Let us consider a $1$-arrow $X$ in $\Sigma^K$ such that $\sharp X$ is greater than each $\sharp u_i$. Then, for any $i$, we know that there exist $2$-arrows $f_i$ and $g_i$ in $\Sigma^K$ such that:
$$
\Phi^X(u_i)\equi{\Delta} f_i\red{\Phi(\alpha_{i+1})} g_i\equi{\Delta} \Phi^X(u_{i+1}).
$$

\noindent Since, for every equation in $E_{\Delta}$, we have in $\Sigma^K$ a $3$-cell going from left to right and a $3$-cell going from right to left, we know that, whenever $f\equi{\Delta}g$ holds, there exists a path in $G(\Sigma^K)$ from $f$ to $g$ that uses only $3$-cells from the family $\Delta$. Hence, we have a path in $G(\Sigma^K)$:
$$
\Phi^X(u_0)\mred{\Delta} f_0 \red{\Phi(\alpha_0)} g_0 \mred{\Delta} \Phi^X(u_1) \mred{\Delta} f_1 \red{\Phi(\alpha_1)} \dots \mred{\Delta} f_{n-1} \red{\Phi(\alpha_n)} g_{n-1} \mred{\Delta} \Phi^X(u_n).
$$

\noindent Since $\pi$ sends each arrow $\mred{\Delta}$ onto an identity, this path is sent by $\pi$ onto the considered arrow of $\mon{G^K}$. Hence $\pi$ is surjective.

\findem\end{dem}

\noindent In order to adapt the vocabulary of proof theory to the $3$-polygraph $\Sigma^K$, we introduce the following:

\begin{defn*}
Let $f$ and $g$ be $2$-arrows of $\Sigma^K$. A \emph{proof from $f$ to $g$} is a path from $f$ to $g$ in the reduction graph $G(\Sigma^K)$. A \emph{complete proof of $f$} is a path from $\top\circ\epsilon_X$ to $f$ in $G(\Sigma^K)$, where $\epsilon_X$ is a generalized eraser from~$X$ to~$\ast$, built as the juxtaposition of elementary erasers.
\end{defn*}

\noindent As a corollary of the previous result, we get:

\begin{thm*}\label{thm:1}
If there exists a SKS proof from $u$ to $v$, then there exists a proof from $\Phi^X(u)$ to $\Phi^X(v)$ in $\Sigma^K$, for every $1$-arrow $X$ such that $\sharp X$ is greater than both $\sharp u$ and $\sharp v$. In particular, if there exists a complete SKS proof of $u$, then there exists a complete proof of every $\Phi^X(u)$, with $X$ such that $\sharp X\geq\sharp u$.

Conversely, if $f$ and $g$ are $2$-arrows with target $F$ or $A$ such that there exists a proof from $f$ to $g$ in~$\Sigma^K$, then there exists a SKS proof from $\pi(f)$ to $\pi(g)$. In particular, if there exists a complete proof of $f$ with target $A$ or $F$, then there exists a complete proof of $\pi(f)$.
\end{thm*}

\noindent To informally summarize this result, one can say that the proof theory of the $3$-polygraph $\Sigma^K$ we have built is the same one as the proof theory of SKS. Hence, we have a polygraphic translation of the system SKS in what we now call its calculus of structures version.

But the $3$-dimensional setting has not really been used for the moment. And, as we are going to see, the unveiling of the three dimensions of proofs allows a direct and simple control on structural bureaucracy.

\section{Three dimensions against structural bureaucracy}\label{sec:3d}

\noindent In [Guglielmi 2005], objects called Formalism A and Formalism B are sketched in order to identify proofs that only differ by structural bureaucracy: this means that the two proofs only differ by the order of application of the same inference rules.

Defining relations that control this bureaucracy may be difficult in the term-like language of the calculus of structures. Indeed, it is much like the classification of branchings generated by a term rewriting system [Baader Nipkow 1998].

Here theorem~\ref{thm:2} states that, once proofs have been translated into $3$-dimensional objects, the equations controlling structural bureaucracy (\ref{sub:bureaucratie}) become really simple to define: they are the equations called \emph{exchange relations} (\ref{sub:echange}).

After the proof of the theorem (\ref{sub:bureaucratie-est-echange}), we conclude the section by a diagram showing the respective positions of the $3$-polygraphs corresponding to SKS and to Formalisms A and B (\ref{sub:geographie}).

\subsection{The two types of structural bureaucracy}\label{sub:bureaucratie}

\noindent Let us start by giving a definition of structural bureaucracy on SKS proofs, which comes in two types, called A and B. The first one is generated by the applications of two inference rules in different subterms. The second one is generated by the application of two inference rules, one inside the other. In both cases, the two rules apply in two zones of the term that do not intersect.

However, this intuitively simple idea is hard to formalize in the term-like setting used by the calculus of structures: it is like the classification of branchings generated by a term rewriting system, involving many tricky notions such as the relative positions of redexes.

On the other hand, the higher-dimensional setting makes the definitions almost trivial: this is mainly due to the facts that, with this point of view, applications of inference rules have been given a name and that both dimensions of the terms are revealed and treated symmetrically.

Here we use the $2$-categorical structure of $\Tb$ to define both bureaucracy relations. The bureaucracy A relation identifies two proofs that differ by the order of application of two rules in two different subterms:

\begin{defn*}
The \emph{bureaucracy type A relation} is the equivalence relation $\equi{A}$ on SKS proofs generated by the rule $\red{A}$ defined, for every two rules $\alpha$ and $\beta$, every three $1$-arrows $X$, $Y$ and $Z$, every two families of terms $u$ and $v$ by the following diagram (when it has a meaning):
$$
\xymatrix{u\circ(X\tens s(\alpha)\tens Y\tens s(\beta)\tens Z)\circ v \ar@{>>}+<3.2cm,0cm>;[rrr]-<3.2cm,0cm>^-{\alpha} \ar@{>>}-<0cm,0.5cm>;[dd]+<0cm,0.5cm>_-{\beta} &&& u\circ(X\tens t(\alpha)\tens Y\tens s(\beta)\tens Z)\circ v \ar@{>>}-<0cm,0.5cm>;[dd]+<0cm,0.5cm>^-{\beta} \ar-<3.5cm,1cm>;[ddlll]+<3.5cm,1cm>^-{A} \\ \\ u\circ(X\tens s(\alpha)\tens Y\tens t(\beta)\tens Z)\circ v \ar@{>>}+<3.2cm,0cm>;[rrr]-<3.2cm,0cm>_-{\alpha} &&& u\circ(X\tens t(\alpha)\tens Y\tens t(\beta)\tens Z)\circ v.}
$$
\end{defn*}

\begin{rem*}
We use the relations $\mred{\alpha}$ and $\mred{\beta}$ instead of $\red{\alpha}$ and $\red{\beta}$. Indeed, when the given factorizations of the terms are projected through $\pi$ onto families of terms, some duplicators or erasers implicitely present in $u$ may duplicate or erase the redexes $s(\alpha)$ and $s(\beta)$. Hence, reducing them may require more or less than one application of either $\alpha$ or $\beta$. The same comment applies to the next definition.
\end{rem*}

\noindent The bureaucracy B relation identifies two proofs that differ by the order of application of two rules, one inside the other:

\begin{defn*}
The \emph{bureaucracy type B relation} is the equivalence relation $\equi{B}$ on SKS proofs generated by the rule $\red{B}$ defined, for every two rules $\alpha$ and $\beta$, every four $1$-arrows $X_1$, $X_2$, $Y_1$ and $Y_2$, every three families of terms $u$, $v$ and $w$ by the following diagram (when it has a meaning):
$$
\xymatrix{*\txt{$u\circ(X_1\tens s(\alpha)\tens X_2)\circ v$\\$\circ(Y_1\tens s(\beta)\tens Y_2)\circ w$} \ar@{>>}+<2.5cm,0cm>;[rrr]-<2.5cm,0cm>^-{\alpha} \ar@{>>}-<0cm,0.6cm>;[dd]+<0cm,0.6cm>_-{\beta} &&& *\txt{$u\circ(X_1\tens t(\alpha)\tens X_2)\circ v$\\$\circ(Y_1\tens s(\beta)\tens Y_2)\circ w$} \ar@{>>}-<0cm,0.6cm>;[dd]+<0cm,0.6cm>^-{\beta} \ar-<2.5cm,1cm>;[ddlll]+<2.5cm,1cm>^-{B} \\ \\ *\txt{$u\circ(X_1\tens s(\alpha)\tens X_2)\circ v$\\$\circ(Y_1\tens t(\beta)\tens Y_2)\circ w$} \ar@{>>}+<2.5cm,0cm>;[rrr]-<2.5cm,0cm>_-{\alpha} &&& *\txt{$u\circ(X_1\tens t(\alpha)\tens X_2)\circ v$\\$\circ(Y_1\tens t(\beta)\tens Y_2)\circ w$}.}
$$

\noindent The \emph{structural bureaucracy relation} is the equivalence relation generated by the union of $\red{A}$ and~$\red{B}$.
\end{defn*}

\begin{rem*}
One observation one can make is that the two bureaucratic relations appear to be different in essence. However, this is an artifact of the term-like notation: in the polygraphic setting, both have the same simple shape. Another observation one does is that these definitions are quite technical (and their version without the two compositions available would be even worse).

Once again, this is due to the term structure, since the bureaucratic relations are really easy to define in the polygraphic setting as we are going to see now. Even better, there they inherit the geometrical interpretation they deserve: they appear as the ability to move blocks representing subproofs one around the other.
\end{rem*}

\subsection{Exchange relations and three-categories}\label{sub:echange}

In this paragraph, we give polygraphic equivalents of the bureaucratic relations. Let us consider the idea behind the definition of bureaucracy A: we want to identify two proofs that only differ by the order of application of two rules in two different subterms. And, in circuits, different subterms are two juxtaposed subcircuits. Hence, bureaucracy A on circuits should identify the two following paths of the reduction graph $G(\Sigma^K)$:
\begin{center}
\begin{picture}(0,0)%
\includegraphics{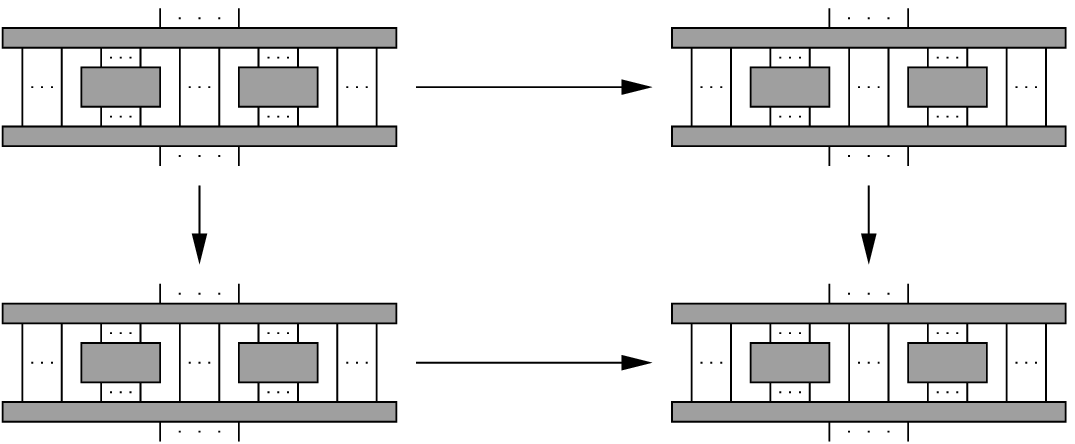}%
\end{picture}%
\setlength{\unitlength}{4144sp}%
\begingroup\makeatletter\ifx\SetFigFont\undefined%
\gdef\SetFigFont#1#2#3#4#5{%
  \reset@font\fontsize{#1}{#2pt}%
  \fontfamily{#3}\fontseries{#4}\fontshape{#5}%
  \selectfont}%
\fi\endgroup%
\begin{picture}(4884,2004)(79,-1333)
\put(4231,-376){\makebox(0,0)[b]{\smash{{\SetFigFont{10}{12.0}{\rmdefault}{\mddefault}{\updefault}{\color[rgb]{0,0,0}$\beta$}%
}}}}
\put(2521,-1186){\makebox(0,0)[b]{\smash{{\SetFigFont{10}{12.0}{\rmdefault}{\mddefault}{\updefault}{\color[rgb]{0,0,0}$\alpha$}%
}}}}
\put(4411,-1006){\makebox(0,0)[b]{\smash{{\SetFigFont{8}{9.6}{\rmdefault}{\mddefault}{\updefault}{\color[rgb]{0,0,0}$t_2(\beta)$}%
}}}}
\put(3691,-1006){\makebox(0,0)[b]{\smash{{\SetFigFont{8}{9.6}{\rmdefault}{\mddefault}{\updefault}{\color[rgb]{0,0,0}$t_2(\alpha)$}%
}}}}
\put(1351,-1006){\makebox(0,0)[b]{\smash{{\SetFigFont{8}{9.6}{\rmdefault}{\mddefault}{\updefault}{\color[rgb]{0,0,0}$t_2(\beta)$}%
}}}}
\put(631,-1006){\makebox(0,0)[b]{\smash{{\SetFigFont{8}{9.6}{\rmdefault}{\mddefault}{\updefault}{\color[rgb]{0,0,0}$s_2(\alpha)$}%
}}}}
\put(2521,389){\makebox(0,0)[b]{\smash{{\SetFigFont{10}{12.0}{\rmdefault}{\mddefault}{\updefault}{\color[rgb]{0,0,0}$\alpha$}%
}}}}
\put(631,254){\makebox(0,0)[b]{\smash{{\SetFigFont{8}{9.6}{\rmdefault}{\mddefault}{\updefault}{\color[rgb]{0,0,0}$s_2(\alpha)$}%
}}}}
\put(1351,254){\makebox(0,0)[b]{\smash{{\SetFigFont{8}{9.6}{\rmdefault}{\mddefault}{\updefault}{\color[rgb]{0,0,0}$s_2(\beta)$}%
}}}}
\put(3691,254){\makebox(0,0)[b]{\smash{{\SetFigFont{8}{9.6}{\rmdefault}{\mddefault}{\updefault}{\color[rgb]{0,0,0}$t_2(\alpha)$}%
}}}}
\put(4411,254){\makebox(0,0)[b]{\smash{{\SetFigFont{8}{9.6}{\rmdefault}{\mddefault}{\updefault}{\color[rgb]{0,0,0}$s_2(\beta)$}%
}}}}
\put(811,-376){\makebox(0,0)[b]{\smash{{\SetFigFont{10}{12.0}{\rmdefault}{\mddefault}{\updefault}{\color[rgb]{0,0,0}$\beta$}%
}}}}
\put(2521,-376){\makebox(0,0)[b]{\smash{{\SetFigFont{10}{12.0}{\rmdefault}{\mddefault}{\updefault}{\color[rgb]{0,0,0}$\equiv$}%
}}}}
\end{picture}%
\end{center}

\vfill\pagebreak
\noindent By removing all unnecessary contexts, we get that this relation is generated on circuits by the following smaller one, indexed by pairs $(\alpha:f\fl f',\beta:g\fl g')$ of arrows in $G(\Sigma^K)$:
\begin{center}
\begin{picture}(0,0)%
\includegraphics{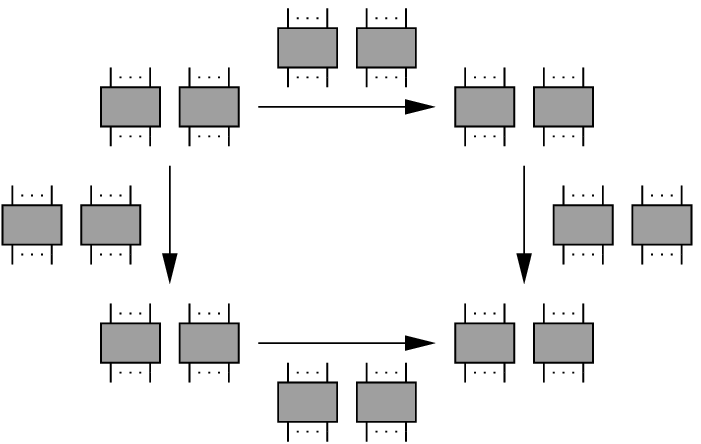}%
\end{picture}%
\setlength{\unitlength}{4144sp}%
\begingroup\makeatletter\ifx\SetFigFont\undefined%
\gdef\SetFigFont#1#2#3#4#5{%
  \reset@font\fontsize{#1}{#2pt}%
  \fontfamily{#3}\fontseries{#4}\fontshape{#5}%
  \selectfont}%
\fi\endgroup%
\begin{picture}(3174,2004)(-371,-973)
\put(1216,-16){\makebox(0,0)[b]{\smash{{\SetFigFont{10}{12.0}{\rmdefault}{\mddefault}{\updefault}{\color[rgb]{0,0,0}$\equiv$}%
}}}}
\put(226,524){\makebox(0,0)[b]{\smash{{\SetFigFont{10}{12.0}{\rmdefault}{\mddefault}{\updefault}{\color[rgb]{0,0,0}$f$}%
}}}}
\put(226,-556){\makebox(0,0)[b]{\smash{{\SetFigFont{10}{12.0}{\rmdefault}{\mddefault}{\updefault}{\color[rgb]{0,0,0}$f$}%
}}}}
\put(-224,-16){\makebox(0,0)[b]{\smash{{\SetFigFont{10}{12.0}{\rmdefault}{\mddefault}{\updefault}{\color[rgb]{0,0,0}$f$}%
}}}}
\put(1846,524){\makebox(0,0)[b]{\smash{{\SetFigFont{10}{12.0}{\rmdefault}{\mddefault}{\updefault}{\color[rgb]{0,0,0}$f'$}%
}}}}
\put(1846,-556){\makebox(0,0)[b]{\smash{{\SetFigFont{10}{12.0}{\rmdefault}{\mddefault}{\updefault}{\color[rgb]{0,0,0}$f'$}%
}}}}
\put(2296,-16){\makebox(0,0)[b]{\smash{{\SetFigFont{10}{12.0}{\rmdefault}{\mddefault}{\updefault}{\color[rgb]{0,0,0}$f'$}%
}}}}
\put(586,524){\makebox(0,0)[b]{\smash{{\SetFigFont{10}{12.0}{\rmdefault}{\mddefault}{\updefault}{\color[rgb]{0,0,0}$g$}%
}}}}
\put(2206,524){\makebox(0,0)[b]{\smash{{\SetFigFont{10}{12.0}{\rmdefault}{\mddefault}{\updefault}{\color[rgb]{0,0,0}$g$}%
}}}}
\put(1396,794){\makebox(0,0)[b]{\smash{{\SetFigFont{10}{12.0}{\rmdefault}{\mddefault}{\updefault}{\color[rgb]{0,0,0}$g$}%
}}}}
\put(586,-556){\makebox(0,0)[b]{\smash{{\SetFigFont{10}{12.0}{\rmdefault}{\mddefault}{\updefault}{\color[rgb]{0,0,0}$g'$}%
}}}}
\put(2206,-556){\makebox(0,0)[b]{\smash{{\SetFigFont{10}{12.0}{\rmdefault}{\mddefault}{\updefault}{\color[rgb]{0,0,0}$g'$}%
}}}}
\put(1396,-826){\makebox(0,0)[b]{\smash{{\SetFigFont{10}{12.0}{\rmdefault}{\mddefault}{\updefault}{\color[rgb]{0,0,0}$g'$}%
}}}}
\put(1036,794){\makebox(0,0)[b]{\smash{{\SetFigFont{10}{12.0}{\rmdefault}{\mddefault}{\updefault}{\color[rgb]{0,0,0}$\alpha$}%
}}}}
\put(1036,-826){\makebox(0,0)[b]{\smash{{\SetFigFont{10}{12.0}{\rmdefault}{\mddefault}{\updefault}{\color[rgb]{0,0,0}$\alpha$}%
}}}}
\put(136,-16){\makebox(0,0)[b]{\smash{{\SetFigFont{10}{12.0}{\rmdefault}{\mddefault}{\updefault}{\color[rgb]{0,0,0}$\beta$}%
}}}}
\put(2656,-16){\makebox(0,0)[b]{\smash{{\SetFigFont{10}{12.0}{\rmdefault}{\mddefault}{\updefault}{\color[rgb]{0,0,0}$\beta$}%
}}}}
\end{picture}%
\end{center}

\noindent If one considers the block-like $3$-dimensional representation of $3$-cells, one gets the following identification, for every pair $(\alpha:f\fl f',\beta:g\fl g')$ of arrows in $G(\Sigma^K)$. The corresponding relation is written below, where $\star$ denotes the composition of paths in $\mon{G(\Sigma^K)}$:
\begin{center}
\begin{picture}(0,0)%
\includegraphics{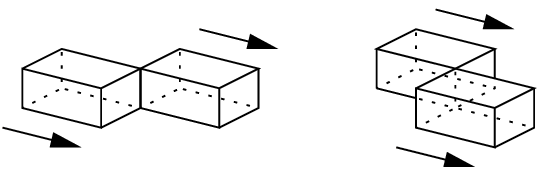}%
\end{picture}%
\setlength{\unitlength}{4144sp}%
\begingroup\makeatletter\ifx\SetFigFont\undefined%
\gdef\SetFigFont#1#2#3#4#5{%
  \reset@font\fontsize{#1}{#2pt}%
  \fontfamily{#3}\fontseries{#4}\fontshape{#5}%
  \selectfont}%
\fi\endgroup%
\begin{picture}(2454,1429)(79,-1325)
\put(1261,-1276){\makebox(0,0)[rb]{\smash{{\SetFigFont{10}{12.0}{\rmdefault}{\mddefault}{\updefault}{\color[rgb]{0,0,0}$(\alpha\tens s_2(\beta))\star(t_2(\alpha)\tens\beta)$}%
}}}}
\put(1531,-466){\makebox(0,0)[b]{\smash{{\SetFigFont{10}{12.0}{\rmdefault}{\mddefault}{\updefault}$\equiv$}}}}
\put(181,-826){\makebox(0,0)[b]{\smash{{\SetFigFont{10}{12.0}{\rmdefault}{\mddefault}{\updefault}{\color[rgb]{0,0,0}$\alpha$}%
}}}}
\put(1981,-916){\makebox(0,0)[b]{\smash{{\SetFigFont{10}{12.0}{\rmdefault}{\mddefault}{\updefault}{\color[rgb]{0,0,0}$\beta$}%
}}}}
\put(1171,-106){\makebox(0,0)[b]{\smash{{\SetFigFont{10}{12.0}{\rmdefault}{\mddefault}{\updefault}{\color[rgb]{0,0,0}$\beta$}%
}}}}
\put(2251,-16){\makebox(0,0)[b]{\smash{{\SetFigFont{10}{12.0}{\rmdefault}{\mddefault}{\updefault}{\color[rgb]{0,0,0}$\alpha$}%
}}}}
\put(1531,-1276){\makebox(0,0)[b]{\smash{{\SetFigFont{10}{12.0}{\rmdefault}{\mddefault}{\updefault}{\color[rgb]{0,0,0}$\equiv$}%
}}}}
\put(1801,-1276){\makebox(0,0)[lb]{\smash{{\SetFigFont{10}{12.0}{\rmdefault}{\mddefault}{\updefault}{\color[rgb]{0,0,0}$(s_2(\alpha)\tens\beta)\star(\alpha\tens t_2(\beta))$}%
}}}}
\end{picture}%
\end{center}

\noindent Now, let us translate the bureaucracy B relation onto paths in the reduction graph $G(\Sigma^K)$. This relation should identify proofs that only differ by the order of application of two rules, one inside the other one. On circuits, this means that the two rules act on circuits that are vertically composed. Thus, bureaucracy type~B on circuits should identify the following paths of $G(\Sigma^K)$:
\begin{center}
\begin{picture}(0,0)%
\includegraphics{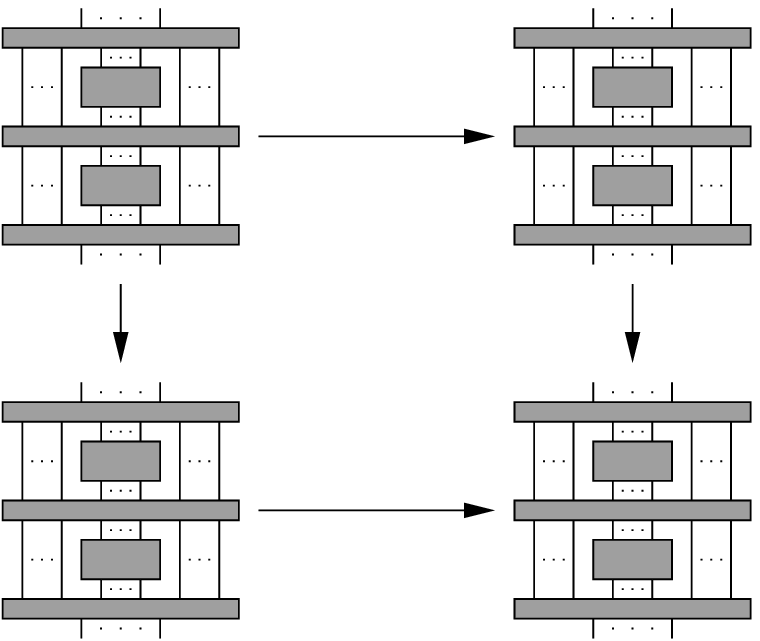}%
\end{picture}%
\setlength{\unitlength}{4144sp}%
\begingroup\makeatletter\ifx\SetFigFont\undefined%
\gdef\SetFigFont#1#2#3#4#5{%
  \reset@font\fontsize{#1}{#2pt}%
  \fontfamily{#3}\fontseries{#4}\fontshape{#5}%
  \selectfont}%
\fi\endgroup%
\begin{picture}(3444,2904)(79,-1693)
\put(631,344){\makebox(0,0)[b]{\smash{{\SetFigFont{8}{9.6}{\rmdefault}{\mddefault}{\updefault}{\color[rgb]{0,0,0}$s_2(\beta)$}%
}}}}
\put(1756,704){\makebox(0,0)[b]{\smash{{\SetFigFont{10}{12.0}{\rmdefault}{\mddefault}{\updefault}{\color[rgb]{0,0,0}$\alpha$}%
}}}}
\put(451,-286){\makebox(0,0)[b]{\smash{{\SetFigFont{10}{12.0}{\rmdefault}{\mddefault}{\updefault}{\color[rgb]{0,0,0}$\beta$}%
}}}}
\put(3151,-286){\makebox(0,0)[b]{\smash{{\SetFigFont{10}{12.0}{\rmdefault}{\mddefault}{\updefault}{\color[rgb]{0,0,0}$\beta$}%
}}}}
\put(1756,-286){\makebox(0,0)[b]{\smash{{\SetFigFont{10}{12.0}{\rmdefault}{\mddefault}{\updefault}{\color[rgb]{0,0,0}$\equiv$}%
}}}}
\put(631,794){\makebox(0,0)[b]{\smash{{\SetFigFont{8}{9.6}{\rmdefault}{\mddefault}{\updefault}{\color[rgb]{0,0,0}$s_2(\alpha)$}%
}}}}
\put(2971,794){\makebox(0,0)[b]{\smash{{\SetFigFont{8}{9.6}{\rmdefault}{\mddefault}{\updefault}{\color[rgb]{0,0,0}$t_2(\alpha)$}%
}}}}
\put(1756,-1276){\makebox(0,0)[b]{\smash{{\SetFigFont{10}{12.0}{\rmdefault}{\mddefault}{\updefault}{\color[rgb]{0,0,0}$\alpha$}%
}}}}
\put(631,-1366){\makebox(0,0)[b]{\smash{{\SetFigFont{8}{9.6}{\rmdefault}{\mddefault}{\updefault}{\color[rgb]{0,0,0}$t_2(\beta)$}%
}}}}
\put(631,-916){\makebox(0,0)[b]{\smash{{\SetFigFont{8}{9.6}{\rmdefault}{\mddefault}{\updefault}{\color[rgb]{0,0,0}$s_2(\alpha)$}%
}}}}
\put(2971,-916){\makebox(0,0)[b]{\smash{{\SetFigFont{8}{9.6}{\rmdefault}{\mddefault}{\updefault}{\color[rgb]{0,0,0}$t_2(\alpha)$}%
}}}}
\put(2971,-1366){\makebox(0,0)[b]{\smash{{\SetFigFont{8}{9.6}{\rmdefault}{\mddefault}{\updefault}{\color[rgb]{0,0,0}$t_2(\beta)$}%
}}}}
\put(2971,344){\makebox(0,0)[b]{\smash{{\SetFigFont{8}{9.6}{\rmdefault}{\mddefault}{\updefault}{\color[rgb]{0,0,0}$s_2(\beta)$}%
}}}}
\end{picture}%
\end{center}

\vfill\pagebreak
\noindent Once again, this relation is generated on circuits by a smaller one, with unnecessary contexts removed, given for every pair $(\alpha:f\fl f',\beta:g\fl g')$ of arrows of $G(\Sigma^K)$:
\begin{center}
\begin{picture}(0,0)%
\includegraphics{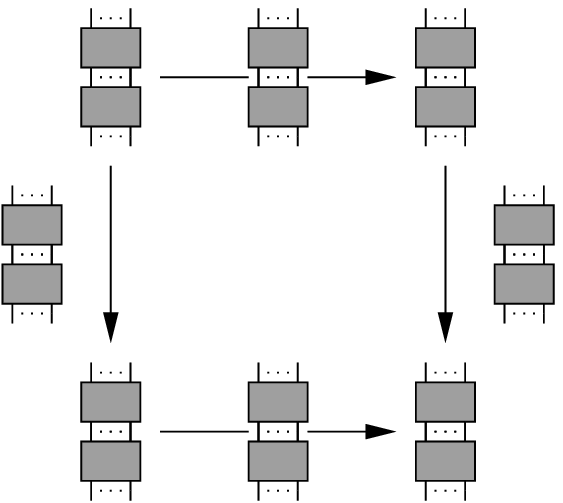}%
\end{picture}%
\setlength{\unitlength}{4144sp}%
\begingroup\makeatletter\ifx\SetFigFont\undefined%
\gdef\SetFigFont#1#2#3#4#5{%
  \reset@font\fontsize{#1}{#2pt}%
  \fontfamily{#3}\fontseries{#4}\fontshape{#5}%
  \selectfont}%
\fi\endgroup%
\begin{picture}(2544,2274)(-101,-1243)
\put(1171,-151){\makebox(0,0)[b]{\smash{{\SetFigFont{10}{12.0}{\rmdefault}{\mddefault}{\updefault}{\color[rgb]{0,0,0}$\equiv$}%
}}}}
\put(406,524){\makebox(0,0)[b]{\smash{{\SetFigFont{10}{12.0}{\rmdefault}{\mddefault}{\updefault}{\color[rgb]{0,0,0}$f$}%
}}}}
\put(406,794){\makebox(0,0)[b]{\smash{{\SetFigFont{10}{12.0}{\rmdefault}{\mddefault}{\updefault}{\color[rgb]{0,0,0}$g$}%
}}}}
\put(406,-826){\makebox(0,0)[b]{\smash{{\SetFigFont{10}{12.0}{\rmdefault}{\mddefault}{\updefault}{\color[rgb]{0,0,0}$g'$}%
}}}}
\put(406,-1096){\makebox(0,0)[b]{\smash{{\SetFigFont{10}{12.0}{\rmdefault}{\mddefault}{\updefault}{\color[rgb]{0,0,0}$f$}%
}}}}
\put( 46,-16){\makebox(0,0)[b]{\smash{{\SetFigFont{10}{12.0}{\rmdefault}{\mddefault}{\updefault}{\color[rgb]{0,0,0}$\beta$}%
}}}}
\put( 46,-286){\makebox(0,0)[b]{\smash{{\SetFigFont{10}{12.0}{\rmdefault}{\mddefault}{\updefault}{\color[rgb]{0,0,0}$f$}%
}}}}
\put(1936,794){\makebox(0,0)[b]{\smash{{\SetFigFont{10}{12.0}{\rmdefault}{\mddefault}{\updefault}{\color[rgb]{0,0,0}$g$}%
}}}}
\put(1936,524){\makebox(0,0)[b]{\smash{{\SetFigFont{10}{12.0}{\rmdefault}{\mddefault}{\updefault}{\color[rgb]{0,0,0}$f'$}%
}}}}
\put(1936,-826){\makebox(0,0)[b]{\smash{{\SetFigFont{10}{12.0}{\rmdefault}{\mddefault}{\updefault}{\color[rgb]{0,0,0}$g'$}%
}}}}
\put(1936,-1096){\makebox(0,0)[b]{\smash{{\SetFigFont{10}{12.0}{\rmdefault}{\mddefault}{\updefault}{\color[rgb]{0,0,0}$f'$}%
}}}}
\put(2296,-286){\makebox(0,0)[b]{\smash{{\SetFigFont{10}{12.0}{\rmdefault}{\mddefault}{\updefault}{\color[rgb]{0,0,0}$f'$}%
}}}}
\put(2296,-16){\makebox(0,0)[b]{\smash{{\SetFigFont{10}{12.0}{\rmdefault}{\mddefault}{\updefault}{\color[rgb]{0,0,0}$\beta$}%
}}}}
\put(1171,-826){\makebox(0,0)[b]{\smash{{\SetFigFont{10}{12.0}{\rmdefault}{\mddefault}{\updefault}{\color[rgb]{0,0,0}$g'$}%
}}}}
\put(1171,-1096){\makebox(0,0)[b]{\smash{{\SetFigFont{10}{12.0}{\rmdefault}{\mddefault}{\updefault}{\color[rgb]{0,0,0}$\alpha$}%
}}}}
\put(1171,524){\makebox(0,0)[b]{\smash{{\SetFigFont{10}{12.0}{\rmdefault}{\mddefault}{\updefault}{\color[rgb]{0,0,0}$\alpha$}%
}}}}
\put(1171,794){\makebox(0,0)[b]{\smash{{\SetFigFont{10}{12.0}{\rmdefault}{\mddefault}{\updefault}{\color[rgb]{0,0,0}$g$}%
}}}}
\end{picture}%
\end{center}

\noindent In block representation, one gets, for every pair $(\alpha,\beta)$ of arrows of $G(\Sigma^K)$ - once again, this representation is only given to favour the geometrical intuition:
\begin{center}
\begin{picture}(0,0)%
\includegraphics{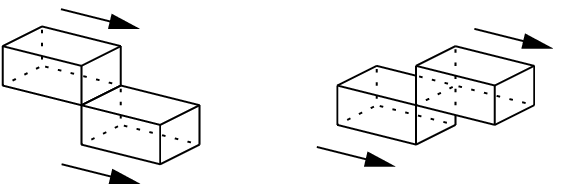}%
\end{picture}%
\setlength{\unitlength}{4144sp}%
\begingroup\makeatletter\ifx\SetFigFont\undefined%
\gdef\SetFigFont#1#2#3#4#5{%
  \reset@font\fontsize{#1}{#2pt}%
  \fontfamily{#3}\fontseries{#4}\fontshape{#5}%
  \selectfont}%
\fi\endgroup%
\begin{picture}(2541,1462)(349,-2000)
\put(1576,-1951){\makebox(0,0)[b]{\smash{{\SetFigFont{10}{12.0}{\rmdefault}{\mddefault}{\updefault}$\equiv$}}}}
\put(1576,-1141){\makebox(0,0)[b]{\smash{{\SetFigFont{10}{12.0}{\rmdefault}{\mddefault}{\updefault}$\equiv$}}}}
\put(721,-1636){\makebox(0,0)[b]{\smash{{\SetFigFont{10}{12.0}{\rmdefault}{\mddefault}{\updefault}{\color[rgb]{0,0,0}$\alpha$}%
}}}}
\put(811,-646){\makebox(0,0)[b]{\smash{{\SetFigFont{10}{12.0}{\rmdefault}{\mddefault}{\updefault}{\color[rgb]{0,0,0}$\beta$}%
}}}}
\put(2701,-736){\makebox(0,0)[b]{\smash{{\SetFigFont{10}{12.0}{\rmdefault}{\mddefault}{\updefault}{\color[rgb]{0,0,0}$\beta$}%
}}}}
\put(1891,-1546){\makebox(0,0)[b]{\smash{{\SetFigFont{10}{12.0}{\rmdefault}{\mddefault}{\updefault}{\color[rgb]{0,0,0}$\alpha$}%
}}}}
\put(1306,-1951){\makebox(0,0)[rb]{\smash{{\SetFigFont{10}{12.0}{\rmdefault}{\mddefault}{\updefault}{\color[rgb]{0,0,0}$(s_2(\alpha)\circ\beta)\star(\alpha\circ t_2(\beta))$}%
}}}}
\put(1891,-1951){\makebox(0,0)[lb]{\smash{{\SetFigFont{10}{12.0}{\rmdefault}{\mddefault}{\updefault}{\color[rgb]{0,0,0}$(\alpha\circ s_2(\beta))\star(t_2(\alpha)\circ\beta)$}%
}}}}
\end{picture}%
\end{center}

\noindent The two families of relations we have exhibited are called \emph{exchange relations}. They are exactly what lies between the free category $\mon{G(\Sigma^K)}$ generated by the reduction graph~$G(\Sigma^K)$ and the \emph{free $3$-category generated by the $3$-polygraph $\Sigma^K$}.

Here we give only a graphical definition of this notion, while a formal one is in [Burroni 1993]. Thereafter, we write $f\star_0 g$, $f\star_1 g$ and~$f\star_2 g$ for~$f\tens g$, $g\circ f$ and $f\star g$ respectively.

\begin{defn*}
Let $\Sigma$ be a $3$-polygraph. The \emph{free $\mathit{3}$-category generated by $\Sigma$} is denoted by $\mon{\Sigma}$ and is made of the $0$, $1$ and $2$-arrows of $\Sigma$, together with a family of \emph{$\mathit{3}$-arrows} which are the paths of the reduction graph~$G(\Sigma)$ \emph{modulo} the congruence generated by the exchange relations:
$$
\begin{aligned}
(\alpha\star_0 s_2(\beta))\star_2(t_2(\alpha)\star_0\beta) &\quad\equi{02}\quad (s_2(\alpha)\star_0\beta)\star_2(\alpha\star_0 t_2(\beta)), \\
(\alpha\star_1 s_2(\beta))\star_2(t_2(\alpha)\star_1\beta) &\quad\equi{12}\quad (s_2(\alpha)\star_1\beta)\star_2(\alpha\star_1 t_2(\beta)).
\end{aligned}
$$

\noindent These equations allow one to extend the compositions $\star_0$ and $\star_1$ on equivalence classes of paths of $G(\Sigma)$ with $\alpha\star_0\beta$ being given by any side of the relation $\equi{02}$ and $\alpha\star_1\beta$ by any side of $\equi{12}$.
\end{defn*}

\noindent Using the deformation relation already defined on arrows of $G(\Sigma)$ together with the two exchange relations $\equi{02}$ and $\equi{12}$, one proves a third exchange relation $\equi{01}$. Note that we need the extensions of the two compositions $\star_0$ and $\star_1$ allowed by $\equi{02}$ and $\equi{12}$ to write this new one.

\begin{lem*}
Let $\Sigma$ be a $3$-polygraph. In the free $3$-category $\mon{\Sigma}$ generated by $\Sigma$, the following exchange relation holds for any $3$-arrows $\alpha$ and $\beta$, both sides being equal to $\alpha\star_0\beta$:
$$
(\alpha\star_0 s_1(\beta))\star_1(t_1(\alpha)\star_0\beta) \quad\equi{01}\quad (s_1(\alpha)\star_0\beta)\star_1(\alpha\star_0 t_1(\beta)).
$$
\end{lem*}

\begin{dem}
Let us consider $3$-arrows $\alpha$ and $\beta$. Then we have:
$$
\begin{aligned}
(\alpha\star_0 s_1\beta)\star_1(t_1\alpha\star_0\beta) \quad&=\quad \big((\alpha\star_0 s_1\beta)\star_1 s_2(t_1\alpha\star_0\beta)\big) \star_2 \big(t_2(\alpha\star_0 s_1\beta)\star_1(t_1\alpha\star_0\beta)\big) \\
&=\quad \big((\alpha\star_0 s_1s_2\beta)\star_1 (t_1\alpha\star_0 s_2\beta)\big) \star_2 \big((t_2\alpha\star_0 s_1\beta)\star_1(t_1t_2\alpha\star_0\beta)\big) \\
&=\quad (\alpha\star_0 s_2\beta)\star_2(t_2\alpha\star_0\beta) \\
&=\quad \alpha\star_0\beta.
\end{aligned}
$$

\noindent The first equality uses the definition of $\star_1$ on $3$-arrows. The second one is due to commutation properties of the sources, targets and compositions operators [Burroni 1993]. Then the deformation relation on arrows of $G(\Sigma)$ yields the third equality. Finally the relation $\equi{02}$ allows one to conclude. A similar computation gives the other part of the seeked relation.

\findem\end{dem}

\medskip
\noindent Alternatively, we can give a more constructive definition of the free $3$-category $\mon{\Sigma}$ generated by a $3$-polygraph $\Sigma$. Its $3$-arrows are generated by the $3$-cells of $\Sigma$ seen as blocks:
\begin{center}
\begin{picture}(0,0)%
\includegraphics{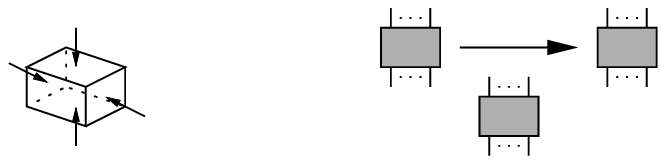}%
\end{picture}%
\setlength{\unitlength}{4144sp}%
\begingroup\makeatletter\ifx\SetFigFont\undefined%
\gdef\SetFigFont#1#2#3#4#5{%
  \reset@font\fontsize{#1}{#2pt}%
  \fontfamily{#3}\fontseries{#4}\fontshape{#5}%
  \selectfont}%
\fi\endgroup%
\begin{picture}(3225,868)(118,-146)
\put(1126, 29){\makebox(0,0)[b]{\smash{{\SetFigFont{10}{12.0}{\rmdefault}{\mddefault}{\updefault}{\color[rgb]{0,0,0}$g$}%
}}}}
\put(676,614){\makebox(0,0)[b]{\smash{{\SetFigFont{10}{12.0}{\rmdefault}{\mddefault}{\updefault}{\color[rgb]{0,0,0}$m$}%
}}}}
\put(676,-106){\makebox(0,0)[b]{\smash{{\SetFigFont{10}{12.0}{\rmdefault}{\mddefault}{\updefault}{\color[rgb]{0,0,0}$n$}%
}}}}
\put(2206,434){\makebox(0,0)[b]{\smash{{\SetFigFont{10}{12.0}{\rmdefault}{\mddefault}{\updefault}$f$}}}}
\put(3196,434){\makebox(0,0)[b]{\smash{{\SetFigFont{10}{12.0}{\rmdefault}{\mddefault}{\updefault}$g$}}}}
\put(2656,119){\makebox(0,0)[b]{\smash{{\SetFigFont{10}{12.0}{\rmdefault}{\mddefault}{\updefault}$F$}}}}
\put(226,389){\makebox(0,0)[b]{\smash{{\SetFigFont{10}{12.0}{\rmdefault}{\mddefault}{\updefault}{\color[rgb]{0,0,0}$f$}%
}}}}
\end{picture}%
\end{center}

\noindent On these generators, one can use the three following constructors, called compositions:
\begin{center}
\begin{picture}(0,0)%
\includegraphics{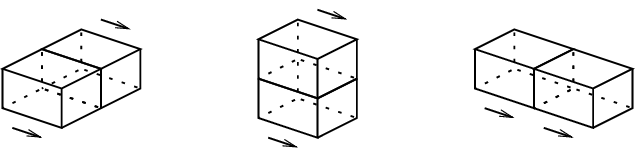}%
\end{picture}%
\setlength{\unitlength}{4144sp}%
\begingroup\makeatletter\ifx\SetFigFont\undefined%
\gdef\SetFigFont#1#2#3#4#5{%
  \reset@font\fontsize{#1}{#2pt}%
  \fontfamily{#3}\fontseries{#4}\fontshape{#5}%
  \selectfont}%
\fi\endgroup%
\begin{picture}(2904,1249)(169,-470)
\put(2701,-421){\makebox(0,0)[b]{\smash{{\SetFigFont{10}{12.0}{\rmdefault}{\mddefault}{\updefault}{\color[rgb]{0,0,0}$\alpha\star_2\beta$}%
}}}}
\put(766,614){\makebox(0,0)[b]{\smash{{\SetFigFont{10}{12.0}{\rmdefault}{\mddefault}{\updefault}{\color[rgb]{0,0,0}$\beta$}%
}}}}
\put(226,-106){\makebox(0,0)[b]{\smash{{\SetFigFont{10}{12.0}{\rmdefault}{\mddefault}{\updefault}{\color[rgb]{0,0,0}$\alpha$}%
}}}}
\put(1396,-151){\makebox(0,0)[b]{\smash{{\SetFigFont{10}{12.0}{\rmdefault}{\mddefault}{\updefault}{\color[rgb]{0,0,0}$\beta$}%
}}}}
\put(2386,-16){\makebox(0,0)[b]{\smash{{\SetFigFont{10}{12.0}{\rmdefault}{\mddefault}{\updefault}{\color[rgb]{0,0,0}$\alpha$}%
}}}}
\put(2656,-106){\makebox(0,0)[b]{\smash{{\SetFigFont{10}{12.0}{\rmdefault}{\mddefault}{\updefault}{\color[rgb]{0,0,0}$\beta$}%
}}}}
\put(1756,659){\makebox(0,0)[b]{\smash{{\SetFigFont{10}{12.0}{\rmdefault}{\mddefault}{\updefault}{\color[rgb]{0,0,0}$\alpha$}%
}}}}
\put(496,-421){\makebox(0,0)[b]{\smash{{\SetFigFont{10}{12.0}{\rmdefault}{\mddefault}{\updefault}{\color[rgb]{0,0,0}$\alpha\star_0\beta$}%
}}}}
\put(1576,-421){\makebox(0,0)[b]{\smash{{\SetFigFont{10}{12.0}{\rmdefault}{\mddefault}{\updefault}{\color[rgb]{0,0,0}$\alpha\star_1\beta$}%
}}}}
\end{picture}%
\end{center}

\noindent If they are sliced, these compositions appear this way:
\begin{center}
\begin{picture}(0,0)%
\includegraphics{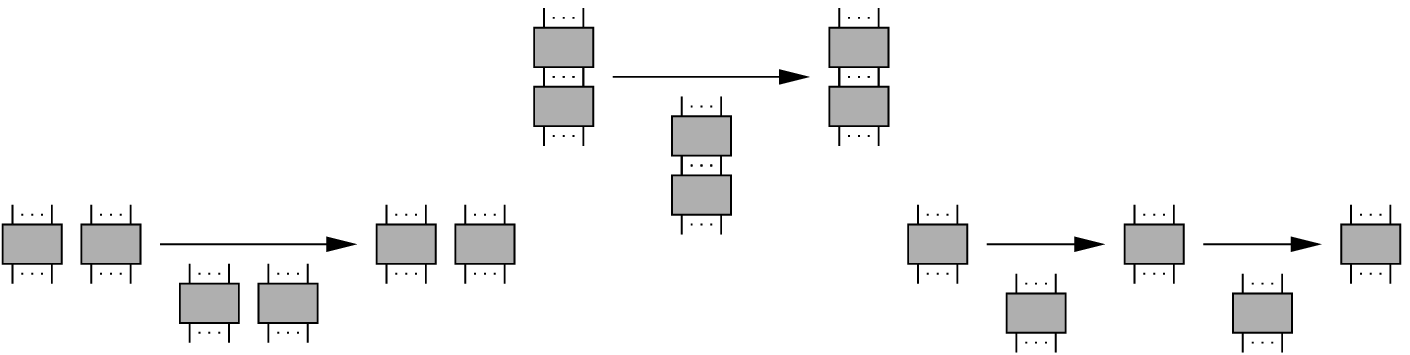}%
\end{picture}%
\setlength{\unitlength}{4144sp}%
\begingroup\makeatletter\ifx\SetFigFont\undefined%
\gdef\SetFigFont#1#2#3#4#5{%
  \reset@font\fontsize{#1}{#2pt}%
  \fontfamily{#3}\fontseries{#4}\fontshape{#5}%
  \selectfont}%
\fi\endgroup%
\begin{picture}(6414,1599)(79,-838)
\put(2656,254){\makebox(0,0)[b]{\smash{{\SetFigFont{10}{12.0}{\rmdefault}{\mddefault}{\updefault}$g$}}}}
\put(4006,524){\makebox(0,0)[b]{\smash{{\SetFigFont{10}{12.0}{\rmdefault}{\mddefault}{\updefault}$f'$}}}}
\put(4006,254){\makebox(0,0)[b]{\smash{{\SetFigFont{10}{12.0}{\rmdefault}{\mddefault}{\updefault}$g'$}}}}
\put(4366,-376){\makebox(0,0)[b]{\smash{{\SetFigFont{10}{12.0}{\rmdefault}{\mddefault}{\updefault}$f$}}}}
\put(5356,-376){\makebox(0,0)[b]{\smash{{\SetFigFont{10}{12.0}{\rmdefault}{\mddefault}{\updefault}$g$}}}}
\put(6346,-376){\makebox(0,0)[b]{\smash{{\SetFigFont{10}{12.0}{\rmdefault}{\mddefault}{\updefault}$h$}}}}
\put(1036,-646){\makebox(0,0)[b]{\smash{{\SetFigFont{10}{12.0}{\rmdefault}{\mddefault}{\updefault}$\alpha$}}}}
\put(1396,-646){\makebox(0,0)[b]{\smash{{\SetFigFont{10}{12.0}{\rmdefault}{\mddefault}{\updefault}$\beta$}}}}
\put(3286,119){\makebox(0,0)[b]{\smash{{\SetFigFont{10}{12.0}{\rmdefault}{\mddefault}{\updefault}$\alpha$}}}}
\put(3286,-151){\makebox(0,0)[b]{\smash{{\SetFigFont{10}{12.0}{\rmdefault}{\mddefault}{\updefault}$\beta$}}}}
\put(4816,-691){\makebox(0,0)[b]{\smash{{\SetFigFont{10}{12.0}{\rmdefault}{\mddefault}{\updefault}$\alpha$}}}}
\put(5851,-691){\makebox(0,0)[b]{\smash{{\SetFigFont{10}{12.0}{\rmdefault}{\mddefault}{\updefault}$\beta$}}}}
\put(226,-376){\makebox(0,0)[b]{\smash{{\SetFigFont{10}{12.0}{\rmdefault}{\mddefault}{\updefault}$f$}}}}
\put(586,-376){\makebox(0,0)[b]{\smash{{\SetFigFont{10}{12.0}{\rmdefault}{\mddefault}{\updefault}$g$}}}}
\put(1936,-376){\makebox(0,0)[b]{\smash{{\SetFigFont{10}{12.0}{\rmdefault}{\mddefault}{\updefault}$f'$}}}}
\put(2296,-376){\makebox(0,0)[b]{\smash{{\SetFigFont{10}{12.0}{\rmdefault}{\mddefault}{\updefault}$g'$}}}}
\put(2656,524){\makebox(0,0)[b]{\smash{{\SetFigFont{10}{12.0}{\rmdefault}{\mddefault}{\updefault}$f$}}}}
\end{picture}%
\end{center}

\vfill\pagebreak
\noindent All the constructions are identified modulo the following moves:
\begin{center}
\begin{picture}(0,0)%
\includegraphics{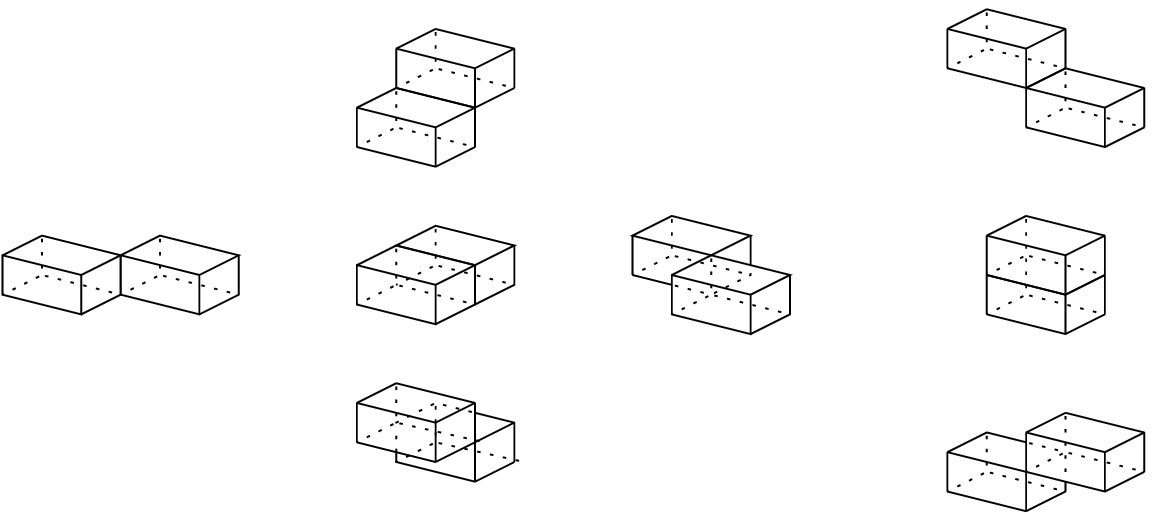}%
\end{picture}%
\setlength{\unitlength}{4144sp}%
\begingroup\makeatletter\ifx\SetFigFont\undefined%
\gdef\SetFigFont#1#2#3#4#5{%
  \reset@font\fontsize{#1}{#2pt}%
  \fontfamily{#3}\fontseries{#4}\fontshape{#5}%
  \selectfont}%
\fi\endgroup%
\begin{picture}(5244,2319)(169,-1513)
\put(2791,-466){\makebox(0,0)[b]{\smash{{\SetFigFont{10}{12.0}{\rmdefault}{\mddefault}{\updefault}{\color[rgb]{0,0,0}$\equiv$}%
}}}}
\put(2161,-106){\makebox(0,0)[b]{\smash{{\SetFigFont{10}{12.0}{\rmdefault}{\mddefault}{\updefault}{\color[rgb]{0,0,0}$\equiv$}%
}}}}
\put(2161,-826){\makebox(0,0)[b]{\smash{{\SetFigFont{10}{12.0}{\rmdefault}{\mddefault}{\updefault}{\color[rgb]{0,0,0}$\equiv$}%
}}}}
\put(4951,-61){\makebox(0,0)[b]{\smash{{\SetFigFont{10}{12.0}{\rmdefault}{\mddefault}{\updefault}{\color[rgb]{0,0,0}$\equiv$}%
}}}}
\put(4951,-916){\makebox(0,0)[b]{\smash{{\SetFigFont{10}{12.0}{\rmdefault}{\mddefault}{\updefault}{\color[rgb]{0,0,0}$\equiv$}%
}}}}
\put(1531,-466){\makebox(0,0)[b]{\smash{{\SetFigFont{10}{12.0}{\rmdefault}{\mddefault}{\updefault}{\color[rgb]{0,0,0}$\equiv$}%
}}}}
\end{picture}%
\end{center}

\medskip
\noindent This picture contains three families of moves, one for each exchange relation $\equi{02}$, $\equi{12}$ and $\equi{01}$. In the next paragraph, we prove that, in $\mon{\Sigma^K}$, the first two relations correspond to the two types of structural bureaucracy on SKS proofs.

As we have seen, the relation $\equi{01}$ is induced by the other two, together with the deformation relation on arrows of $G(\Sigma^K)$: thus it can be seen as a ghost exchange relation, generated by identification of circuits \emph{modulo} homeomorphic deformation.

\subsection{Structural bureaucracy is exchange}\label{sub:bureaucratie-est-echange}

\noindent Here, we prove that exchange relations are the polygraphic version of structural bureaucracy:

\begin{thm*}\label{thm:2}
For every two arrows $\alpha$ and $\beta$ in $G(\Sigma^K)$, the following two equations hold:
$$
\begin{aligned}
\pi\big((\alpha\star_0 s_2(\beta))\star_2(t_2(\alpha)\star_0\beta)\big) &\quad\equi{A}\quad \pi\big((s_2(\alpha)\star_0\beta)\star_2(\alpha\star_0 t_2(\beta))\big), \\
\pi\big((\alpha\star_1 s_2(\beta))\star_2(t_2(\alpha)\star_1\beta)\big) &\quad\equi{B}\quad \pi\big((s_2(\alpha)\star_1\beta)\star_2(\alpha\star_1 t_2(\beta))\big).
\end{aligned}
$$

\noindent Conversely, let us consider two SKS proofs $D:a\mred{\alpha}b\mred{\beta}c$ and $D':a\mred{\beta}b'\mred{\alpha}c$ such that $D\red{A}D'$ (resp. $D\red{B}D'$). Then there exist $2$-arrows $f$, $g$, $g'$ and $h$ in $\Sigma^K$ and arrows $\alpha_1$, $\alpha_2$, $\beta_1$ and $\beta_2$ in $G(\Sigma^K)$ such that the following conditions hold:
\begin{enumerate}
\item[-] The following two diagrams are paths in $G(\Sigma^K)$:
$$
f\red{\alpha_1}g\red{\beta_2}h \et f\red{\beta_1}g'\red{\alpha_2}h.
$$

\item[-] The following relation hold in $\mon{G(\Sigma^K)}$ with $i=0$ (resp. $i=1$):
$$
f\red{\alpha_1}g\red{\beta_2}h \quad\equi{i2}\quad f\red{\beta_1}g'\red{\alpha_2}h.
$$

\item[-] The following two equalities hold in $\mon{G^K}$:
$$
\pi(f\red{\alpha_1}g\red{\beta_2}h)\:=\:D \et \pi(f\red{\beta_1}g'\red{\alpha_2}h)\:=\:D'.
$$
\end{enumerate}
\end{thm*}

\begin{dem}
Let us fix two arrows $\alpha$ and $\beta$ in $G(\Sigma^K)$. By construction, these arrows are of the form $(f,\Phi(\alpha_0),g)$ and $(h,\Phi(\beta_0),k)$. Therefore, they form a diagram of the following shape in the reduction graph $G(\Sigma^K)$:

\begin{center}
\begin{picture}(0,0)%
\includegraphics{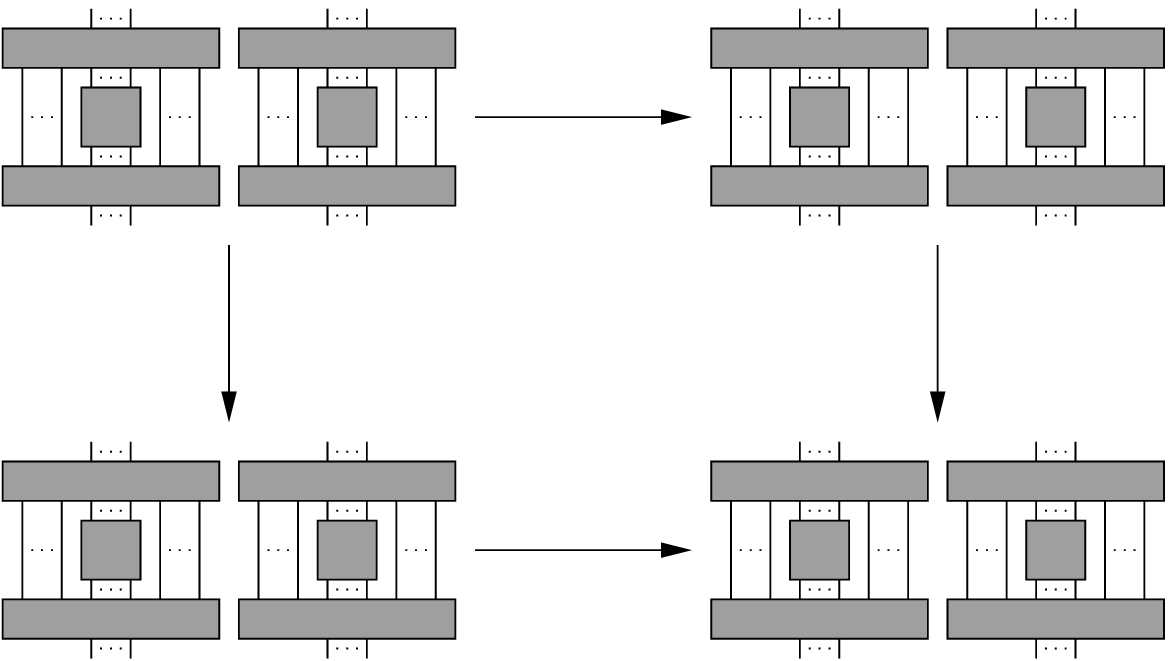}%
\end{picture}%
\setlength{\unitlength}{4144sp}%
\begingroup\makeatletter\ifx\SetFigFont\undefined%
\gdef\SetFigFont#1#2#3#4#5{%
  \reset@font\fontsize{#1}{#2pt}%
  \fontfamily{#3}\fontseries{#4}\fontshape{#5}%
  \selectfont}%
\fi\endgroup%
\begin{picture}(5334,2994)(79,-2233)
\put(4906,-1771){\makebox(0,0)[b]{\smash{{\SetFigFont{10}{12.0}{\rmdefault}{\mddefault}{\updefault}{\color[rgb]{0,0,0}$t'$}%
}}}}
\put(2701,389){\makebox(0,0)[b]{\smash{{\SetFigFont{10}{12.0}{\rmdefault}{\mddefault}{\updefault}{\color[rgb]{0,0,0}$\alpha\star_0 s_2(\beta)$}%
}}}}
\put(2701,-1951){\makebox(0,0)[b]{\smash{{\SetFigFont{10}{12.0}{\rmdefault}{\mddefault}{\updefault}{\color[rgb]{0,0,0}$\alpha\star_0 t_2(\beta)$}%
}}}}
\put(676,-736){\makebox(0,0)[b]{\smash{{\SetFigFont{10}{12.0}{\rmdefault}{\mddefault}{\updefault}{\color[rgb]{0,0,0}$s_2(\alpha)\star_0\beta$}%
}}}}
\put(4816,-736){\makebox(0,0)[b]{\smash{{\SetFigFont{10}{12.0}{\rmdefault}{\mddefault}{\updefault}{\color[rgb]{0,0,0}$t_2(\alpha)\star_0\beta$}%
}}}}
\put(3826,524){\makebox(0,0)[b]{\smash{{\SetFigFont{10}{12.0}{\rmdefault}{\mddefault}{\updefault}{\color[rgb]{0,0,0}$f$}%
}}}}
\put(3826,-106){\makebox(0,0)[b]{\smash{{\SetFigFont{10}{12.0}{\rmdefault}{\mddefault}{\updefault}{\color[rgb]{0,0,0}$g$}%
}}}}
\put(1666,-1456){\makebox(0,0)[b]{\smash{{\SetFigFont{10}{12.0}{\rmdefault}{\mddefault}{\updefault}{\color[rgb]{0,0,0}$h$}%
}}}}
\put(1666,-2086){\makebox(0,0)[b]{\smash{{\SetFigFont{10}{12.0}{\rmdefault}{\mddefault}{\updefault}{\color[rgb]{0,0,0}$k$}%
}}}}
\put(3826,-1456){\makebox(0,0)[b]{\smash{{\SetFigFont{10}{12.0}{\rmdefault}{\mddefault}{\updefault}{\color[rgb]{0,0,0}$f$}%
}}}}
\put(3826,-2086){\makebox(0,0)[b]{\smash{{\SetFigFont{10}{12.0}{\rmdefault}{\mddefault}{\updefault}{\color[rgb]{0,0,0}$g$}%
}}}}
\put(4906,-1456){\makebox(0,0)[b]{\smash{{\SetFigFont{10}{12.0}{\rmdefault}{\mddefault}{\updefault}{\color[rgb]{0,0,0}$h$}%
}}}}
\put(4906,-2086){\makebox(0,0)[b]{\smash{{\SetFigFont{10}{12.0}{\rmdefault}{\mddefault}{\updefault}{\color[rgb]{0,0,0}$k$}%
}}}}
\put(586,524){\makebox(0,0)[b]{\smash{{\SetFigFont{10}{12.0}{\rmdefault}{\mddefault}{\updefault}{\color[rgb]{0,0,0}$f$}%
}}}}
\put(586,-106){\makebox(0,0)[b]{\smash{{\SetFigFont{10}{12.0}{\rmdefault}{\mddefault}{\updefault}{\color[rgb]{0,0,0}$g$}%
}}}}
\put(1666,524){\makebox(0,0)[b]{\smash{{\SetFigFont{10}{12.0}{\rmdefault}{\mddefault}{\updefault}{\color[rgb]{0,0,0}$h$}%
}}}}
\put(1666,-106){\makebox(0,0)[b]{\smash{{\SetFigFont{10}{12.0}{\rmdefault}{\mddefault}{\updefault}{\color[rgb]{0,0,0}$k$}%
}}}}
\put(4906,524){\makebox(0,0)[b]{\smash{{\SetFigFont{10}{12.0}{\rmdefault}{\mddefault}{\updefault}{\color[rgb]{0,0,0}$h$}%
}}}}
\put(4906,-106){\makebox(0,0)[b]{\smash{{\SetFigFont{10}{12.0}{\rmdefault}{\mddefault}{\updefault}{\color[rgb]{0,0,0}$k$}%
}}}}
\put(586,-1456){\makebox(0,0)[b]{\smash{{\SetFigFont{10}{12.0}{\rmdefault}{\mddefault}{\updefault}{\color[rgb]{0,0,0}$f$}%
}}}}
\put(586,-2086){\makebox(0,0)[b]{\smash{{\SetFigFont{10}{12.0}{\rmdefault}{\mddefault}{\updefault}{\color[rgb]{0,0,0}$g$}%
}}}}
\put(586,209){\makebox(0,0)[b]{\smash{{\SetFigFont{10}{12.0}{\rmdefault}{\mddefault}{\updefault}{\color[rgb]{0,0,0}$s$}%
}}}}
\put(1666,209){\makebox(0,0)[b]{\smash{{\SetFigFont{10}{12.0}{\rmdefault}{\mddefault}{\updefault}{\color[rgb]{0,0,0}$s'$}%
}}}}
\put(3826,209){\makebox(0,0)[b]{\smash{{\SetFigFont{10}{12.0}{\rmdefault}{\mddefault}{\updefault}{\color[rgb]{0,0,0}$t$}%
}}}}
\put(1666,-1771){\makebox(0,0)[b]{\smash{{\SetFigFont{10}{12.0}{\rmdefault}{\mddefault}{\updefault}{\color[rgb]{0,0,0}$t'$}%
}}}}
\put(586,-1771){\makebox(0,0)[b]{\smash{{\SetFigFont{10}{12.0}{\rmdefault}{\mddefault}{\updefault}{\color[rgb]{0,0,0}$s$}%
}}}}
\put(4906,209){\makebox(0,0)[b]{\smash{{\SetFigFont{10}{12.0}{\rmdefault}{\mddefault}{\updefault}{\color[rgb]{0,0,0}$s'$}%
}}}}
\put(3826,-1771){\makebox(0,0)[b]{\smash{{\SetFigFont{10}{12.0}{\rmdefault}{\mddefault}{\updefault}{\color[rgb]{0,0,0}$t$}%
}}}}
\end{picture}%
\end{center}

\bigskip
\noindent where $s=s_2(\Phi(\alpha_0))$, $s'=s_2(\Phi(\beta_0))$, $t=t_2(\Phi(\alpha_0))$ and $t'=t_2(\Phi(\beta_0))$. When $\pi$ is applied to this diagram, we get two SKS proofs that satisfy the definition of $\red{A}$. The proof of the second equality is handled similarly.

Conversely, let us consider two SKS proofs $D$ and $D'$ such that $D\red{A}D'$. Let us consider the $1$-arrows $Y_1$, $Y_2$ and~$Y_3$ and the families of terms $u$ and $v$ such that these two SKS proofs form the following diagram:

$$
\xymatrix{u\circ(Y_1\tens s(\alpha)\tens Y_2\tens s(\beta)\tens Y_3)\circ v \ar@{>>}+<3.2cm,0cm>;[rrr]-<3.2cm,0cm>^-{\alpha} \ar@{>>}-<0cm,0.5cm>;[dd]+<0cm,0.5cm>_-{\beta} &&& u\circ(Y_1\tens t(\alpha)\tens Y_2\tens s(\beta)\tens Y_3)\circ v \ar@{>>}-<0cm,0.5cm>;[dd]+<0cm,0.5cm>^-{\beta} \ar-<3.5cm,1cm>;[ddlll]+<3.5cm,1cm>^-{A} \\ \\ u\circ(Y_1\tens s(\alpha)\tens Y_2\tens t(\beta)\tens Y_3)\circ v \ar@{>>}+<3.2cm,0cm>;[rrr]-<3.2cm,0cm>_-{\alpha} &&& u\circ(Y_1\tens t(\alpha)\tens Y_2\tens t(\beta)\tens Y_3)\circ v,}
$$

\medskip
\noindent with $a$, $b$, $b'$ and $c$ being the terms at the corners, from left to right and top to bottom. Then, from the decomposition of $a$, we know that there exists a $2$-arrow $f$ in $\Sigma^K$ such that $\Phi^X(a)\equi{\Delta}f$ and $f$ has the following shape:
\begin{center}
\begin{picture}(0,0)%
\includegraphics{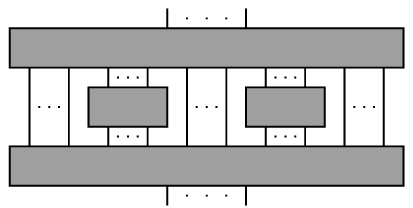}%
\end{picture}%
\setlength{\unitlength}{4144sp}%
\begingroup\makeatletter\ifx\SetFigFont\undefined%
\gdef\SetFigFont#1#2#3#4#5{%
  \reset@font\fontsize{#1}{#2pt}%
  \fontfamily{#3}\fontseries{#4}\fontshape{#5}%
  \selectfont}%
\fi\endgroup%
\begin{picture}(2145,924)(-242,-163)
\put(-89,254){\makebox(0,0)[b]{\smash{{\SetFigFont{10}{12.0}{\rmdefault}{\mddefault}{\updefault}{\color[rgb]{0,0,0}$f=$}%
}}}}
\put(631,254){\makebox(0,0)[b]{\smash{{\SetFigFont{10}{12.0}{\rmdefault}{\mddefault}{\updefault}{\color[rgb]{0,0,0}$s$}%
}}}}
\put(1351,254){\makebox(0,0)[b]{\smash{{\SetFigFont{10}{12.0}{\rmdefault}{\mddefault}{\updefault}{\color[rgb]{0,0,0}$s'$}%
}}}}
\put(991,-16){\makebox(0,0)[b]{\smash{{\SetFigFont{10}{12.0}{\rmdefault}{\mddefault}{\updefault}{\color[rgb]{0,0,0}$k'$}%
}}}}
\put(991,524){\makebox(0,0)[b]{\smash{{\SetFigFont{10}{12.0}{\rmdefault}{\mddefault}{\updefault}{\color[rgb]{0,0,0}$k$}%
}}}}
\end{picture}%
\end{center}

\noindent Then, the following diagram represents two paths in $G(\Sigma^K)$ which are equivalent \emph{modulo} $\equi{02}$, with $s=s_2(\Phi(\alpha))$, $s'=s_2(\Phi(\beta))$, $t=t_2(\Phi(\alpha))$, $t'=t_2(\Phi(\beta))$ and the arrows $\alpha_1$, $\alpha_2$, $\beta_1$ and $\beta_2$ defined implicitely:

\begin{center}
\begin{picture}(0,0)%
\includegraphics{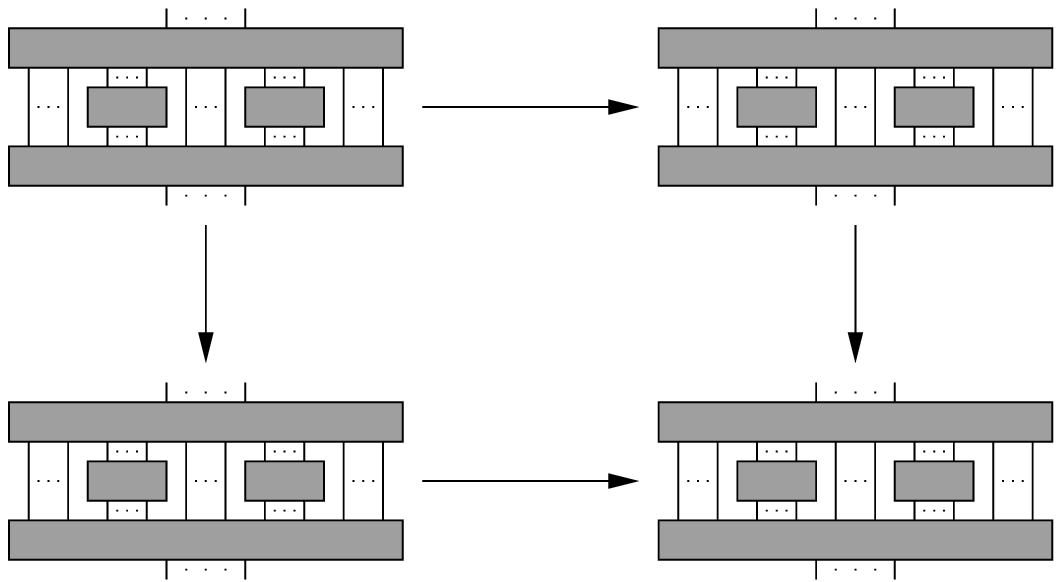}%
\end{picture}%
\setlength{\unitlength}{4144sp}%
\begingroup\makeatletter\ifx\SetFigFont\undefined%
\gdef\SetFigFont#1#2#3#4#5{%
  \reset@font\fontsize{#1}{#2pt}%
  \fontfamily{#3}\fontseries{#4}\fontshape{#5}%
  \selectfont}%
\fi\endgroup%
\begin{picture}(5467,2634)(-265,-1873)
\put(3601,-1456){\makebox(0,0)[b]{\smash{{\SetFigFont{10}{12.0}{\rmdefault}{\mddefault}{\updefault}{\color[rgb]{0,0,0}$t$}%
}}}}
\put(631,254){\makebox(0,0)[b]{\smash{{\SetFigFont{10}{12.0}{\rmdefault}{\mddefault}{\updefault}{\color[rgb]{0,0,0}$s$}%
}}}}
\put(1351,254){\makebox(0,0)[b]{\smash{{\SetFigFont{10}{12.0}{\rmdefault}{\mddefault}{\updefault}{\color[rgb]{0,0,0}$s'$}%
}}}}
\put(991,-16){\makebox(0,0)[b]{\smash{{\SetFigFont{10}{12.0}{\rmdefault}{\mddefault}{\updefault}{\color[rgb]{0,0,0}$k'$}%
}}}}
\put(991,524){\makebox(0,0)[b]{\smash{{\SetFigFont{10}{12.0}{\rmdefault}{\mddefault}{\updefault}{\color[rgb]{0,0,0}$k$}%
}}}}
\put(-89,254){\makebox(0,0)[b]{\smash{{\SetFigFont{10}{12.0}{\rmdefault}{\mddefault}{\updefault}{\color[rgb]{0,0,0}$f=$}%
}}}}
\put(2431,389){\makebox(0,0)[b]{\smash{{\SetFigFont{10}{12.0}{\rmdefault}{\mddefault}{\updefault}{\color[rgb]{0,0,0}$\alpha_1$}%
}}}}
\put(-89,-1456){\makebox(0,0)[b]{\smash{{\SetFigFont{10}{12.0}{\rmdefault}{\mddefault}{\updefault}{\color[rgb]{0,0,0}$g'=$}%
}}}}
\put(5041,254){\makebox(0,0)[b]{\smash{{\SetFigFont{10}{12.0}{\rmdefault}{\mddefault}{\updefault}{\color[rgb]{0,0,0}$=g$}%
}}}}
\put(5041,-1456){\makebox(0,0)[b]{\smash{{\SetFigFont{10}{12.0}{\rmdefault}{\mddefault}{\updefault}{\color[rgb]{0,0,0}$=h$}%
}}}}
\put(4141,-601){\makebox(0,0)[b]{\smash{{\SetFigFont{10}{12.0}{\rmdefault}{\mddefault}{\updefault}{\color[rgb]{0,0,0}$\beta_2$}%
}}}}
\put(811,-601){\makebox(0,0)[b]{\smash{{\SetFigFont{10}{12.0}{\rmdefault}{\mddefault}{\updefault}{\color[rgb]{0,0,0}$\beta_1$}%
}}}}
\put(2431,-1636){\makebox(0,0)[b]{\smash{{\SetFigFont{10}{12.0}{\rmdefault}{\mddefault}{\updefault}{\color[rgb]{0,0,0}$\alpha_2$}%
}}}}
\put(991,-1726){\makebox(0,0)[b]{\smash{{\SetFigFont{10}{12.0}{\rmdefault}{\mddefault}{\updefault}{\color[rgb]{0,0,0}$k'$}%
}}}}
\put(991,-1186){\makebox(0,0)[b]{\smash{{\SetFigFont{10}{12.0}{\rmdefault}{\mddefault}{\updefault}{\color[rgb]{0,0,0}$k$}%
}}}}
\put(3961,524){\makebox(0,0)[b]{\smash{{\SetFigFont{10}{12.0}{\rmdefault}{\mddefault}{\updefault}{\color[rgb]{0,0,0}$k$}%
}}}}
\put(3961,-16){\makebox(0,0)[b]{\smash{{\SetFigFont{10}{12.0}{\rmdefault}{\mddefault}{\updefault}{\color[rgb]{0,0,0}$k'$}%
}}}}
\put(3961,-1726){\makebox(0,0)[b]{\smash{{\SetFigFont{10}{12.0}{\rmdefault}{\mddefault}{\updefault}{\color[rgb]{0,0,0}$k'$}%
}}}}
\put(3961,-1186){\makebox(0,0)[b]{\smash{{\SetFigFont{10}{12.0}{\rmdefault}{\mddefault}{\updefault}{\color[rgb]{0,0,0}$k$}%
}}}}
\put(631,-1456){\makebox(0,0)[b]{\smash{{\SetFigFont{10}{12.0}{\rmdefault}{\mddefault}{\updefault}{\color[rgb]{0,0,0}$s$}%
}}}}
\put(4321,254){\makebox(0,0)[b]{\smash{{\SetFigFont{10}{12.0}{\rmdefault}{\mddefault}{\updefault}{\color[rgb]{0,0,0}$s'$}%
}}}}
\put(3601,254){\makebox(0,0)[b]{\smash{{\SetFigFont{10}{12.0}{\rmdefault}{\mddefault}{\updefault}{\color[rgb]{0,0,0}$t$}%
}}}}
\put(1351,-1456){\makebox(0,0)[b]{\smash{{\SetFigFont{10}{12.0}{\rmdefault}{\mddefault}{\updefault}{\color[rgb]{0,0,0}$t'$}%
}}}}
\put(4321,-1456){\makebox(0,0)[b]{\smash{{\SetFigFont{10}{12.0}{\rmdefault}{\mddefault}{\updefault}{\color[rgb]{0,0,0}$t'$}%
}}}}
\end{picture}%
\end{center}

\bigskip
\noindent These two paths satisfy the relation $\equi{02}$. Furthermore, from the decompositions of $a$, $b$, $b'$ and $c$ on one hand, from the ones of $f$, $g$, $g'$ and $h$ on the other hand, we have:
$$
f\equi{\Delta}\Phi^X(a), \quad g\equi{\Delta}\Phi^X(b), \quad g'\equi{\Delta}\Phi^X(b'), \quad h\equi{\Delta}\Phi^X(c).
$$

\medskip
\noindent Finally, it is straightforward to check that $\pi(\alpha_1\star_2\beta_2)=D$ and $\pi(\beta_1\star_2\alpha_2)=D'$. The proof in the case $D\red{B}D'$ follows the same scheme.

\findem\end{dem}

\subsection{Some geography}\label{sub:geographie}

\noindent The formalism SKS had two main offsprings, called Formalism A and Formalism B [Guglielmi 2005]: in the former, proofs are identified \emph{modulo} bureaucracy A and, in the latter, \emph{modulo} both types A and B. However, there is not much freedom in the construction of more formalisms, depending on what proofs one wants to identify.

But, in the higher-dimensional setting, there are $3$-polygraphs corresponding to each of these three formalisms, among many others that are linked by the following \emph{categorical map} - a diagram in the category of families of $3$-polygraphs over the $2$-polygraph $\Sigma^F$ of SKS formulas:
$$
\xymatrixnocompile{ \Sigma^K \ar@{ >->}[r] \ar@{ >->}[dr] & G(\Sigma^K) \ar@{ >->}[r] \ar@{ >->}[d] & \mon{G(\Sigma^K)} \ar@{ >->}[d] \\ & \Gr(\Sigma^K) \ar@{ >->}[r] & \mon{\Gr(\Sigma^K)} \ar@{->>}[r] \ar@{->>}[d] & \mon{\Gr(\Sigma^K)}/\equi{02} \ar@{->>}[d] \\ && \mon{\Gr(\Sigma^K)}/\equi{12} \ar@{->>}[r] & \mon{\Sigma^K}}
$$

\noindent Let us give a description of all these objects. One starts with the $3$-polygraph $\Sigma^K$: its $3$-cells are the inference rules, the structural rules and the resources management rules. From this object, one can consider all the rules, applied in any context, which yields the reduction graph $G(\Sigma^K)$: its arrows are all the one-step sequential reductions. Alternatively, one can consider all the rules applied in any existing context and possibly in parallel to build $\Gr(\Sigma^K)$, a graph which arrows are the one-step parallel reductions.

Then, one considers the paths generated by $G(\Sigma^K)$: this produces the free-category $\mon{G(\Sigma^K)}$ which arrows correspond to SKS proofs. This is the polygraphic equivalent of the calculus of structures version of SKS. Alternatively, the paths generated by $\Gr(\Sigma^K)$ give the free category $\mon{\Gr(\Sigma^K)}$. There, arrows correspond to SKS proofs generalized with the possible application in parallel of inference rules. Here, bureaucracy is at its highest level, since all the described proofs differing by the order of application of subproofs are distinguished; furthermore, there is at each time a third possible proof, consisting in the simultaneous application of both subproofs. Hence, this is the biggest object of this classification.

There one starts the quotients of $\mon{\Gr(\Sigma^K)}$ by the exchange relations. The first possibility is to quotient by the first family of exchange relations, corresponding to bureaucracy type A. This yields the object $\mon{\Gr(\Sigma^K)}/\equi{02}$, which is the polygraphic version of Formalism A. As an alternative, one can instead quotient $\mon{\Gr(\Sigma^K)}$ by the exchange relations corresponding to bureaucracy type B, to get $\mon{\Gr(\Sigma^K)}/\equi{12}$, which has no equivalent in SKS derived formalisms. Finally, doing both quotients, one gets the free $3$-category $\mon{\Sigma^K}$ generated by $\Sigma^K$, where all the bureaucracy is killed. This is the polygraphic equivalent of Formalism B.

Hence, this diagram localizes the polygraphic equivalents of the known formalisms: $\Sigma^K$ for the signature of SKS, $\mon{G(\Sigma^K)}$ for the calculus of structures version of SKS, $\mon{\Gr(\Sigma^K)}/\equi{02}$ for Formalism~A and $\mon{\Sigma^K}$ for Formalism~B. But the diagram also encompasses still unknown formalisms that could prove to be useful, like the biggest one $\mon{\Gr(\Sigma^K)}$, where parallel applications of rules are allowed and distinguished from sequential ones, or $\mon{\Gr(\Sigma^K)}/\equi{12}$, where only bureaucracy B is killed. This is an example of the freedom the higher-dimensional setting lets to the user in the exact design of the (equivalence classes of) proofs he wants to consider. Another example of freedom is given in section \ref{sec:4d} about the possibilities offered for handling the equations between formulas.

\section{Representing proofs in three dimensions}\label{sec:dessins-3d}

\noindent This section is a first attempt at representing proofs in $3$ dimensions, so that one can view them as the genuine $3$-dimensional objects they are.

In order to represent $2$-arrows, Penrose diagrams are really convenient; they make $2$-arrows appear as circuits, using the following scheme: each $2$-cell is pictured as a vertice in a graph, each $1$-cell as an edge and each $0$-cell as a part of the plane which boundaries are the edges of the graph. Thus, each $k$-cell is pictured as a $(2-k)$-dimensional object. Then, the produced vertices and edges are thickened until they are $2$-dimensional; note that in the circuit representation, wires are not thickened to make drawing easier, but they should be for sake of coherence.

The application of a similar process to a $3$-dimensional arrow gives that each $k$-cell is represented as a $(3-k)$-dimensional object. In details: each $3$-cell is pictured as a point; each $2$-dimensional cell is a line (either open or between two points); each $1$-dimensional cell is a surface (either open or with a line as a boundary); each $0$-dimensional cell is a volume lying between surfaces. Finally, every object is thickened, if necessary, until it gets $3$-dimensional.

\noindent Let us draw a $3$-dimensional proof, using \emph{$\mathit{3}$-dimensional Penrose diagrams}. First, let us make a Penrose diagram for the following rewriting-style rule:
\begin{center}
\includegraphics{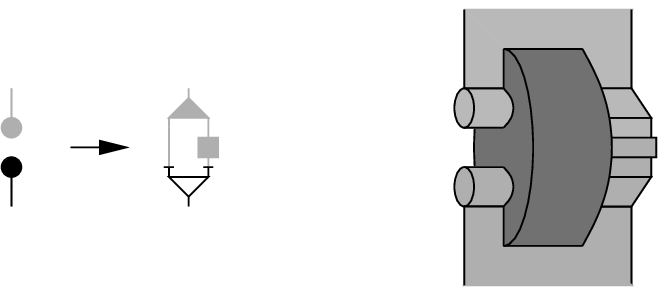}
\end{center}

\noindent In order to make pictures simpler, we do not distinguish the two sorts $A$ and $F$ anymore, the $2$-cells~$\wedge$ and~$\vee$ are drawn the same way and the $2$-cell $\iota$ disappears: these are only temporary choices, until we get easier ways to draw $3$-dimensional Penrose diagrams. When each $3$-cell has been given a $3$-dimensional representation, proofs can be drawn as pastings of these $3$-dimensional blocks, such as the following one:

\begin{center}
\includegraphics{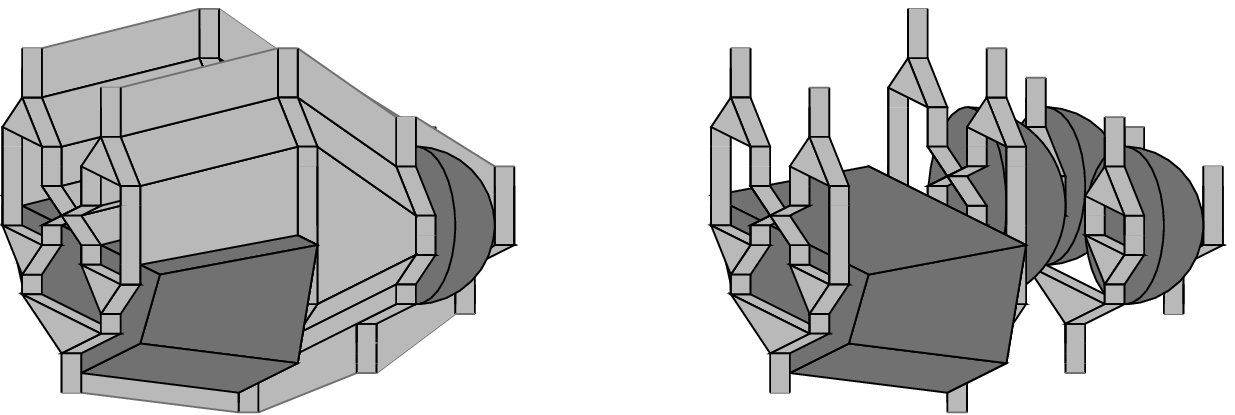}
\end{center}

\noindent Except for the aforementionned simplifications, the left-hand side picture is an accurate $3$-dimensional representation of a proof that the implication $(a\wedge b)\vee(a\wedge b)\:\implies\:a\wedge b$ holds for any atoms $a$ and~$b$. In the right-hand side picture, surfaces corresponding to $3$-dimensional identities have been removed in order to see internal parts of the proof. For a better understanding of how this object is built (and what lies behind some opaque volumes), one can make vertical slices of this object, to produce the following rewriting-style proof:
\begin{center}
\begin{picture}(0,0)%
\includegraphics{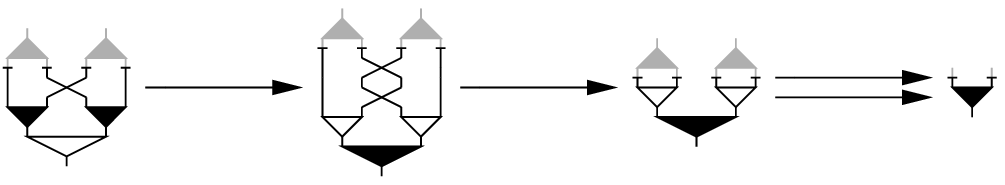}%
\end{picture}%
\setlength{\unitlength}{4144sp}%
\begingroup\makeatletter\ifx\SetFigFont\undefined%
\gdef\SetFigFont#1#2#3#4#5{%
  \reset@font\fontsize{#1}{#2pt}%
  \fontfamily{#3}\fontseries{#4}\fontshape{#5}%
  \selectfont}%
\fi\endgroup%
\begin{picture}(4568,1047)(462,-376)
\put(766,-331){\makebox(0,0)[b]{\smash{{\SetFigFont{10}{12.0}{\rmdefault}{\mddefault}{\updefault}{\color[rgb]{0,0,0}$(a\wedge b)\vee(a\wedge b)$}%
}}}}
\put(2206,-331){\makebox(0,0)[b]{\smash{{\SetFigFont{10}{12.0}{\rmdefault}{\mddefault}{\updefault}{\color[rgb]{0,0,0}$(a \vee a)\wedge(b\vee b)$}%
}}}}
\put(3646,-331){\makebox(0,0)[b]{\smash{{\SetFigFont{10}{12.0}{\rmdefault}{\mddefault}{\updefault}{\color[rgb]{0,0,0}$(a \vee a)\wedge(b\vee b)$}%
}}}}
\put(4906,-331){\makebox(0,0)[b]{\smash{{\SetFigFont{10}{12.0}{\rmdefault}{\mddefault}{\updefault}{\color[rgb]{0,0,0}$a\wedge b$}%
}}}}
\end{picture}%
\end{center}

\noindent Since the given representation uses only a fake third dimension, one could prefer to use a software dedicated to $3$-dimensional pictures. This has many advantages, such as being able to turn around the object and make snapshots from different points of view. For example, the following views of the same proof were generated using the software POV-Ray, a ray-tracer, freely available on http://www.povray.org. \\

\begin{center}
\hspace*{\stretch{1}} \includegraphics[scale=0.4]{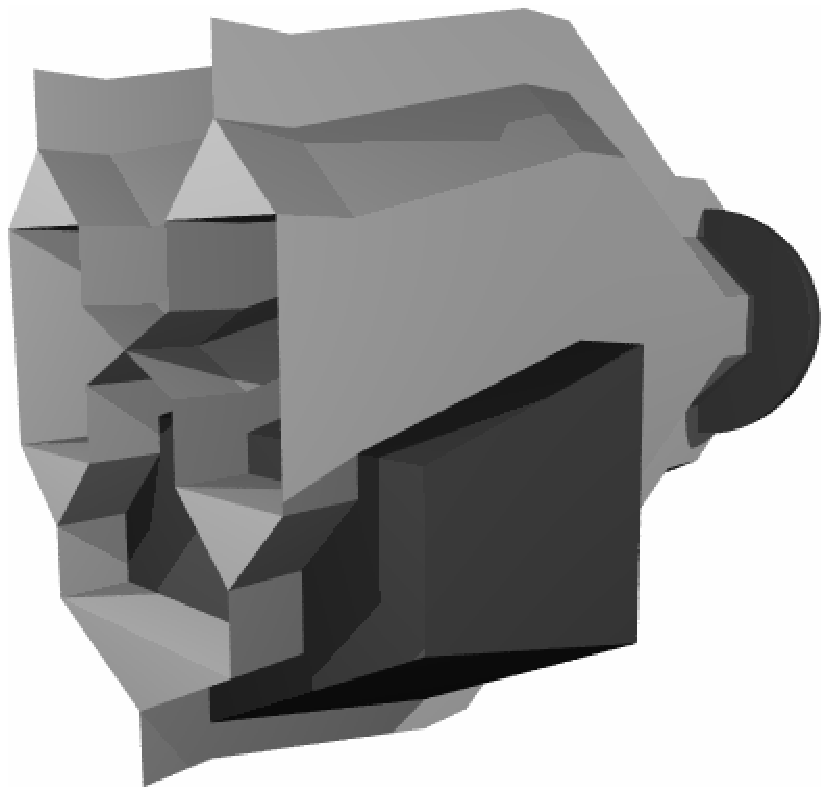} \hspace*{\stretch{1}} \includegraphics[scale=0.4]{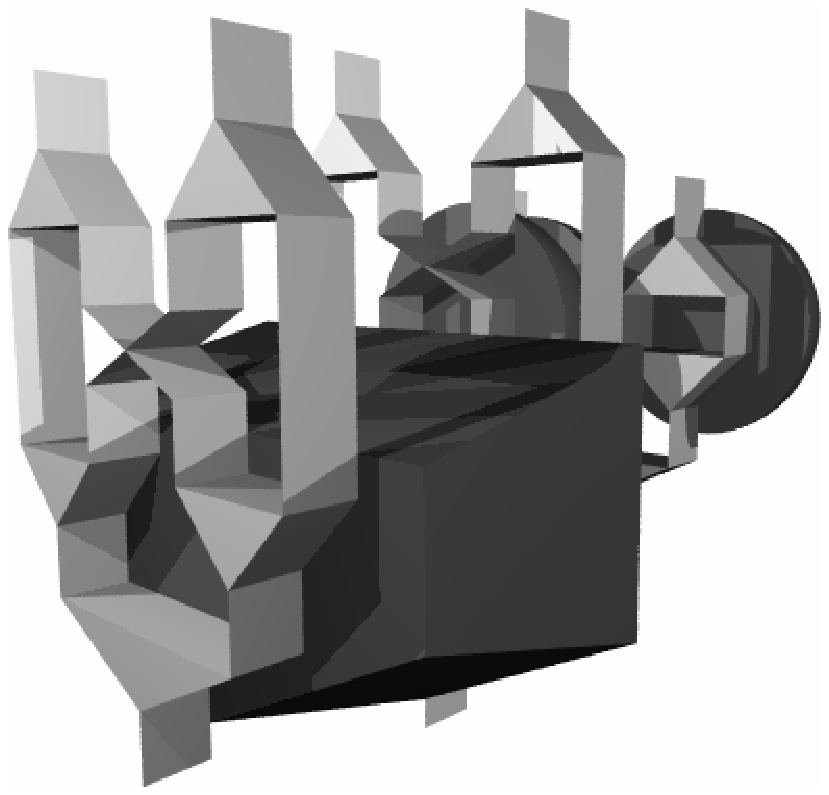} \hspace*{\stretch{1}} \strut
\end{center}

\begin{center}
\hspace*{\stretch{1}} \includegraphics[scale=0.4]{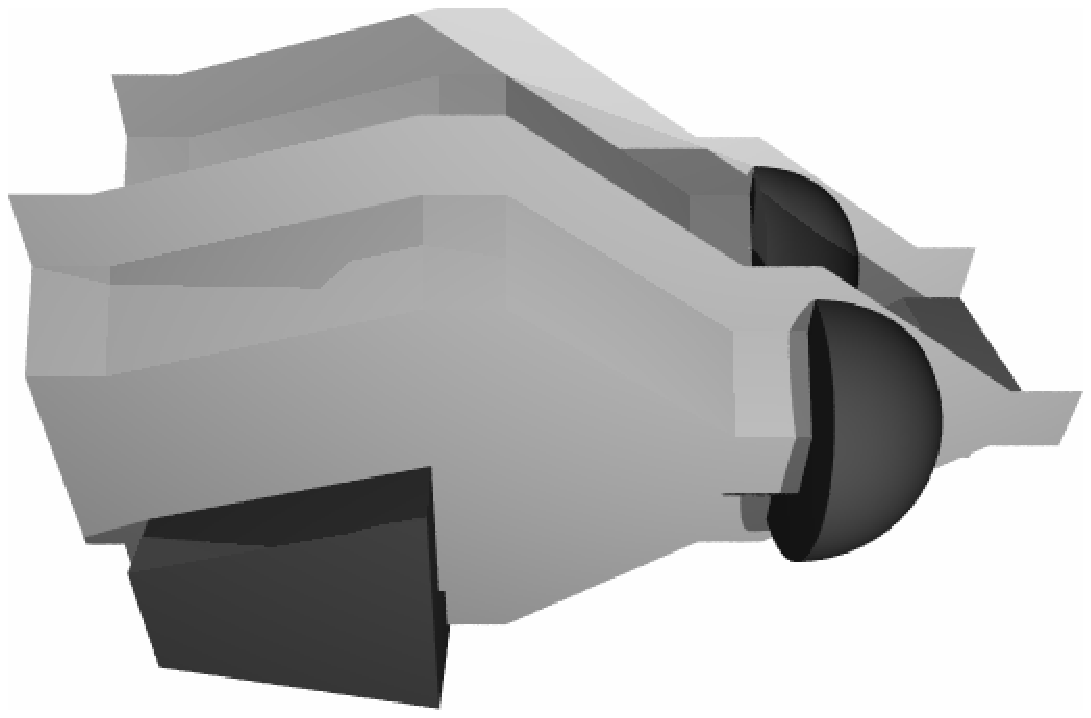} \hspace*{\stretch{1}} \includegraphics[scale=0.4]{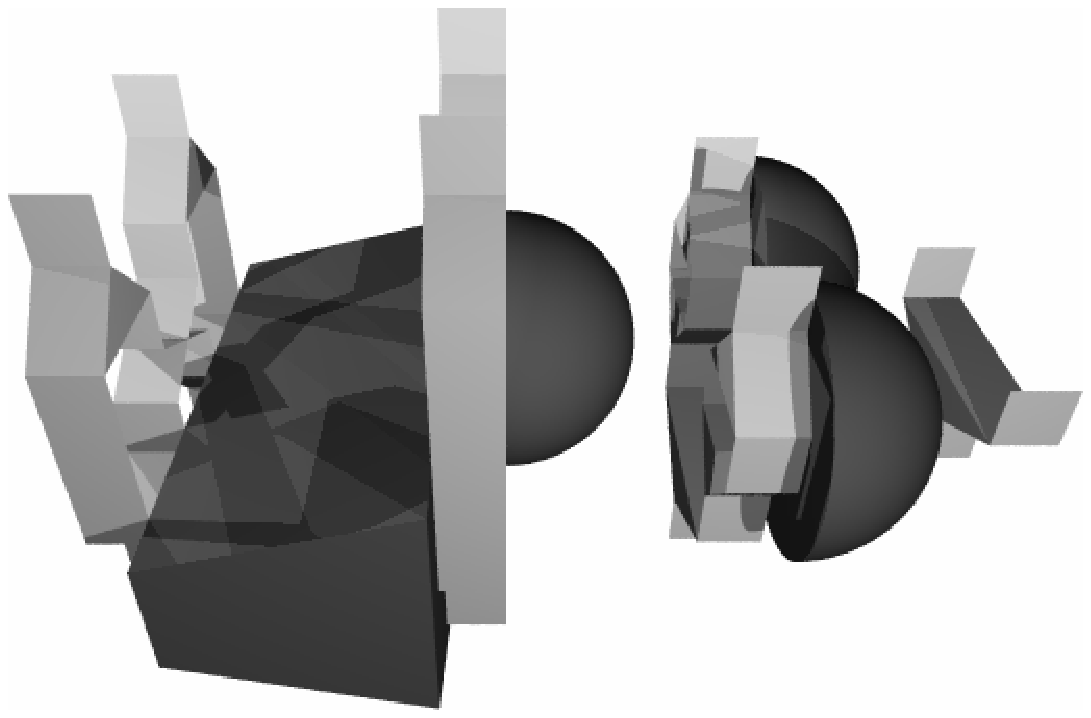} \hspace*{\stretch{1}} \strut
\end{center}

\bigskip\bigskip
\noindent Once again, the left-hand side pictures are the correct ones, while the right-hand side ones show internal parts of the proof. This part is quite new and some work will be necessary to easily produce nicer, more usable representations, so that the third dimension can provide more insight on what kind of objects proofs are.

\section{Normalization of proofs}\label{sec:4d}

When the third dimension gets involved, one can ask whether this dimensional increase will stop. The answer is quite simple: no. Indeed there are, at least, two good reasons to proceed to the fourth dimension.

The first one is total abstract nonsense - which does not mean that it is not a good reason. In category theory, there is a proverb saying: \emph{when one wants to study some objects, one should rather study their morphisms}. In higher-dimensional rewriting, there is something similar: \emph{when one wants to study some objects modulo some congruence relation, one should rather replace equations by rewriting rules} (this operation is called categorification in [Baez Dolan 1998]).

The second, more concrete reason is that there are two kinds of examples that give rise to $4$-dimensional arrows: equations between formulas and local transformations on proofs. This section is about a short glance at these two issues.

\subsection{Equations between formulas}\label{sub:equations-4d}

\noindent Previously, structural equations between formulas have been treated as pairs of inverse rules. But this is just one possibility, the higher-order rewriting framework allowing one to choose between many possible considerations. Here are three of them, but one can at least take any desired combination of them. \\

\noindent\textbf{Equations are equations.} The first possibility is, as stated before, to translate equations between formulas into equations between circuits. In that case, one considers circuits \emph{modulo} two families of equations. The first one is a faithful translation of the equations on formulas, so that, for example, one can recognize associativity of $\wedge$ and $\vee$ among them:
\begin{center}
\begin{picture}(0,0)%
\includegraphics{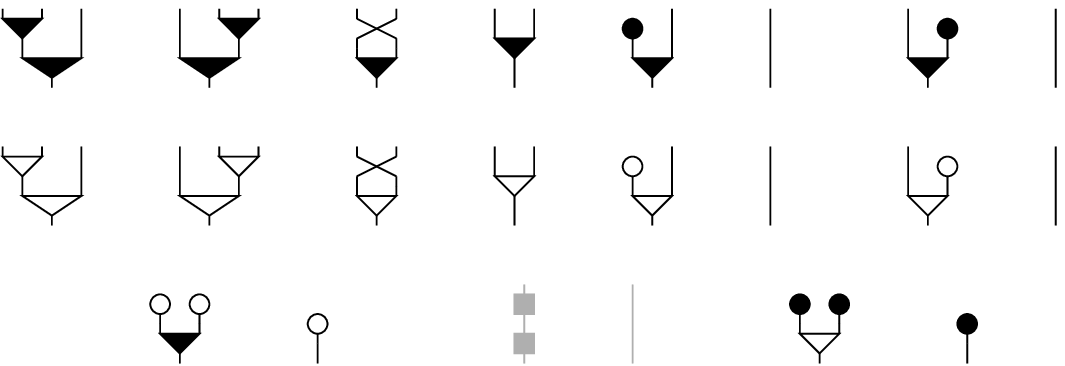}%
\end{picture}%
\setlength{\unitlength}{4144sp}%
\begingroup\makeatletter\ifx\SetFigFont\undefined%
\gdef\SetFigFont#1#2#3#4#5{%
  \reset@font\fontsize{#1}{#2pt}%
  \fontfamily{#3}\fontseries{#4}\fontshape{#5}%
  \selectfont}%
\fi\endgroup%
\begin{picture}(4839,1644)(79,-883)
\put(4186,-781){\makebox(0,0)[b]{\smash{{\SetFigFont{12}{14.4}{\rmdefault}{\mddefault}{\updefault}{\color[rgb]{0,0,0}$\equiv$}%
}}}}
\put(676,479){\makebox(0,0)[b]{\smash{{\SetFigFont{12}{14.4}{\rmdefault}{\mddefault}{\updefault}{\color[rgb]{0,0,0}$\equiv$}%
}}}}
\put(2116,479){\makebox(0,0)[b]{\smash{{\SetFigFont{12}{14.4}{\rmdefault}{\mddefault}{\updefault}{\color[rgb]{0,0,0}$\equiv$}%
}}}}
\put(3376,479){\makebox(0,0)[b]{\smash{{\SetFigFont{12}{14.4}{\rmdefault}{\mddefault}{\updefault}{\color[rgb]{0,0,0}$\equiv$}%
}}}}
\put(4681,479){\makebox(0,0)[b]{\smash{{\SetFigFont{12}{14.4}{\rmdefault}{\mddefault}{\updefault}{\color[rgb]{0,0,0}$\equiv$}%
}}}}
\put(676,-151){\makebox(0,0)[b]{\smash{{\SetFigFont{12}{14.4}{\rmdefault}{\mddefault}{\updefault}{\color[rgb]{0,0,0}$\equiv$}%
}}}}
\put(2116,-151){\makebox(0,0)[b]{\smash{{\SetFigFont{12}{14.4}{\rmdefault}{\mddefault}{\updefault}{\color[rgb]{0,0,0}$\equiv$}%
}}}}
\put(3376,-151){\makebox(0,0)[b]{\smash{{\SetFigFont{12}{14.4}{\rmdefault}{\mddefault}{\updefault}{\color[rgb]{0,0,0}$\equiv$}%
}}}}
\put(4681,-151){\makebox(0,0)[b]{\smash{{\SetFigFont{12}{14.4}{\rmdefault}{\mddefault}{\updefault}{\color[rgb]{0,0,0}$\equiv$}%
}}}}
\put(2746,-781){\makebox(0,0)[b]{\smash{{\SetFigFont{12}{14.4}{\rmdefault}{\mddefault}{\updefault}{\color[rgb]{0,0,0}$\equiv$}%
}}}}
\put(1261,-781){\makebox(0,0)[b]{\smash{{\SetFigFont{12}{14.4}{\rmdefault}{\mddefault}{\updefault}{\color[rgb]{0,0,0}$\equiv$}%
}}}}
\end{picture}%
\end{center}

\noindent The second family purpose is to give the resource management operators their real meaning, so that, for example, $\delta_A$ really is a local duplicator of atoms; among others, one gets the following equations:
\begin{center}
\begin{picture}(0,0)%
\includegraphics{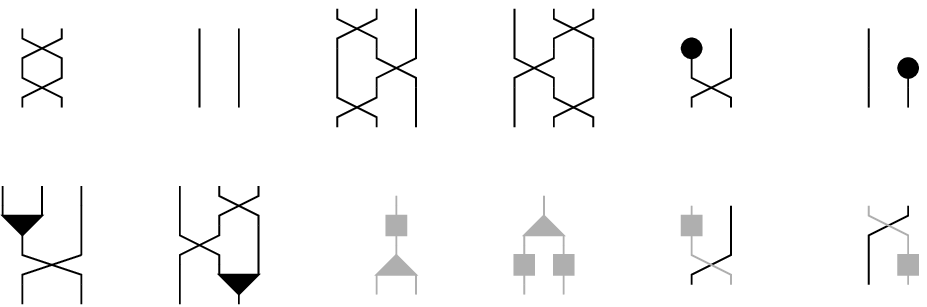}%
\end{picture}%
\setlength{\unitlength}{4144sp}%
\begingroup\makeatletter\ifx\SetFigFont\undefined%
\gdef\SetFigFont#1#2#3#4#5{%
  \reset@font\fontsize{#1}{#2pt}%
  \fontfamily{#3}\fontseries{#4}\fontshape{#5}%
  \selectfont}%
\fi\endgroup%
\begin{picture}(4209,1374)(79,-613)
\put(2206,-421){\makebox(0,0)[b]{\smash{{\SetFigFont{12}{14.4}{\rmdefault}{\mddefault}{\updefault}{\color[rgb]{0,0,0}$\equiv$}%
}}}}
\put(3736,389){\makebox(0,0)[b]{\smash{{\SetFigFont{12}{14.4}{\rmdefault}{\mddefault}{\updefault}{\color[rgb]{0,0,0}$\equiv$}%
}}}}
\put(676,389){\makebox(0,0)[b]{\smash{{\SetFigFont{12}{14.4}{\rmdefault}{\mddefault}{\updefault}{\color[rgb]{0,0,0}$\equiv$}%
}}}}
\put(2206,389){\makebox(0,0)[b]{\smash{{\SetFigFont{12}{14.4}{\rmdefault}{\mddefault}{\updefault}{\color[rgb]{0,0,0}$\equiv$}%
}}}}
\put(676,-421){\makebox(0,0)[b]{\smash{{\SetFigFont{12}{14.4}{\rmdefault}{\mddefault}{\updefault}{\color[rgb]{0,0,0}$\equiv$}%
}}}}
\put(3736,-421){\makebox(0,0)[b]{\smash{{\SetFigFont{12}{14.4}{\rmdefault}{\mddefault}{\updefault}{\color[rgb]{0,0,0}$\equiv$}%
}}}}
\end{picture}%
\end{center}

\medskip
\noindent\textbf{From equations to $\mathbf{3}$-dimensional isomorphisms.} Rather than considering equations on formulas as equations on circuits, one can treat them as invertible computations. Indeed, equations are often clashing with computational considerations, so that, whenever possible, they are replaced by local computations. Hence, one could replace the two aforementionned families of equations by two families of invertible $3$-cells. For example, the equation enforcing the associativity of $\wedge$ is split into two $3$-cells:

\begin{center}
\includegraphics{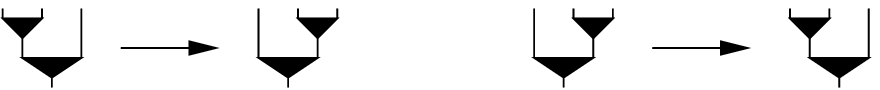}
\end{center}

\medskip
\noindent Then, in order to ensure that they are $3$-dimensional isomorphisms, one adds equations between proofs: both possible composites are equal to the corresponding identity. Hence, this leaves no equation between objects of dimension $2$, while two of them appear between objects of dimension $3$ for each equation on formulas. \\

\noindent\textbf{From equations to $\mathbf{4}$-dimensional computations.} There is no reason to stop the process of lifting up equations. In order to achieve this, the pairs of $3$-dimensional cells replacing equations are keeped, but equations between $3$-dimensional composites are lifted up. Hence, instead of considering commutative diagrams between $3$-dimensional arrows, one defines $4$-dimensional cells. Each one represents a \emph{computation} from one composite to the identity $3$-cell, such as in the following diagram:
\begin{center}
\begin{picture}(0,0)%
\includegraphics{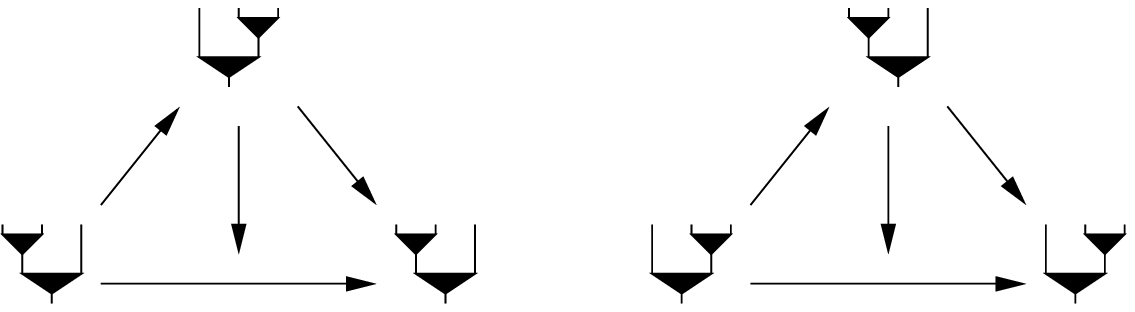}%
\end{picture}%
\setlength{\unitlength}{4144sp}%
\begingroup\makeatletter\ifx\SetFigFont\undefined%
\gdef\SetFigFont#1#2#3#4#5{%
  \reset@font\fontsize{#1}{#2pt}%
  \fontfamily{#3}\fontseries{#4}\fontshape{#5}%
  \selectfont}%
\fi\endgroup%
\begin{picture}(5154,1493)(169,1428)
\put(4231,1469){\makebox(0,0)[b]{\smash{{\SetFigFont{10}{12.0}{\rmdefault}{\mddefault}{\updefault}$=$}}}}
\put(1261,1469){\makebox(0,0)[b]{\smash{{\SetFigFont{10}{12.0}{\rmdefault}{\mddefault}{\updefault}$=$}}}}
\end{picture}%
\end{center}

\noindent There, the equation about the associativity of $\wedge$ is finally replaced by two $3$-cells, together with the above pair of $4$-cells. When this transformation is done, there is no more equations between $2$-arrows (formulas) or $3$-arrows (proofs). Only computations between proofs remain, in the form of pairs of $4$-cells.

\subsection{Local computations on proofs}\label{sub:calcul-4d}

The next example of $4$-dimensional cells is in fact a generalization of the former one. Indeed, it arises whenever one wants to compute normal forms for proofs, \emph{modulo} some specified equations.

This encompasses the former example, since these equations can be the ones stating that two $3$-cells are inverse one another. As an example of generalized computation, the following $4$-cell can be introduced in order to simplify proofs with a weakening followed by a contraction, both acting on the same atom:
\begin{center}
\includegraphics{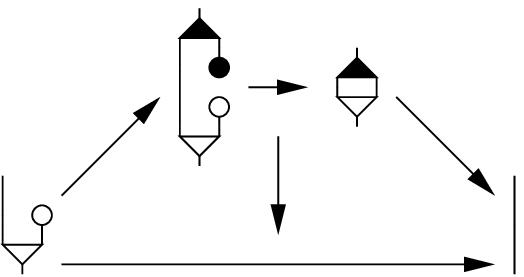}
\end{center}

\medskip
\noindent In fact, any local computation on proofs can be replaced by a $4$-cell. All the $4$-cells being given, the computations they generate are the $4$-dimensional arrows of a free $4$-category.

\begin{rem*}
Let us make an immediate remark on this $4$-cell. A weakening followed by a contraction is some kind of higher-dimensional version of the composition of a duplicator, followed by an eraser. Hence, this $4$-cell should be part of a family of $4$-dimensional resources management cells, an higher-dimensional version of the family $E_{\Delta}$ of $3$-dimensional ones.

We need to explore this potential family and, for example, check if it is automatically produced by its $3$-dimensional version. Another topic is to study its rewriting properties. A conjecture is that this family controls another form of bureaucracy, called type C in [Stra\ss{}burger 2005], which is not geometric like the other two.
\end{rem*}


\subsection{A word on cut-elimination}\label{sub:coupure}

We have not discussed cut-elimination, though it is the most known and studied computation on proofs. This is due to the fact that there is no known cut-elimination procedure on the system SKS which is generated by local rules between parallel proofs, unlike the ones known for various sequent calculi.

Indeed, the known procedure is a global algorithm, which takes into account the whole proof in order to eliminate the cuts [Brünnler 2004]. Hence, at least for the time being, there is no $4$-dimensional interpretation for the cut-elimination mechanism.

\subsection{Some temporary relief}\label{sub:commentaires}

The $4$-dimensional point of view immediately arises the following question: how can one use the fact that these computations are $4$-dimensional objects? This comes with the subsidiary question: how can one represent $4$-dimensional objects? In fact, this is not necessary at this point.

To explain this answer, let us step back by one dimension. Term rewriting is about some properties (termination and confluence) of computations on $2$-dimensional objects. While considering the whole $2$-dimensional structure of terms is really useful, the computations need not be seen as genuine $3$-dimensional objects: the only purpose of doing so would be to identify reduction paths \emph{modulo} bureaucracy. But term rewriting is not concerned with the classification of reduction paths (only their existence) and neither termination nor confluence are modified by bureaucracy.

Then comes proof theory which, with the higher-dimensional point of view, studies $3$-dimensional objects, or rather computations between them. Hence, with the same arguments as above, considering the whole $3$-dimensional structure of proofs shall prove to be useful. But the four dimensions of computations on proofs are not involved if one only wants to prove termination or confluence of proof normalization processes.

In conclusion, if it is only about (normalization of) proofs, then one can live with rewriting paths on $3$-dimensional arrows. But when times will come when the classification of rewriting paths on proofs is concerned, then the fourth dimension will be useful.

In order to manage the six types of geometric bureaucracy lurking in dimension $4$, for example\dots

\section{A polygraphic translation of SLLS}\label{sec:slls}

This calculus of structures-style formalism is presented in [Stra\ss{}burger 2003] and describes proofs of propositional linear logic [Girard 1987]. Since its structure is very similar to the one of SKS, we present here a polygraph which is (strongly) conjectured to satisfy the same properties with respect to SLLS as~$\Sigma^K$ does with respect to SKS.

In term-like version, the signature of SLLS has two sorts $A$ and $F$ and the following constructors:
$$
\xymatrix{&&& A \ar@(ul,ur)^{\nu} \ar[d]^{\iota} \\ \ast \ar[rrr]_{\top,\:\bot,\:1,\:0} &&& F \ar@(dr,dl)^{!,\:?} &&& F\times F \ar[lll]^{\tens,\:\oplus,\:\parll,\:\etll}}
$$

\noindent Terms are equipped with the structural congruence generated by the following rules, where $(\mu,\eta)$ is any pair among $(\parll,\bot)$, $(\tens,1)$, $(\oplus,0)$ and $(\etll,\top)$:
$$
\begin{array}{r c l c r c l}
&&\hfill\mu(\mu(x_1,x_2),x_3) &\lfl& \mu(x_1,\mu(x_2,x_3))\hfill\strut \\
&&\hfill\mu(\eta,x_1) &\lfl& x_1\hfill\strut \\
&&\hfill\mu(x_1,x_2) &\lfl& \mu(x_2,x_1)\hfill\strut\\
?(\bot) &\lfl &\bot && !(1) &\lfl& 1 \\
\bot\oplus\bot &\lfl &\bot && 1\etll 1 &\lfl& 1\\
&&\hfill\nu(\nu(a_1)) &\lfl& a_1.\hfill\strut
\end{array}
$$

\medskip
\noindent The same argumentation as the one developped for SKS throughout the section leads to the replacement of the set of formulas by the free $2$-category generated by the following $2$-polygraph with one cell in dimension $0$, two cells in dimension $1$ and twenty cells in dimension $2$:
\begin{center}
\begin{picture}(0,0)%
\includegraphics{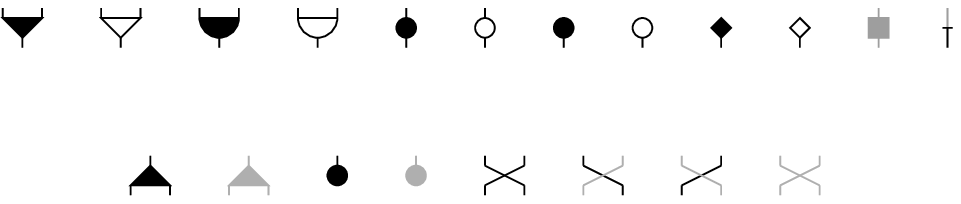}%
\end{picture}%
\setlength{\unitlength}{4144sp}%
\begingroup\makeatletter\ifx\SetFigFont\undefined%
\gdef\SetFigFont#1#2#3#4#5{%
  \reset@font\fontsize{#1}{#2pt}%
  \fontfamily{#3}\fontseries{#4}\fontshape{#5}%
  \selectfont}%
\fi\endgroup%
\begin{picture}(4367,1186)(349,-515)
\put(4681,209){\makebox(0,0)[b]{\smash{{\SetFigFont{10}{12.0}{\rmdefault}{\mddefault}{\updefault}{\color[rgb]{0,0,0}$\iota$}%
}}}}
\put(4366,209){\makebox(0,0)[b]{\smash{{\SetFigFont{10}{12.0}{\familydefault}{\mddefault}{\updefault}{\color[rgb]{0,0,0}$\nu$}%
}}}}
\put(4006,-466){\makebox(0,0)[b]{\smash{{\SetFigFont{10}{12.0}{\rmdefault}{\mddefault}{\updefault}{\color[rgb]{0,0,0}$\tau_{AA}$}%
}}}}
\put(3556,-466){\makebox(0,0)[b]{\smash{{\SetFigFont{10}{12.0}{\rmdefault}{\mddefault}{\updefault}{\color[rgb]{0,0,0}$\tau_{AF}$}%
}}}}
\put(3106,-466){\makebox(0,0)[b]{\smash{{\SetFigFont{10}{12.0}{\rmdefault}{\mddefault}{\updefault}{\color[rgb]{0,0,0}$\tau_{FA}$}%
}}}}
\put(2656,-466){\makebox(0,0)[b]{\smash{{\SetFigFont{10}{12.0}{\rmdefault}{\mddefault}{\updefault}{\color[rgb]{0,0,0}$\tau_{FF}$}%
}}}}
\put(1486,-466){\makebox(0,0)[b]{\smash{{\SetFigFont{10}{12.0}{\rmdefault}{\mddefault}{\updefault}{\color[rgb]{0,0,0}$\delta_A$}%
}}}}
\put(2251,-466){\makebox(0,0)[b]{\smash{{\SetFigFont{10}{12.0}{\rmdefault}{\mddefault}{\updefault}{\color[rgb]{0,0,0}$\epsilon_A$}%
}}}}
\put(1891,-466){\makebox(0,0)[b]{\smash{{\SetFigFont{10}{12.0}{\rmdefault}{\mddefault}{\updefault}{\color[rgb]{0,0,0}$\epsilon_F$}%
}}}}
\put(1036,-466){\makebox(0,0)[b]{\smash{{\SetFigFont{10}{12.0}{\rmdefault}{\mddefault}{\updefault}{\color[rgb]{0,0,0}$\delta_F$}%
}}}}
\put(451,209){\makebox(0,0)[b]{\smash{{\SetFigFont{10}{12.0}{\familydefault}{\mddefault}{\updefault}{\color[rgb]{0,0,0}$\parll$}%
}}}}
\put(901,209){\makebox(0,0)[b]{\smash{{\SetFigFont{10}{12.0}{\familydefault}{\mddefault}{\updefault}{\color[rgb]{0,0,0}$\tens$}%
}}}}
\put(1351,209){\makebox(0,0)[b]{\smash{{\SetFigFont{10}{12.0}{\familydefault}{\mddefault}{\updefault}{\color[rgb]{0,0,0}$\oplus$}%
}}}}
\put(1801,209){\makebox(0,0)[b]{\smash{{\SetFigFont{10}{12.0}{\familydefault}{\mddefault}{\updefault}{\color[rgb]{0,0,0}$\etll$}%
}}}}
\put(2206,209){\makebox(0,0)[b]{\smash{{\SetFigFont{10}{12.0}{\familydefault}{\mddefault}{\updefault}{\color[rgb]{0,0,0}$!$}%
}}}}
\put(2566,209){\makebox(0,0)[b]{\smash{{\SetFigFont{10}{12.0}{\familydefault}{\mddefault}{\updefault}{\color[rgb]{0,0,0}$?$}%
}}}}
\put(2926,209){\makebox(0,0)[b]{\smash{{\SetFigFont{10}{12.0}{\familydefault}{\mddefault}{\updefault}{\color[rgb]{0,0,0}$\bot$}%
}}}}
\put(3286,209){\makebox(0,0)[b]{\smash{{\SetFigFont{10}{12.0}{\familydefault}{\mddefault}{\updefault}{\color[rgb]{0,0,0}$1$}%
}}}}
\put(3646,209){\makebox(0,0)[b]{\smash{{\SetFigFont{10}{12.0}{\familydefault}{\mddefault}{\updefault}{\color[rgb]{0,0,0}$0$}%
}}}}
\put(4006,209){\makebox(0,0)[b]{\smash{{\SetFigFont{10}{12.0}{\familydefault}{\mddefault}{\updefault}{\color[rgb]{0,0,0}$\top$}%
}}}}
\end{picture}%
\end{center}

\noindent Then, we consider the family of 48 resources management equations, given through the following schemes:
\begin{center}
\begin{picture}(0,0)%
\includegraphics{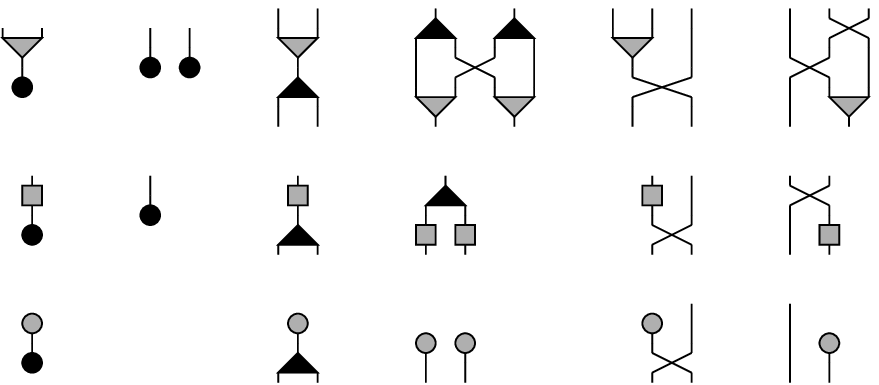}%
\end{picture}%
\setlength{\unitlength}{4144sp}%
\begingroup\makeatletter\ifx\SetFigFont\undefined%
\gdef\SetFigFont#1#2#3#4#5{%
  \reset@font\fontsize{#1}{#2pt}%
  \fontfamily{#3}\fontseries{#4}\fontshape{#5}%
  \selectfont}%
\fi\endgroup%
\begin{picture}(3984,1734)(-461,-973)
\put(2926,434){\makebox(0,0)[b]{\smash{{\SetFigFont{10}{12.0}{\rmdefault}{\mddefault}{\updefault}{\color[rgb]{0,0,0}$\equiv$}%
}}}}
\put(2926,-826){\makebox(0,0)[b]{\smash{{\SetFigFont{10}{12.0}{\rmdefault}{\mddefault}{\updefault}{\color[rgb]{0,0,0}$\equiv$}%
}}}}
\put(1216,-826){\makebox(0,0)[b]{\smash{{\SetFigFont{10}{12.0}{\rmdefault}{\mddefault}{\updefault}{\color[rgb]{0,0,0}$\equiv$}%
}}}}
\put(-44,-826){\makebox(0,0)[b]{\smash{{\SetFigFont{10}{12.0}{\rmdefault}{\mddefault}{\updefault}{\color[rgb]{0,0,0}$\equiv$}%
}}}}
\put(-44,-241){\makebox(0,0)[b]{\smash{{\SetFigFont{10}{12.0}{\rmdefault}{\mddefault}{\updefault}{\color[rgb]{0,0,0}$\equiv$}%
}}}}
\put(1216,-241){\makebox(0,0)[b]{\smash{{\SetFigFont{10}{12.0}{\rmdefault}{\mddefault}{\updefault}{\color[rgb]{0,0,0}$\equiv$}%
}}}}
\put(2926,-241){\makebox(0,0)[b]{\smash{{\SetFigFont{10}{12.0}{\rmdefault}{\mddefault}{\updefault}{\color[rgb]{0,0,0}$\equiv$}%
}}}}
\put(-44,434){\makebox(0,0)[b]{\smash{{\SetFigFont{10}{12.0}{\rmdefault}{\mddefault}{\updefault}{\color[rgb]{0,0,0}$\equiv$}%
}}}}
\put(1216,434){\makebox(0,0)[b]{\smash{{\SetFigFont{10}{12.0}{\rmdefault}{\mddefault}{\updefault}{\color[rgb]{0,0,0}$\equiv$}%
}}}}
\end{picture}%
\end{center}

\noindent Finally, the structural rules are translated into the following family of 17 rules on parallel $2$-arrows:
\begin{center}
\includegraphics{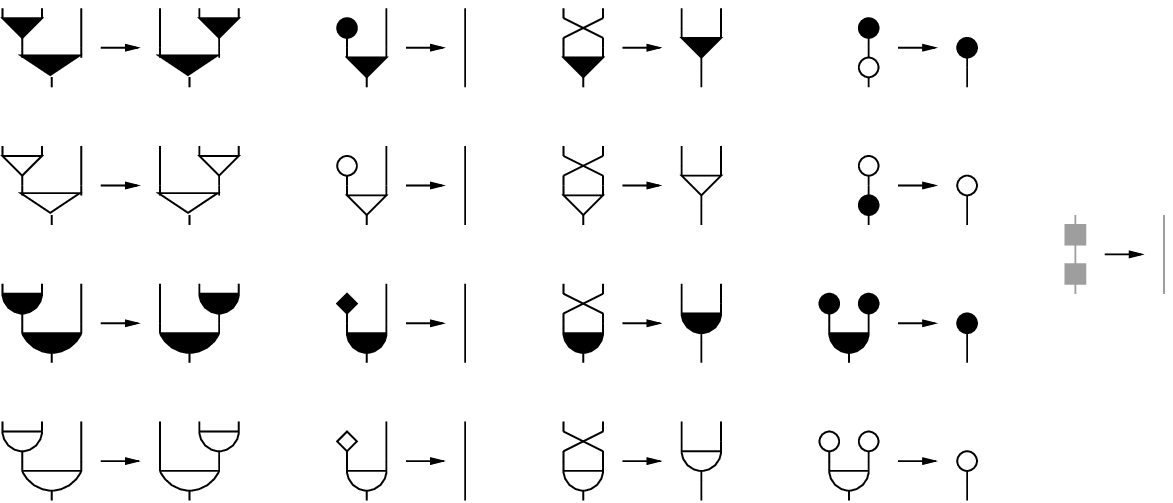}
\end{center}

\noindent In [Stra\ss{}burger 2003], the inference rules of system SLLS are given in a term-rewriting style. We do not recall them from there and instead directly give the corresponding $3$-cells, placed in the same order as in the original manuscript so that each one can be recognized:

\begin{center}
\includegraphics{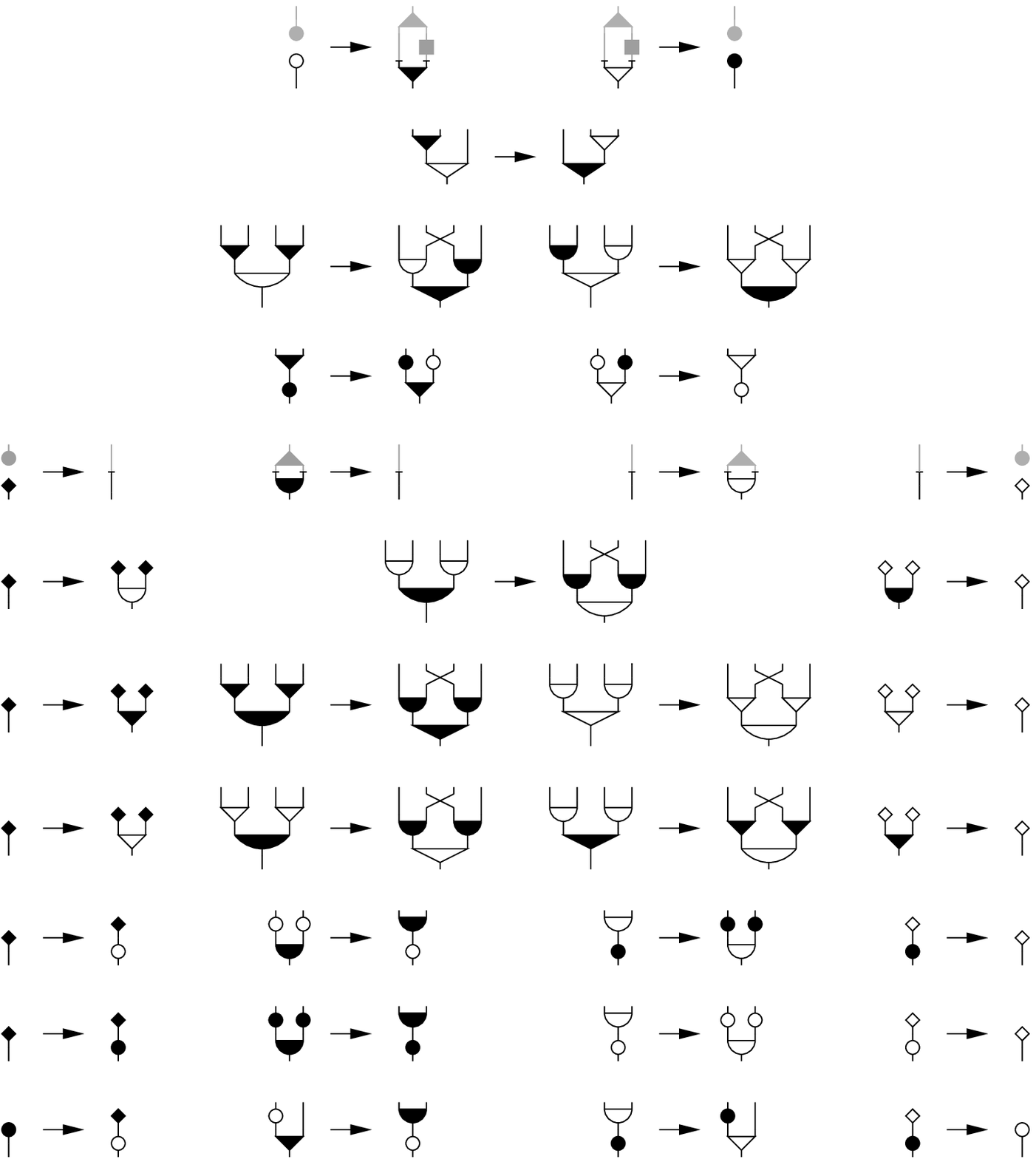}
\end{center}

\vfill\pagebreak
\section*{Comments and future directions}
\emptysectionmark{Comments and future directions}

This documents presents a higher-dimensional rewriting point of view for the deep inference system named SKS and, similarly, for the one called SLLS. One benefit of this setting is to provide a uniform theory for many possible systems, depending on what the user wants to emphasize. Indeed, much freedom is left on how to consider bureaucracy or how to see equations. Moreover, bureaucracy of geometrical nature is easily handled through the exchange relations. However, there is a bureaucracy type~C [Stra\ss{}burger 2005] which seems to come from a higher-dimensional version of the resources management rules. This type of bureaucracy must be studied to check if it can be described easily in the polygraphic language.

Higher-dimensional rewriting provides a common view on equations and computations between proofs: they are seen as $4$-dimensional cells between proofs. So one just has to choose the local computations he wants to study, then the $4$-dimensional rewriting theory can be used to see if the generated calculus is terminating or not, confluent or not. Yet, some work will be necessary here to provide the required tools, such as a recipe to craft termination orders like the one in [Guiraud 2004] for $3$-dimensional rewriting. Another tool will concern the study of $4$-dimensional critical pairs. This will be an adaptation of one that is still under development for $3$-dimensional critical pairs and will be described in a subsequent paper. In particular, these $3$-dimensional and $4$-dimensional tools will give answer on the existence of a finite and convergent family of $3$-cells which is equivalent to the union of the resources management relations and of the structural equations, for SKS and for SLLS. If there is no such convergent family, we should seek a finite equivalent family of $3$-cells, with a $4$-dimensional finite convergent rewriting calculus on it.

Aside from these computational issues, proofs seen as $3$-dimensional objects are naturally equipped with a graphical representation, using $3$-dimensional Penrose diagrams. The links between these pictures and proof nets still have to be explored. For the moment, we can at least say that the proposed $3$-dimensional representations provide a completely different way to look at proofs.

Another comment is that the $3$-dimensional translation of proofs relies on the calculus of structures version of the considered logic. One consequence is that we have to consider $2$-sorted terms and, thus, polygraphs with two generating $1$-cells. As we have seen, this always makes notations and constructions (much) tougher. Another negative point is the fact that cut-elimination cannot be described locally, which is disturbing for such an important relation between parallel proofs. A future work will propose a direct $3$-dimensional version of proofs, based on equivalences in the theory of boolean algebras, in which cut-elimination will be a $4$-dimensional computation.

The final comment concerns binders: for the time being, higher-dimensional rewriting is unable to handle them. This is a major issue which is to be solved, either by proposing a polygraphic account of the $\lambda$-calculus or by extending the higher-dimensional setting to encompass it. The main step to reach this goal is to check if there exists a polygraphic presentation of the structure of cartesian closed category, like the one that was found for cartesian categories and layed the bases of the field [Burroni 1993]. 

\vfill
\begin{flushright}
\begin{minipage}{130mm}
\emph{I wish to thank the referee for remarks that have greatly helped to improve the paper. I also wish to thank all the people from Marseille and (formerly) from Dresden who have (patiently) listened to these results and commented them.}
\end{minipage}
\end{flushright}

\section*{References}
\emptysectionmark{References}

\bib{Franz Baader, Tobias Nipkow}{Term rewriting and all that}{Cambridge University Press, 1998}

\medskip
\bib{John Carlos Baez, James Dolan}{Categorification}{ArXiv preprint, 1998}

\medskip
\bib{Kai Br\"unnler}{Deep inference and symmetry in classical proofs}{Logos Verlag, 2004}

\medskip
\bib{Albert Burroni}{Higher-dimensional word problems with applications to equational logic}{\\ \indent\indent Theoretical Computer Science 115(1), 1993}

\medskip
\bib{Eugenia Chang, Aaron Lauda}{Higher-dimensional categories: an illustrated guide book}{2004}

\medskip
\bib{Jean-Yves Girard}{Linear logic}{Theoretical Computer Science 50(1), 1987}

\medskip
\bibdeux{Alessio Guglielmi}{The problem of bureaucracy and identity of proofs from the perspective of deep inference}{\\ \indent\indent Proceedings of Structures and deduction ICALP workshop, 2005}{A system of interaction and structure}{\\ \indent\indent ACM Transactions on Computational Logic, to be published (2004)}

\medskip
\bibtrois{Yves Guiraud}{Présentations d'opérades et systèmes de réécriture}{thèse de doctorat, 2004(T)}{Termination orders for $3$-dimensional rewriting}{\\ \indent\indent Journal of Pure and Applied Algebra, to be published (2004)}{Two polygraphic presentations of Petri nets}{submitted preprint, 2005}

\medskip
\bib{Yves Lafont}{Towards an algebraic theory of boolean circuits}{Journal of Pure and Applied Algebra 184, 2003}

\medskip
\bib{Saunders MacLane}{Categories for the working mathematician}{Springer, 1998}

\medskip
\bib{François Métayer}{Resolutions by polygraphs}{Theory and Applications of Categories 11(7), 2003}

\medskip
\bibdeux{Lutz Stra\ss{}burger}{Linear logic and noncommutativity in the calculus of structures}{PhD thesis, 2003}{From deep inference to proof nets}{Structures and Deduction ICALP worshop, 2005}

\end{document}